\documentclass[12pt]{amsart}

\usepackage[T1]{fontenc}
\usepackage[utf8]{inputenc}
\usepackage[english]{babel}
\usepackage{epsfig}
\usepackage{amsmath}
\usepackage{amssymb}
\usepackage{amscd}
\usepackage{latexsym}
\usepackage{tabularx}
\usepackage{a4wide}
\usepackage[usenames]{color}
\usepackage{enumerate}
\usepackage{subfigure}
\usepackage{url}
\usepackage{verbatim}

\usepackage{circuitikz}
\ctikzset{bipoles/resistor/height=0.15}
\ctikzset{bipoles/resistor/width=0.4}

\usepackage[
colorlinks, citecolor=blue,
pdfauthor={Omid Amini, Eduardo Esteves},
pdftitle={Voronoi tilings, toric arrangements and degenerations of line bundles I},
pdfstartview ={FitV},
]{hyperref}

\usepackage{tikz, float}
\usetikzlibrary {positioning}

\usetikzlibrary{calc,decorations.markings}
\usetikzlibrary{shapes,snakes}

\numberwithin{equation}{section}

\tikzstyle{Cwhite}=[scale = .8,circle, fill = white, minimum size=3mm]
\tikzstyle{Cgray}=[scale = .4,circle, fill = gray, minimum size=3mm]
\tikzstyle{Cblack2}=[scale = .4,circle, fill = black, minimum size=5mm]
\tikzstyle{Cblack}=[scale = .7,circle, fill = black, minimum size=3mm]
\tikzstyle{C0}=[scale = .9,circle, fill = black!0, inner sep = 0pt, minimum size=3mm]
\tikzstyle{C1}=[scale = .7,circle, fill = black!0, inner sep = 0pt, minimum size=3mm]
\tikzstyle{Cred}=[scale = .4,circle, fill = red, minimum size=3mm]

\usepackage[matrix,arrow,tips,curve]{xy}
\usepackage{pb-diagram,pb-xy}
\usepackage{verbatim}
\usepackage{mathrsfs}
\usepackage{color}
\usepackage[leqno]{amsmath}

\newcommand{\Vor}{\mathrm{Vor}}
\newcommand{\FP}{\mathcal{FP}}

\newcommand{\E}{\mathbb E}
\newcommand{\im}{\textrm{im}}
\newcommand{\AC}{\mathcal{CAC}}
\newcommand{\supp}{\mathrm{supp}}
\newcommand{\Bo}{\Upsilon_1}

\newcommand{\m}{\mathfrak m}

\small\normalsize

\newcommand{\wt}{\widetilde}
\renewcommand{\:}{\colon}

\newcommand{\I}{\mathcal I}

\newcommand{\ol}{\overline}

\newtheorem{thm}{Theorem}[section]

\newtheorem{lemma}[thm]{Lemma}

\newtheorem{propdefi}[thm]{Proposition - Definition}
\newtheorem{prop}[thm]{Proposition}

\newtheorem{claim}[thm]{Claim}
\newtheorem{cor}[thm]{Corollary}

\newtheorem{remark}[thm]{Remark}

\newtheorem{defi}[thm]{Definition}

\newcommand{\dl}{\mathfrak d}

\def\X{\mathcal X}

\newcommand{\Z}{\mathbb{Z}}

\newcommand{\R}{\mathbb{R}}

\newcommand{\N}{\mathbb{N}}

\newcommand{\C}{\mathfrak{C}}

\newcommand{\ord}{\mathrm{ord}}

\newcommand{\Div}{\operatorname{Div}}

\newcommand{\Prin}{\operatorname{Prin}}
\renewcommand{\div}{\mathrm{div}}

\renewcommand{\k}{\kappa}

\newcommand{\f}{\mathfrak f}

\newcommand{\F}{F}

\newcommand{\Cay}{\mathrm{Cay}}

\numberwithin{equation}{section}

\newcommand{\te}{\text{h}}  %%%%target of e
\newcommand{\he}{\text{t}}  %%%%source of e

\usepackage{scalerel}

\newcommand{\ind}[1]{{_{\scaleto{#1}{4.8pt}}}}
\newcommand{\indm}[1]{{_{\scaleto{#1}{3.2pt}}}} %%% take care of indices like \chi_e \chi_v
\newcommand{\indmbar}[1]{{_{\scaleto{#1}{4.6pt}}}} %%% take care of indices like \chi_{\ol e}
\newcommand{\indM}[1]{{_{\scaleto{#1}{9.2pt}}}} %%% take care of indices like \chi_{\E()}
\newcommand{\indbi}[2]{_{{\scaleto{#1}{4.8pt}}_{\hspace{-.05cm}#2}}} %%%%take care of indices like \chi_{S_j}
\newcommand{\indmbi}[2]{_{{\scaleto{#1}{3.2pt}}_{\hspace{-.03cm}#2}}} %%%%take care of indices like \chi_{v_j}

\makeatletter
\newsavebox\myboxA
\newsavebox\myboxB
\newlength\mylenA

\newcommand*\overbar[2][0.75]{%
    \sbox{\myboxA}{$\m@th#2$}%
    \setbox\myboxB\null
    \ht\myboxB=\ht\myboxA%
    \dp\myboxB=\dp\myboxA%
    \wd\myboxB=#1\wd\myboxA%
    \sbox\myboxB{$\m@th\overline{\copy\myboxB}$}%  Overlined phantom
    \setlength\mylenA{\the\wd\myboxA}%   calc width diff
    \addtolength\mylenA{-\the\wd\myboxB}%
    \ifdim\wd\myboxB<\wd\myboxA%
       \rlap{\hskip 1\mylenA\usebox\myboxB}{\usebox\myboxA}%
    \else
        \hskip -0.5\mylenA\rlap{\usebox\myboxA}{\hskip 0.5\mylenA\usebox\myboxB}%
    \fi}
\makeatother

\renewcommand{\o}{\mathfrak o}
\renewcommand{\P}{\mathbb P}
\title{Voronoi tilings, toric arrangements and degenerations of line bundles I}
 \author{Omid Amini}
 \author{Eduardo Esteves}
 \date{December 31, 2020}
  \address{CNRS - Centre de math\'ematiques Laurent Schwartz, \'Ecole Polytechnique,  France}
\email{omid.amini@polytechnique.edu}

\address{Instituto Nacional de Matem\'atica Pura e Aplicada, Estrada Dona Castorina 110,
22460-320 Rio de Janeiro RJ, Brazil}
\email{esteves@impa.br}
\begin{document}

\maketitle
\begin{abstract}  We describe limits of line bundles on 
nodal curves in terms of toric arrangements associated to 
Voronoi tilings of Euclidean spaces.
These tilings encode information on the relationship between 
the possibly infinitely many limits, and ultimately give rise to a new
definition of \emph{limit linear series}. This paper and its second 
and third companion parts are the first in a series aimed to explore
this new approach. In the present article, we set up the combinatorial 
framework and show how graphs with integer lengths
associated to the edges provide tilings of Euclidean spaces by 
certain polytopes associated to the graph itself and to certain of its
subgraphs. We further provide a description of the combinatorial 
structure of these polytopes and the way they are glued together in the tiling.

In the second part of the series, we describe the 
\emph{arrangements of  toric varieties} associated to these tilings. 
These results will be of use in the third part to achieve our goal 
of describing all \emph{stable limits} of
a family of line bundles along a degenerating family of curves.
  \end{abstract}

\vskip0.5cm

\tableofcontents

\addtocontents{toc}{\protect\setcounter{tocdepth}{1}}

\section{Introduction}
The aim of this work and its sequels~\cite{AE2, AE3} is to describe stable limits of line bundles on nodal curves in
terms of combinatorial geometric properties of graphs. In this introduction, after briefly mentioning how graphs enter in the picture for such a purpose, we mainly
concentrate on the combinatorial results related to graphs and their geometry presented in this paper, which we hope should be of  independent interest.
A reader interested in algebraic geometry shall find a more detailed exposition
and a more clear picture of the link to algebraic geometry in the introductions to the second and third parts.

Algebraic curves are the most studied and the most understood objects in algebraic geometry.
They were initially studied as embedded objects, first in the projective plane, then in the projective three dimensional space, and
later in higher-dimensional spaces. However, nowadays, it has become clear that the right point of view for
many applications is to
view algebraic curves abstractly and study their various embeddings in
terms of \emph{linear series}, which are simply the data of a line bundle on the curve and a linear subspace of the
space of sections of that line bundle. The line bundle and the chosen subspace of its global sections
give rise to a morphism of the curve into a projective space.

 Abstractly, we can study how curves vary and that leads naturally to the construction of the moduli space of smooth curves of a given genus. That space is not compact and though different compactifications of it can be pursued, the one that has established
itself is that by Deligne and Mumford~\cite{DM69}, which is obtained by adding all the \emph{stable curves} to the boundary. It is thus
natural to study how line bundles and their linear series vary, not just along families of
smooth curves but along families that include also stable curves. This study has turned
out to be a complex and multifaceted problem.

Stable curves are algebraic curves (proper, one-dimensional, reduced,
connected but not necessarily irreducible schemes) that fail to be smooth in the weakest possible form: 
the singularities of a stable curve are
all normal crossings, that is, they are ordinary nodes. Apart from this condition, they
are characterized as being those which have an ample canonical bundle, and this is in
fact the property that allows for the construction of their moduli
space. One of the key properties of stable curves that allows for the
construction of a decent moduli space is the stable reduction theorem,
which says that a family of stable curves parameterized by
a punctured smooth curve can be completed, after a finite base change, in a unique way to
a family over the whole smooth curve. Unfortunately, no such thing
holds, in general, for line bundles and linear series.

The main combinatorial object associated to a nodal curve $X$ (a curve with only normal crossings singularities) is its \emph{dual
graph}, and it has been now understood that the combinatorics of the dual graphs plays a crucial role 
in understanding various questions on the (limiting) behavior of families of curves, see e.g.~\cite{Ami-W, Ami-hdr, AN20, BJ, CC, CKV, CKV2, EM, ES07, OS79}. 

The dual graph $G=(V,E)$ of a nodal curve $X$ consists of a vertex set $V$ in one-to-one
correspondence with the set of components of $X$ and an edge set $E$ in
one-to-one correspondence with the set of all the nodes. These
bijections are also compatible, meaning that two vertices are
connected by an edge if the node corresponding to the edge lies on
both the components corresponding to the vertices.

The geometric study of graphs is linked in many ways to polyhedral geometry. In particular, one can associate to a graph several types of polytopes, arising also quite naturally in applications~\cite{Sch}. This includes the famous Edmonds' \emph{matching polytope}, see~\cite[Section 7]{LP} or \cite[Section 25]{Sch}, a fundamental object in combinatorial optimization  with plenty of
diverse applications, the \emph{stable polytope} defined by Chv\'atal~\cite{Chv} and further studied by 
Nemhauser and Trotter~\cite{NT} and Padberg~\cite{pad1,pad2}, which is a natural generalization of the matching polytope, and the  \emph{flow polytope} associated to an oriented graph~\cite{Sch}, used by Altmann and Hille~\cite{AH99} in connection to some moduli problems and by~\cite{BZ, BD, BS, Chen} for the study of graph polynomials (from the point of view of Erhard Theory of integral points in polytopes).

There are two other natural polytopes we can associate to a given graph which are relevant to algebraic geometry. Recall first that a full rank lattice $L$ in a real vector space $V$ of dimension $n$ is a discrete subgroup of $V$ which  has rank  $n$. If the vector space $V$ comes with an Euclidean norm, the lattice $L$ gives rise to a tiling of $V$ by polytopes, which is called the \emph{Voronoi tiling}. Each Voronoi cell is a polytope centred at a point $\lambda$ of $L$, and by definition consists of all the points of the vector space $V$ which have $\lambda$ as their closest point in the lattice, see e.g. the pioneering work of Conway and Sloane for the discussion of the geometry of Voronoi cells of certain lattices~\cite{CS, CS2}.  

 A given graph $G$ gives naturally rise to two different lattices: that of \emph{integral flows}, the \emph{flow lattice}, and that of \emph{integral cuts}, the \emph{cut lattice}.  The study of these lattices was pioneered already in connection with algebraic geometry by the work of Bacher, de la Harpe and Nagnibeda~\cite{BHN97}, in which they studied several basic questions about their combinatorial properties and posed several interesting questions on the shape of their Voronoi polytopes. 
 
The \emph{flow Voronoi polytope} is the Voronoi cell associated to
 the flow lattice and appears naturally in the description of the
 local structure of the 
compactified Jacobians~\cite{CKV, CKV2}.  
 By Torelli Theorem for graphs~\cite{CV, Art,  Ger, DG}, it also 
allows the determination of the graphic matroid 
of the graph and gives a classification of the graph up to 
Whitney equivalences~\cite{Whi1, Whi2}.  The structure of the face 
poset of this polytope, conjectured by Caporaso and Viviani in~\cite{CV}, 
 was determined in a previous work by the first named
 author~\cite{Ami}. 
Furthermore, it was shown in~\cite{AM} that Riemann-Roch Theorem 
for graphs proved by Baker and Norine~\cite{BN06} 
  can be captured and reproved by studying the geometry of the Voronoi 
cell of the cut lattice under a certain simplicial distance function 
(this means a nonsymmetric distance function having
 a simplex as its ball of unity).

The main objective of this paper is to study the \emph{cut Voronoi polytope}, which is the Voronoi polytope associated to the cut lattice.  We shall describe the face structure of this polytope in terms of 
certain acyclic orientations of subgraphs of $G$. Moreover, when $G$
comes equipped with an integral edge length function and the extra
data of a \emph{twisting}, a concept which will be introduced later in this paper, we show that the Voronoi polytopes associated to $G$ and certain of its subgraphs, determined explicitly  by the arithmetic of the
length function and the twisting factor, give rise to a periodic tiling of the Euclidean space. 
These particularly esthetic  tilings give rise in~\cite{AE2} to what
we 
call \emph{toric arrangements}, naturally embedded into products of infinite chains of projective lines, and are used as the main tools  in~\cite{AE3} in solving the degeneration problem for line bundles.

\subsection{Overview of the results} 

Let $G=(V,E)$ be a given graph that we assume to be connected. 
For our purposes, we can discard those edges that form loops, so 
we assume $G$ to be loopless. Let $\mathbb E$ be the set of all 
oriented edges (arrows) that can be obtained from $E$: for each
edge we get two arrows pointing each to one of the two different vertices
incident to that edge. For an oriented edge $e\in\mathbb E$ with 
extremities $u$ and $v$, we write $e=uv$ if $e$ is oriented from 
$u$ to $v$, even if $e$ is not unique (in case $G$ has parallel edges). In this case we call
$u$ and $v$ the tail and head of $e$, respectively, and sometimes
denote them by $\he_e$ and $\te_e$.  Also, we let $\ol e$ denote 
the reverse arrow, i.e., $\ol e =vu$. In other words, 
$\te_{\ol e} = \he_{e}$ and $\he_{\ol e} = \te_{e}$.

Given a ring $A$, we associate to $G$ two complexes:
$$
d_A\colon C^0(G,A) \longrightarrow C^1(G,A)\quad\text{and}\quad
\partial_A\colon C_1(G,A) \longrightarrow C_0(G,A).
$$
Here, $C^0(G,A)$ is the $A$-module of functions $V\to A$, and 
$C_0(G,A)$ is the free $A$-module generated by $V$. Also, $C^1(G,A)$
is the $A$-module of functions $f\colon\mathbb E\to A$ subject to
the condition that $f(\ol e)=-f(e)$ for each $e\in\mathbb E$,
and $C_1(G,A)$ is the quotient of the free $A$-module generated by the elements of
$\mathbb E$ modulo the submodule generated by elements of the form $e+\ol e$ for all 
$e\in\mathbb E$. Finally, $d_A(f)(e)=f(v)-f(u)$ and $\partial_A(e)=v-u$ for
each oriented edge $e=uv\in\mathbb E$. 

There are natural isomorphisms: $C_0(G,A)\to C^0(G,A)$, taking $v$ to
the characteristic function $\chi_{\indm v}$; and $C_1(G,A)\to C^1(G,A)$
taking $e$ to $\chi_{\indm e}-\chi_{\indmbar{\ol e}}$. 
(Here $\chi_{\indm e}$ is the characteristic function of $e$ which takes 
value 1 on $e$ and value 0 on every other oriented edge,
including $\ol e$.)  Also, there are
bilinear forms $\langle\,,\rangle$ on $C_0(G,A)$ and $C_1(G,A)$ such that 
$\langle v\,,w\rangle=\delta_{v,w}$ for $v,w\in V$ and
$\langle e\,,f\rangle=\delta_{e,f}-\delta_{e,\overline f}$ for $e,f\in\mathbb E$. The
isomorphisms induce bilinear forms on $C^0(G,A)$ and 
$C^1(G,A)$ as well. Let $d^*_A\colon C^1(G,A)\to C^0(G,A)$ be
the homomorphism corresponding to $\partial_A$ under the
isomorphisms. It is easy to see that $d^*_A$ is the adjoint to $d_A$, that is,
$\langle f,d^*_A(h)\rangle=\langle d_A(f)\,,h\rangle$ 
for all $f\in C^0(G,A)$ and $h\in C^1(G,A)$.

Let $H_{0,A}:=\Bigl\{f\in C^0(G,A)\,|\,\sum f(v)=0\Bigr\}$ and $F_A:=\text{Im}(d_A)$. 
Denote by $\Delta_A:=d^*_Ad_A$ the Laplacian of the graph $G$. The homomorphism 
$d^*_A$ induces an injective map from $F_A$ to $H_{0,A}$. 
For $A=\mathbb R$ the map is an isomorphism. Moreover,  
the bilinear form $\langle\,,\rangle$ on $C^1(G, \mathbb R)$ induces
one on $F_{\mathbb R}$, which corresponds via $d^*_{\mathbb R}$ to 
the bilinear form $(\,,)$ on $H_{0,\R}$ defined by 
$$
\forall f,g\in H_{0, \mathbb R}\, \qquad 
(f\,,g) := \langle h\,, g\rangle,
$$
where $h$ is any element of $C^0(G, \R)$ with $\Delta_\R(h)= f$. 
Denote by $q$ the positive quadratic form corresponding  to $(\,,)$, 
so using the previous notation we have $q(f) := \langle h\,, f\rangle$.

In the present article we will describe an infinite family of tilings of
$H_{0,\mathbb R}$ by polytopes where each tiling in the family is associated  to the choice of an edge length function $\ell: E \rightarrow \mathbb Z_{>0}$ and a twisting factor $\m \in C^1(G,\Z)$.
Here, by a tiling of a real vector space $V$ by polytopes we mean an
(infinite) family $\mathcal F$ consisting of polytopes embedded in $V$
and covering $V$ which verify the following two conditions:
\begin{itemize}
\item Each face of  a polytope in $\mathcal F$ belongs to $\mathcal F$.
\item For each pair of polytopes $P$ and $Q$ in $\mathcal F$, the intersection $P\cap Q$ is either empty or is a common face of both polytopes.
\end{itemize}
 By removing from a polytope in $\mathcal F$ all
the faces of positive codimension, we get the
corresponding open face. The open faces of the 
polytopes in $\mathcal F$ then form a stratification of 
$H_{0,\mathbb R}$ into strata of various dimensions. 

In the case of the \emph{regular edge length function}, that is, 
when all the edges have length equal to one, and $\m$ is trivial, the tiling is 
the usual Voronoi tiling defined as follows: 
Let $\Lambda_A:=\text{Im}(d^*_A)$. Then, as we observed before, we have 
$\Lambda_{\mathbb R}=H_{0,\mathbb R}$, and $\Lambda_{\mathbb Z}$ forms a
sublattice of $H_{0,\mathbb Z}$.  Since the graph $G$ is connected, it
is a consequence of the celebrated \emph{Kirchhoff matrix-tree
  Theorem} that $\Lambda_{\mathbb Z}$ has finite index in
$H_{0,\mathbb Z}$ equal to the number of spanning trees of $G$. (A
spanning subgraph of a graph $G=(V,E)$ is a subgraph $T=(V, F)$ 
with the same vertex set $V$ and any edge set $F \subset E$; it is a
tree if it is a connected graph without cycles.) 
In this case, the tiling of $H_{0,\R}$ is what we call \emph{the
  standard Voronoi tiling} 
of $G$, which is by definition the Voronoi decomposition of
$H_{0,\mathbb R}$ with respect to the lattice $\Lambda_{\mathbb Z}$ 
and the quadratic form $q$: the
tiles are the full dimensional polytopes
$$
\text{Vor}_q(\beta):=\Bigl\{\eta\in H_{0,\mathbb R}\,\Bigl|\,
q(\eta-\beta)\leq q(\eta-\alpha)\text{ for every }\alpha\in 
\Lambda_{\mathbb Z}-\{\beta\}\Bigr\}
$$
for $\beta\in \Lambda_{\mathbb Z}$, called Voronoi cells, and their faces.

In the standard Voronoi tiling each Voronoi cell is a translation of
that centered at the origin, which we call the \emph{cut Voronoi
  polytope} of $G$. It comes with a partial order on its faces. 
One of our results in this paper establishes for any graph $G$ an isomorphism between this
partially ordered set and the partially ordered set of what we call \emph{coherent
acyclic orientations of cut subgraphs} of $G$. 

A \emph{cut subgraph} of $G$ is a spanning subgraph $G'$ of $G$ for which the
vertex set $V$ can be partitioned into non-empty subsets $V_i$, for $i$ in a finite set $I$, in such
a way that the edges of $G'$ are those of $G$ connecting vertices of
$V_i$ to $V_j$ for $i\neq j$. A coherent acyclic orientation of the cut subgraph $G'$ is an
orientation of $G'$ which is induced by the choice of partition of vertex set $V = \sqcup_{i\in I} V_i$ as above and a total order $<$ on $I$ in such a way that  for each pair of elements $i,j\in I$ with $i<j$ in the total order, all edges between $V_i$ and $V_j$ get the orientation with their head in $V_j$.

Finally, we define a partial order on coherent acyclic orientations of cut subgraphs of $G$ by saying that $D_1$ is \emph{smaller} than 
$D_2$, for coherent acyclic orientations $D_1$ and $D_2$ of cut
subgraphs $G_1$ and $G_2$ of $G$, if $\mathbb E(D_2)\subseteq\mathbb
E(D_1)$. Here, $\mathbb E(D_j)$ denotes the set of oriented edges of
$D_j$ in $G$ for $j=1,2$. 

We may now state one of the theorems we prove in this paper, see Theorem~\ref{thm:vor}:
\medskip

\noindent \textbf{Theorem A.} \emph{The face poset of the cut Voronoi
  polytope of $G$ is isomorphic to the poset of coherent
acyclic orientations of cut subgraphs of $G$}.

\medskip

The standard Voronoi tiling is one of the tilings we
consider. In the case of a more general length function $\ell$, 
the corresponding tiling might have
non-isomorphic top-dimensional cells, and the tiling in that case 
is made of a gluing of various cut Voronoi polytopes 
associated to certain connected
subgraphs of $G$, which may not include $G$ itself if the twisting
$\m$ is nontrivial! These subgraphs
are determined by basic arithmetic properties concerning divisibility
by the edge lengths. In order to motivate the consideration of these
tilings, we need to come back to our original algebro-geometric motivation.

Let thus $X$  be a connected nodal curve defined over an algebraically closed
field $\k$ such that $G$ is the graph obtained from the dual graph of
$X$ by removing the loops. Let $\pi\:\mathcal X\to B$ be a (one-parameter)
smoothing of $X$. Here, $B$ is the spectrum of $\k[[t]]$ and 
$\pi$ is a projective flat morphism whose generic fiber is
smooth and whose special fiber is isomorphic to $X$. We fix such an
isomorphism.  
The total space $\mathcal X$ is regular except
possibly at the nodes of
$X$. For each $e\in E$, the completion of the local ring of $\mathcal X$ at
the corresponding node $N_e$ is isomorphic to
$\k[[u,v,t]]/(uv-t^{\ell_e})$ for a certain integer $\ell_e>0$, called the
\emph{singularity degree} or the \emph{thickness} of $\pi$ at
$N_e$. If all $\ell_e=1$, then for each $v\in V$ the corresponding 
irreducible components $X_v$ of $X$ is a Cartier divisor of $\mathcal
X$; we call $\pi$ \emph{Cartier} in this case. 

A finite base change $B\to B$ is obtained by
sending $t$ to $t^n$ for a given $n$. The resulting family
$\pi^n\:\mathcal X^n\to B$ is similar to
the original one: the special fiber is the same, the generic fiber is
a base field extension of the original one, but the
singularity degrees $\ell_e$ change to $n\ell_e$.

Let now $L_\eta$ be an invertible sheaf on the generic fiber of
$\pi^m$ for some $m$. If $\mathcal X^m$ is regular (only if $m=1$),
then $L_\eta$ extends to an
invertible sheaf $\mathcal L$ on $\X^m$. More generally, $L_\eta$ extends
to a relatively torsion-free, rank-one sheaf $\mathcal L$ on $\X^m/B$,
that is, a $B$-flat coherent sheaf on $\mathcal X^m$ whose fibers over
$B$ are torsion-free, rank-one. The extensions are not unique. If
$\pi^m$ is Cartier, then for each $f\in C^0(G,\mathbb
Z)$, the sheaf $\mathcal L\otimes\mathcal O_{\mathcal X^m}(\sum
f(v)X_v)$ is another extension, and these are all the extensions. More
generally, in~\cite{Esteves01}, a similar procedure is described to change from one extension
to another. Or one could apply a sequence of blowups at nodes to get a
partial resolution $\widetilde\X^m\to\X^m$ such that the composition 
$\widetilde\pi^m\:\widetilde\X^m\to B$ is Cartier. The special fiber
of $\widetilde\pi^m$ is the curve $X^{m\ell}$ obtained from $X$ by
splitting each node $N_e$ apart and adding a chain $C^{(m\ell_e)}_e$ of $m\ell_e-1$
smooth rational curves connecting the branches. The
generic fiber of $\widetilde\pi^m$ is the same as that of $\pi^m$, so
$L_\eta$ extends to relatively
torsion-free, rank-one sheaves $\mathcal L$ on $\X^m/B$ that differ
from each other as
described above. Pushforwarding under the resolution map
$\widetilde\X^m\to\X^m$ gives us a bijection from the set of those
extensions $\mathcal L$ which are
\emph{admissible} (see below) to that of all the relatively torsion-free, rank-one
extensions of $L_\eta$  on $\X^m/B$.

Furthermore, one could do yet another finite base change, this time to
$\pi^m$, extend $L_\eta$ to the new generic fiber and consider its
torsion-free, rank-1 extensions. Of course, they will be extensions on a different total
space. But the special fibers of all the $\pi^n$ are the same. The
restrictions to $X$ of all these extensions for all integers $n$
divisbile by $m$ are torsion-free, rank-one sheaves that
we call the \emph{stable limits} of $L_\eta$.

Let $\mathcal L$ be a  relatively torsion-free, rank-one sheaf on
$\widetilde\X^m/B$ extending $L_\eta$. Let $H:=H^{m\ell}$ be the subgraph
of the dual graph of $X^{m\ell}$ obtained by removing the loops; it is
obtained directly from $G$ by 
subdividing  $m\ell_e-1$ times each edge $e$ of $G$. Since
$\widetilde\pi^m$ is Cartier, $\mathcal L$ restricts to torsion-free,
rank-one sheaves on the components of $X^{m\ell}$; their degrees form a
divisor $\widetilde D$ of $H$. We say that $\mathcal L$ and
$\widetilde D$ are admissible if for each edge $e$ of $G$ the value of
$\widetilde D$ is 0 at all but at most one vertex of $V(H)-V(G)$ on
$e$, where the value is 1. In Section~\ref{sec:admissible} we explain
how the various divisors $\wt D$ associated to admissible extensions
$\mathcal L$ on $\widetilde\X^m/B $ are related.
We show that they differ by the (principal)
divisors of certain functions
in $C^0(H,\Z)$ uniquely determined by their restrictions to
$G$; see Theorem~\ref{thm:admissible2}.

Furthermore, fix one extension $\mathcal L$ on $\widetilde\X^m/B$,
with associated divisor $\widetilde D$ on $H$. Fix an
orientation $\mathfrak o: E \to \E$ of the edges of
$G$, and denote by $E^{\mathfrak o} := \mathrm{Im}(\mathfrak o)$
the set of oriented edges in the orientation. Let $\mathfrak m\in
C^1(G,\Z)$ whose value at each $e\in E^{\mathfrak o}$ 
is equal to the number of edges between $\te_e$
and the vertex of $H$ on the subdvision of $e$ for which
$\widetilde D$ takes value 1. Also, for each
$f\in C^0(G,\Z)$ and each $e\in\E(G)$ put
$$
\delta_e^{\m}(f):=\Big\lfloor \frac{f(\te_e)-f(\he_e)+\mathfrak
  m_e}{m\ell_e}\Big\rfloor
\quad\text{and}\quad \dl^{\m}_f(e):=
\frac 12\Big( \delta_e^{\m}(f)+\delta_{\ol e}^{\m}(f)\Big).
$$
Then the $\dl^{\m}_f(e)$ are the values of a function
$\dl^{\m}_f\in C^1(G,\R)$. They are half integers; let
$G^{\m}_f$ be the spanning subgraph of $G$
whose edges are those supporting
$e\in\E(G)$ with $\dl^{\m}_f(e)\in\Z$. Let
$$
\Vor^{\m}_H(f):=d^*(\dl^{\m}_f)+\Vor_{G^{\m}_f}(O).
$$
In the present paper we prove that the set of polytopes
$\Vor^{\m}_H(f)$ with $G^{\m}_f$ connected, as $f$ ranges in
$C^0(G,\R)$, provide a tiling of $H_{0,\R}$, denoted $\Vor^{\m}_H$,
the most general tiling we consider; see
Theorem~\ref{thm:projection3}.

The tiling $\Vor^{\m}_H$ is obtained from purely
combinatorial data: $G$, $m\ell$ and
$\m$ chosen at will. But if it arises from algebro-geometric
degeneration data as above, we will see in \cite{AE3} that there is a
natural bijection between the set of all the tiles of
$\Vor^{\m}_H$, top-dimensional or not, and the set of all
stable limits of $L_\eta$. Furthermore, the toric arrangement
$Y^{\mathrm{st}}_{\ell,\m}$ associated to $\Vor^{\m}_{G,\ell}$, 
as we describe in \cite{AE2}, Section~4,
parameterizes a relatively torsion-free, rank-one sheaf on 
$X\times Y^{\mathrm{st}}_{\ell,\m}$ over $Y^{\mathrm{st}}_{\ell,\m}$ 
whose fibers are all the stable limits of $L_\eta$, giving a structure
to this set; see \cite{AE3}, Theorem~5.7.

\subsection{Content of the paper} After a brief presentation of the
divisor theory on graphs, in Section~\ref{sec:admissible} we define
admissible divisors and give a classification of those equivalent to a
given divisor in terms of principal $G$-admissible divisors. They are
then used to define the collection of all pairs $(D_I, G_I)$
associated to the stable limits in the above picture, although here we
only use the combinatorial definition and leave the link to the actual
stable limits to the third part \cite{AE3} of the work .  In Section~\ref{sec:tiling1} we treat the case of uniform edge lengths equal to one, and establish Theorem A. Sections~\ref{sec:generaltiling} and~\ref{sec:moregeneraltiling} are then devoted to the generalization of the picture given in Section 4 to the case of arbitrary edge lengths and the presence of a twisting. 

\subsection{Convention} For ease of reading, we gather here some of the terminology and notation we use in this paper, which are mostly standard and can be found
in any text book on graph theory, see e.g.~\cite{Bol, BM}. All our graphs will
be finite. Our graphs can have multiple edges but we will not allow
loops. The vertex and edge set of a graph
$G$ are denoted by $V(G)$ and $E(G)$, respectively (or just $V$ and
$E$, if the graph $G$ is clear from the context). A subgraph $G'$ of
$G$ is called spanning if $V(G')=V(G)$; the data of the spanning
subgraph $G'$ is then equivalent to the data of a subset $E'$ of
$E$. For a subset $X$ of $V(G)$, the induced subgraph on $X$ is the
subgraph whose vertex set is $X$ and whose edge set consists 
of all the edges $e$ of $G$ whose both endpoints are in $X$. 

For a graph $G=(V,E)$, we denote by $\E$ the set consisting of the
edges of $E$ with a choice of an orientation. Thus for each edge
$\{u,v\}$ in $E$, there is one oriented edge $uv$ in $\E$, having $u$
as its tail and $v$ as its head, as well as an oriented edge $vu$ in
$\E$, having $v$ as its tail and $u$ as its head in $\E$. If $e=uv$ is
an oriented edge in $\E$, the same edge with the reverse orientation
$vu$ is denoted by $\ol e$. Also, we sometimes
denote the tail of $e$ by $\he_e$ and its head by $\te_e$. 

An orientation $\o$ of the graph $G$ is a map from $E$ to $\E$ which to any edge $\{u,v\}$ associates one of the two oriented edges $uv$ or $vu$ in $\E$. 
An oriented graph is by definition a graph $G=(V,E)$ and a choice of
an orientation for the edges of $G$; the oriented graph associated to
the orientation is sometimes denoted by $G_{\mathfrak o}$. If $D$ is
an oriented graph, we also use $\E(D)$ to denote the oriented edges of
$D$. An oriented cycle  in an oriented graph $D$ is a sequence $v_1,
e_1, v_2, e_2, \dots, v_l, e_l, v_1$ such that $v_1, \dots, v_l$ are
distinct vertices and $e_1, \dots, e_l$ are oriented edges in $D$
satisfying $e_i = v_{i}v_{i+1}$ for each $i=1,\dots, l-1$ and 
$e_l = v_{l}v_1$. An oriented graph  $D$ is called acyclic if there is no oriented cycle in $D$; if $D$ is associated to an orientation $\o$ of $G$, the orientation $\o$ is called acyclic in that case.  

\section{Admissible Divisors on graphs}\label{sec:admissible}

\subsection{Divisors on graphs} In this subsection, we briefly recall the formalism of divisor theory on graphs which parallels the theory on algebraic curves. More on the subject can be found in~\cite{BHN97, BN06, Ami-hdr, BJ}. We then present in the next section a variant of the setup for graphs with the presence of edge-lengths. 

Let $G=(V,E)$ be a finite loopless and connected graph. Multiple edges are allowed. The group of divisors on $G$ denoted by $\Div(G)$ is by definition the free abelian group generated by vertices in $V$. We write $(v)$ for the generator associated to $v\in V$, so
 \[\Div(G) := \Bigl\{\sum_{v\in V} n_v(v)\,\, |\,\, n_v\in \mathbb Z\, \Bigr\}.\]
 
For a divisor $D\in \Div(G)$, the coefficient of $(v)$ in $D$ is
denoted $D(v)$. Its support, $\mathrm{Supp}(D)$, is the set of vertices $v$
with $D(v)\neq 0$. And its degree, $\mathrm{deg}(D)$, is defined as 
 \[\forall\,\, D\in \Div(G), \qquad \mathrm{deg}(D):= \sum_{v\in V} D(v).\]
 
Denote as before $C^0(G, \Z):=\{f: V \rightarrow \mathbb Z\}$, the set of all integer valued functions on vertices of  $G$. The functions in $C^0(G,\Z)$ play a role analogous to the role of  rational functions in the theory of algebraic curves. 

We can define the orders of vanishing of these rational functions at vertices as follows.  For each vertex $v\in V$, denote by $\ord_v : C^0(G,\Z) \rightarrow \mathbb Z$ the function of order of vanishing at $v$ which on $f\in C^0(G, \Z)$ takes the value 
 \[\ord_v(f) := \sum_{\substack{e\in\mathbb E\,\\ \te_e=v}} f(\he_e) -f(\te_e).\]
 
To any $f$ we associate the divisor $\div(f)$ defined by 
 \[\div(f) := \sum_{v\in V} \ord_v(f) (v).\]
  Elements of this form in $\Div(G)$ are called principal, and  the subgroup of $\Div(G)$ formed by principal divisors is denoted by $\Prin(G)$.

A divisor $D_1$ is called linearly equivalent to a divisor $D_2$ and we write $D_1 \sim D_2$ if the difference $D_1 -D_2$ is principal, i.e., there exists $f\in C^0(G,\Z)$ such that  $D_1 =D_2+ \div(f)$.

 \subsection{Chip firing game} There is a close connection between the theory of divisors on graphs and a game called chip-firing played on the graph~\cite{BTW, Dhar, Gab, Gab2, Biggs, BLS}.
 
Consider the following game played on a connected graph $G$. Vertices represent people in a group, where edges represent friendship. Each person $v\in V$ has a certain number of chips $n_v\in \mathbb Z$, where $n_v <0$ means $v$ is in debt. 
 The aim of the group is to achieve a situation where no one is in
 debt. The rule of the game is that at each step one person can decide
 to give one chip to each and all of its friends. The question is the existence of a winning strategy for the group. One can represent each configuration of the game with the corresponding divisor $\sum_{v\in V} n_v (v)$. A configuration can be reached from another if the two corresponding divisors are linearly equivalent. The existence of a winning strategy is then  equivalent to the existence of a function $f$ such that adding $\div(f)$ to the divisor associated to the configuration gives an effective divisor, i.e., a divisor with only non-negative coefficients.

\subsection{Admissible divisors on graphs} 
Let $G=(V, E)$ be a finite connected graph with integer edge lengths 
$\ell: E \to \mathbb N$. The length of an edge $e$ is denoted by 
$\ell_e$. We  extend $\ell$ to $\mathbb E$ via the forgetful map $\mathbb E\to E$.

Let $H$ be the graph obtained by subdividing $\ell_e-1$ 
times each edge $e$ of $G$, that is, 
by replacing each oriented edge $e=uv$ in $G$ by  a path 
$P_e = ux_1^ex_2^e\dots x_{\ell_e-1}^ev$ for new vertices
$x_1^e,\dots, x^e_{\ell_e-1}$. 
It will be convenient to define $x_0^e:=u$ and 
$x_{\ell_e}^e:=v$. Note that, with our convention, 
$x_i^{e} = x_{\ell_e-i}^{\ol e}$ for $i=0,\dots,\ell_e$.

In this section we define
$G$-admissible divisors on $H$ and characterize those that are
linearly equivalent to a given $D\in\Div(H)$. 

\begin{defi}[$G$-admissible divisors]\label{Gadmdiv}\rm
 A divisor $D$ on $H$  is called \emph{$G$-admissible} if for each
 oriented edge $e$ of $G$ the value of $D$ is $0$ at all but at 
most one vertex among the $x^e_j$ for 
$j=1,\dots,\ell_e-1$, where the value is $1$.
\end{defi}

\subsection{Characterization of $G$-admissible divisors: basic case} Our aim in this subsection and the following one is to characterize all admissible divisors which are linearly equivalent to a given divisor $D \in \mathrm{Div}(H)$. To simplify the presentation, we first treat the case where $D$ has support in the vertices of $G$, that is, 
$\mathrm{Supp}(D) \subseteq V(G) \subseteq V(H)$.  

Let $f\colon V(G) \rightarrow \mathbb Z$ be an integer valued function on the vertices of $G$. 
For each oriented edge $e=uv$ of $G$, define 
$$
\delta_{e}(f) := \lfloor \frac{f(v)-f(u)}{\ell_{e}}\rfloor.
$$
Note that the value of $\delta_{e}(f) + \delta_{\ol e}(f)$ is either 0 or $-1$, 
 depending on whether $\frac{f(v)-f(u)}{\ell_{e}}$ is an integer or not, respectively.

 \begin{defi}[Principal $G$-admissible divisor]\rm
The \emph{principal $G$-admissible} divisor $\div_\ell(f)$ associated to a 
function $f\colon V(G) \rightarrow \mathbb Z$ 
is the $G$-admissible divisor in $\Div(H)$ defined as follows:
\begin{itemize}
 \item For each vertex $u$ of $G$, define the coefficient of $\div_\ell(f)$ at $u$ by:
 $$\div_\ell(f)(u) := \sum_{\substack{e\in\mathbb E\,\\ \he_e=u}}\delta_e(f).$$ 
 \item For each oriented edge $e=uv$ of $G$ with difference
   $f(v)-f(u)$ not divisible by $\ell_e$, define
$$\div_\ell(f)(x_{r_e}^{\ol e}) :=1,$$
where $r_e$ is the remainder of the Euclidean division of $f(v)-f(u)$ by $\ell_e$.
 \item For each other vertex $z$ of $H$, define $\div_\ell(f)(z):=0$.
\end{itemize}
  \end{defi}
  
By definition, it is clear that $\div_\ell(f)$ is $G$-admissible.
We give now an alternative description of $\div_\ell(f)$, which shows
that it is indeed a 
principal divisor on $H$. For this, consider the
extension $\widetilde f: V(H) \rightarrow \mathbb Z$ of 
$f$ defined as follows: For each oriented edge $e=uv$ of $G$
and each $j=1,\dots,\ell_e-1$, set
\[
\widetilde f(x_j^e):=
\begin{cases}
    f(u) + j \lfloor \frac{f(v)-f(u)}{\ell_{e}}\rfloor &
    \qquad\textrm{if $j\leq\ell_e-r_e$},\\ 
    f(u) + j \lfloor \frac{f(v)-f(u)}{\ell_{e}}\rfloor + (j-\ell_e+r_e)
    & 
    \qquad \textrm{if $j\geq\ell_e-r_e$},
   \end{cases}  
\]
where $r_e$ is the remainder of the division of 
$f(v)-f(u)$ by $\ell_e$. (Notice that the value of $\widetilde f$ at
$x_j^e$ is independent of the orientation of $e$.) 
The following straightforward proposition justifies the name given to
$\div_\ell(f)$.

\begin{prop} For every function $f\colon V(G) \rightarrow \mathbb Z$, the two divisors
 $\div_\ell(f)$ and $\div(\widetilde f)$ are the same. In particular, $\div_\ell(f)$ is a principal divisor on $H$.
\end{prop}

We have the following theorem.

\begin{thm}\label{thm:admissible} Let $D\in \Div(H)$ with support in
$V(G) \subseteq V(H)$. Then:
\begin{itemize}
\item[($i$)] For every $f\colon V(G) \rightarrow \mathbb Z$, 
the divisor $D+ \div_\ell(f)$ is $G$-admissible and linearly
equivalent to $D$ on $H$.
\item[($ii$)] Every $G$-admissible divisor $D' $on $H$
linearly equivalent to $D$ is of the form 
$D+ \div_\ell(f)$ for some $f\colon V(G) \rightarrow \mathbb Z$.
\end{itemize}
\end{thm}

\begin{proof} Let $f\colon V(G) \rightarrow \mathbb Z$ and 
$\widetilde f \colon V(H) \rightarrow \mathbb Z$ the extension of $f$ described above. 
We have $\div_\ell(f) = \div(\widetilde f)$, which
shows that $D+ \div_\ell (f) \sim D$ on $H$. In addition, since
$\div_\ell(f)$ is $G$-admissible, and
$\mathrm{Supp}(D) \subseteq V(G)$,  the divisor
$D+\div_\ell(f)$ is $G$-admissible. This proves the first statement.
 
To prove the second statement, let $D'$ be
a $G$-admissible divisor linearly equivalent to $D$, and let 
$F: V(H)\rightarrow \mathbb Z$ be a function so that
$D' = D + \div(F)$.
Denote by $f$ the restriction of $F$ to  $V(G)$. 
We claim that $F$ coincides with the extension 
$\widetilde f$ of $f$ to $V(H)$, which will obviously prove
$(ii)$. Indeed, let $e=uv$ be an oriented edge of $G$, and for 
$j=1, \dots, \ell_e$, put 
$$
s_j := F(x_j^e) - F(x_{j-1}^e).
$$
Now,  $s_{j+1}-s_j = D'(x_j^e)$ for $j=1,\dots,\ell_e-1$ because $D$ has support in
 $V(G)$. Since $D'$ is 
$G$-admissible, the difference 
 $s_{j+1}-s_j$ is either 0 or 1, being 1 for at most one value of
 $j$. In other words, there are integers $s$ and $\rho$, with
 $1\leq\rho\leq \ell_e$, such that
$$
s_1=s_2=\dots=s_{\rho}=s, \qquad \textrm{and} \qquad s_{\rho+1}=
\dots=s_{\ell_e}=s+1.
$$
Note that
$f(v) - f(u) =F(v) -F(u) = \sum_{j=1}^{\ell_e} s_j = \ell_e s + \ell_e-\rho.$
It follows that $\ell_e-\rho$ is the remainder of the division of 
$f(v)-f(u)$ by $\ell_e$, and a simple verification
shows that the two functions $F$ and $\widetilde f$ take the
same value on all the vertices $x^e_j$.
\end{proof}

\subsection{Characterization of $G$-admissible divisors: general case} 
 We treat now the general
case of a divisor $D \in \Div(H)$ whose support is not 
necessarily in $V(G)$.  
To do so, we need to introduce a ``twisted'' version
of $\div_\ell(f)$, taking into account the values of 
$D$ in $V(H) - V(G)$.

First, to each divisor $D\in \Div(H)$, we associate 
the function $t^D\colon \E(G) \rightarrow \Z$ which takes the value
$t^D_{e}$ on each oriented edge 
$e=uv$ of $G$ given by:
$$t^D_{e} := \sum_{j=1}^{\ell_e-1} (\ell_e-j) D(x^e_j).$$

\begin{prop}\label{tde}
 For every oriented edge $e=uv$ in $\E(G)$, we have
 $$
 t^D_{e} + t^D_{\ol e} = \ell_e \sum_{j=1}^{\ell_e-1} D(x_j^e).
 $$
\end{prop}

\begin{proof}
 Clear from the identity $x_j^e = x_{\ell_e-j}^{\ol e}$.
\end{proof}

\begin{defi}\label{twdr}\rm Let $t\colon \E(G) \rightarrow \Z$
be an integer valued function taking value $t_e$ 
on the oriented edge $e\in \E(G)$. For each $f\colon V(G) \rightarrow \Z$
and $e\in\E(G)$, define 
$$
\delta_{e}(f;t):=\lfloor \frac{f(v)-f(u)+t_{e}}{\ell_e}\rfloor.
$$
\end{defi}

\begin{propdefi}[Canonical extension of functions with respect to a divisor]\label{propdef:1}
 Let $D$ be a divisor on $H$.
 For each $f\colon V(G) \rightarrow \Z$, there is a unique extension 
$\widetilde f\colon V(H) \rightarrow \Z$ such that 
$D+\div(\widetilde f)$ is $G$-admissible. The function $\widetilde f$ 
is called the
\emph{canonical extension} of $f$ with respect to $D$, and is
alternatively characterized by the following properties:
\begin{enumerate}
\item[(1)] For each oriented edge $e=uv$
 of $G$, and each $j=1,\dots,\ell_e-1$, the divisor $D+\div(\widetilde
 f)$ takes value $0$ at $x_j^e$, unless $\ell_e-j$ is the remainder of 
the division of $f(v) -f(u)+ t^D_{e}$ by $\ell_e$, in which case it takes
value $1$.
\item[(2)] $\widetilde f(x_1^e)=f(u)+\delta_e(f;t^D)$.
\end{enumerate}
\end{propdefi}

\begin{proof}  It is enough to prove the existence and uniqueness of a
  sequence of integers $s_1^e,\dots, s_{\ell_e}^e$ for
  each oriented edge $e=uv$ of $G$ satisfying the following two conditions:
\begin{enumerate}
  \item[(i)] For each $j=1, \dots, \ell_e-1$,  the quantity $\rho_j^e$ defined by
$$\rho^e_j:=D(x_{j}^e) + s^e_{j+1} - s^e_j$$
is either 0 or 1, and we have $\rho_j^e=1$ for at most one value of $j$. 
\item[(ii)] $f(v) = f(u) + \sum_{j=1}^{\ell_e} s_j^e$.
\end{enumerate}
Indeed, if $\widetilde f$ exists such that $D+\div(\widetilde f)$ is
$G$-admissible, then put, for each oriented edge $e=uv$ of $G$,
$$
s_j^e:= \widetilde f(x_j^e) - \widetilde f(x_{j-1}^e)\textrm{ for
}j=1,\dots,\ell_e.
$$
Then the sequences $(s_j^e)$ satisfy (i)--(ii). 
Also, the uniqueness of the sequences $(s^e_j)$ implies that of
$\widetilde f$. Conversely, the existence of the $(s^e_j)$ implies
that of $\widetilde f$, as it is enough to put, for each oriented edge
$e=uv$ of $G$, 
$$
\widetilde f(x_j^e):=f(u)+\sum_{i=1}^js_i^e \textrm{ for } j=1,\dots,\ell_e-1.
$$
It is actually necessary to verify that the value of $\widetilde
f$ at the vertex $x_j^e$ does not depend on the orientation of
$e$ for all $e\in\mathbb E$. 
This follows from Condition (ii) and the equality
\begin{enumerate}
\item[(iii)] $s^e_j = -s^{\ol e}_{\ell_e-j+1}$ for $j=1,\dots,\ell_e$,
\end{enumerate}
which we will prove to hold in a moment as a consequence of Conditions (i) and
(ii). Finally, the uniqueness of $\widetilde f$ implies that of the sequences $(s^e_j)$.

Let $e=uv$ be an oriented edge of $G$. Given $s_1^e$, the equations in
Condition (i) express a bijective relation between
the sequences $(s_j^e)$ and  $(\rho_j^e)$. The relation between the  $s^e_j$ and
the $\rho^e_j$ can also be expressed as
\[
s^e_{j+1} = s^e_1-\sum_{i=1}^j D(x_i^e) + \sum_{i=1}^j
\,\rho_i^e\,\,\textrm{ for }j=1,\dots,\ell_e-1.
\]
In particular,
\[
s^e_{\ell_e} =s^e_1 - \sum_{j=1}^{\ell_e-1}D(x_j^e) +
\sum_{j=1}^{\ell_e-1} \rho^e_j.
\]

Let thus $s_1^e,\dots,s^e_{\ell_e}$ be a sequence of integers and 
$\rho_1^e,\dots,\rho_{\ell_e-1}^e$ the corresponding sequence given by
the equations in Condition (i). Then
\begin{align*}
 f(u)+\sum_{j=1}^{\ell_e} s^e_j &= f(u) + \sum_{j=1}^{\ell_e} \Bigl(s_1^e - \sum_{i=1}^{j-1} D(x_i^e) + \sum_{i=1}^{j-1}\rho^e_i\Bigr)\\
 &= f(u) + \ell_e s_1^e - \sum_{j=1}^{\ell_e-1} (\ell_e-j) D(x_j^e) +
   \sum_{j=1}^{\ell_e-1}(\ell_e-j) \rho^e_j\\
&= f(u) + \ell_e(s_1^e-\delta_e(f;t^D))+\ell_e \lfloor \frac{f(v)-f(u)+t_{e}^D}{\ell_e} \rfloor - t_{e}^D + \sum_{j=1}^{\ell_e-1}(\ell_e-j) \rho^e_j\\
 &=f(v)-r_e + \ell_e(s_1^e-\delta_e(f;t^D)) + \sum_{j=1}^{\ell_e-1}(\ell_e-j)\rho^e_j,
 \end{align*}
where $r_e$ is the remainder of the division of $f(v)-f(u)+t_{e}^D$ by
$\ell_e$. Therefore, Condition (ii) is verified for the sequence
$(s_j^e)$ if and only if
$$
\ell_e\delta_e(f;t^D)+r_e=\ell_es_1^e+\sum_{j=1}^{\ell_e-1}(\ell_e-j)\rho^e_j.
$$
We conclude that if Conditions (i) and (ii) are verified, then so are the following two properties:
\begin{enumerate}
\item[(1')] For each oriented edge $e=uv$ of $G$ and each
  $j=1,\dots,\ell_e-1$, the value of $\rho^e_j$ is zero, 
unless $\ell_e-j$ is the remainder of 
the division of $f(v) -f(u)+ t^D_{e}$ by $\ell_e$, in which case the
value is 1;
\item[(2')] $s_1^e=\delta_e(f;t^D)$.
\end{enumerate}
They are the counterpart to Properties (1) and (2) of the corresponding
$\widetilde f$ mentioned in the statement of the proposition. 
They prove the uniqueness of the
$\rho^e_j$ and the $s^e_j$.

On the other hand, defining the $\rho^e_j$ by Property (1') above, and the $s_j^e$
such that Property (2') holds and the equations in Condition (i)
are verified, we get Condition (ii). 

It remains to prove that Condition (iii) follows from Conditions (i)
and (ii). Notice first that, 
for each oriented edge $e=uv$ of $G$ and each
$j=1,\dots,\ell_e-1$, we have $\rho^e_j=\rho^{\ol e}_{\ell_e-j}$, 
a direct consequence of Property (1'), describing the value of $\rho^e_j$, and Proposition~\ref{tde}. This combined with the definition of $\rho^e_j$ given in (i), gives the equation 
$s_{j+1}^e-s_j^e=s^{\ol e}_{\ell_e-j+1}-s^{\ol e}_{\ell_e-j}$ as
well. It will be thus enough to verify that 
$s^e_{\ell_e} = -s_1^{\ol e}$. And, indeed, first note that, by
Property (2') and 
Proposition \ref{tde},
\begin{align*}
- s^{\ol e}_{1}=- \delta_{\ol e}(f;t^D) &= - \lfloor
                                          \frac{f(u)-f(v)+t^D_{\ol
                                          e}}{\ell_e} \rfloor = \lceil
                                          \frac{f(v)-f(u)-t^D_{\ol e}}{\ell_e}\rceil\\
 &=\lceil \frac{f(v)-f(u)+t^D_{e} - \ell_e (\sum_{j=1}^{\ell_e-1}D(x_j^e))}{\ell_e}\rceil\\
 &= \lceil \frac{f(v) -f(u) + t_{e}^D}{\ell_e}\rceil - \sum_{j=1}^{\ell_e-1} D(x_j^e)\\
 &= s^e_1- \sum_{j=1}^{\ell_e-1} D(x_j^e) + \lceil \frac{f(v) -f(u) + t_{e}^D}{\ell_e}\rceil - \lfloor  \frac{f(v) -f(u) + t_{e}^D}{\ell_e}\rfloor.
\end{align*}
On the other hand, from  (1'), we get the equation 
\[
\lceil \frac{f(v) -f(u) + t_{e}^D}{\ell_e}\rceil - \lfloor  \frac{f(v)
  -f(u) + t_{e}^D}{\ell_e}\rfloor = \sum_{j=1}^{\ell_{e}-1}\rho_j^e .
\]
Combining these two equations, we get 
\[
- s^{\ol e}_{1} = s_1^e - \sum_{j=1}^{\ell_e-1} (D(x_j^e) -\rho_j^e) =
s_1^e - \sum_{j=1}^{\ell_e-1} (s_j^e - s_{j+1}^e) =s_{\ell_e}^e,
\]
and (iii) follows. 
\end{proof}

\begin{defi}\rm Let $D$ be a divisor on $H$. For each function 
$f\colon V(G)\to\Z$,  we denote by $\div_\ell(f; D)$ the principal divisor 
$\div(\widetilde f)$ associated to the canonical extension $\widetilde f$ of $f$ with respect to $D$.
\end{defi}

\begin{prop}\label{prop:adm}
Let $D$ be a divisor on $H$ and $f\colon V(G)\to\Z$. Then $D':= D+
\div_\ell(f;D)$ has the following properties: 
\begin{itemize}
 \item[(1)] For each vertex $u$ of $G$, we have 
$D'(u) = D(u) + \sum_{e\in\E\, ;\, \text{\rm t}_e=u}\,\delta_{e}(f;t^D)$. 
 \item[(2)] For each oriented edge $e=uv$ of $G$ such that the
   division of $f(v)-f(u)+t^D_{e}$ by $\ell_e$ has positive remainder
   $r$, we have $D'(x_{\ell_e-r}^e) =1$.
 \item[(3)] For any other vertex $z$ of $H$, we have $D'(z)=0$.
\end{itemize}
\end{prop}
\begin{proof} This is a reformulation of Properties (1)
  and (2) in Proposition-Definition~\ref{propdef:1}.
\end{proof}

The following is a refinement of Theorem~\ref{thm:admissible}. 

\begin{thm}\label{thm:admissible2} Let $D$ be a divisor on $H$. 
\begin{itemize}
 \item[$(i)$] For every $f\colon V(G)\to\Z$, 
the divisor $D+\div_\ell(f;D)$ is $G$-admissible.
 \item[$(ii)$] Every $G$-admissible divisor $D'$ on $H$ linearly 
equivalent to $D$ 
 is of the form $D+\div_\ell(f;D)$ for some $f\colon V(G)\to\Z$.
\end{itemize}
\end{thm}

\begin{proof}
 By Proposition~\ref{prop:adm}, all divisors of the form
 $D+\div_\ell(f;D)$ are $G$-admissible, whence the first statement. 
Conversely, to prove the second statement, 
write $D' = D +\div(F)$ for a function $F\colon V(H)\to\Z$, and denote
by 
$f\colon V(G)\to\Z$ the restriction of $F$ to the vertices of
$G$. Then $F$ is an extension of $f$ such that $D+\div(F)$ is
$G$-admissible. It follows from Proposition-Definition~\ref{propdef:1}
that $F$ is the canonical
extension of $f$ with respect to $D$, and thus $D'=D+\div_\ell(f;D)$.
\end{proof}

For later use, we state the following consequence of
Proposition-Definition~\ref{propdef:1}.

\begin{prop}\label{cor:sumcan} Let $D$ be a divisor on $H$ and $f_1,
  f_2\colon V(G)\to\Z$. Let $f :=f_1+f_2$.
Denote by $\widetilde f$ and $\widetilde f_1$ the canonical
extensions of $f$ and $f_1$ with respect to $D$.
Denote by $\widetilde f_2$ the canonical
extension of $f_2$ with respect to the divisor $D_1:=D +\div_\ell(f_1; D)$.
Then we have 
\[
\widetilde f_1+\widetilde f_2=\widetilde f.
\]
 In particular, 
 \[
\div_\ell(f;D) = \div_\ell(f_1; D)+\div_\ell(f_2, D_1).
 \]
\end{prop}

\begin{proof} 
Note that
\begin{align*}
D+ \div(\widetilde f_1+ \widetilde f_2) =
D+ \div(\widetilde f_1)+ \div(\widetilde f_2) =
D_1+ \div(\widetilde f_2).
\end{align*}
It follows that $D+ \div(\widetilde f_1+ \widetilde f_2)$ is
$G$-admissible. 
Since $\widetilde f_1+ \widetilde f_2$ restricts to
$f_1+f_2$ on the vertices of $G$, it follows from
Proposition-Definition~\ref{propdef:1} that $\widetilde f_1+
\widetilde f_2=\widetilde f$. The second statement is immediate.
\end{proof}

 \subsection{$G$-admissible chip firing} In this section we define the 
notion of $G$-admissible chip firing, and show that any two
$G$-admissible divisors in the same linear equivalence class are
connected via a 
sequence of $G$-admissible chip firing moves. 
 
Let $D$ be a $G$-admissible divisor on $H$, and let $v$ be a vertex of $G$. 
For each oriented edge $e$ of $G$ with $\he_e=v$, define
\begin{equation}\label{eq:2}
 j_e(v) := \begin{cases} j & \qquad \textrm{ if }  D(x_j^e)=1 \textrm{
     for a certain } j\in\{1,\dots,\ell_e-1\},\\
0 & \qquad \textrm{otherwise}.
        \end{cases}
\end{equation}
Define the cut $C_v=C_v(D)$ in $H$ as the subset of
 all the vertices $x_i^e$ for $0\leq i\leq j_e(v)$ for all oriented
 edges $e$ of $G$ with $\he_e=v$.

 \begin{defi}[$G$-admissible chip-firing]\rm Let $D$ be a $G$-admissible divisor, and let $v$ be a vertex of $V(G)$. The  \emph{$G$-admissible chip firing} 
 move of $D$ at $v$, denoted $M_v(D)$, is the divisor on $H$ 
obtained from $D$ after all the vertices in the cut $C_v(D)$ fire.    
 \end{defi}

The terminology is justified by the following  proposition. For each vertex subset $A\subseteq V(G)$, let
$\chi_{\ind A}\colon V(G)\to\Z$ be the characteristic function of $A$, taking
value 1 at $u\in A$ 
and 0 elsewhere. The characteristic function of $A = \{v\}$ is simply
denoted by $\chi_{\indm v}$.

\begin{prop}\label{prop:2}Let $D$ be a $G$-admissible divisor on
   $H$. Then for each vertex $v$, the divisor $M_v(D)$ is
   $G$-admissible. Moreover, $M_v(D) = D + \div_\ell(\chi_{\indm v};D)$.
\end{prop}

\begin{proof}
The first statement is clear. For the second, observe 
first that $M_v(D)=D+\div(\chi_{\indbi{C}{v}})$. Since $M_v(D)$ is admissible,
and $\chi_{\indbi{C}{v}}$ restricts to $\chi_{\indm v}$ on the vertices of $G$, it
follows from Proposition-Definition~\ref{propdef:1} that $\chi_{\indbi{C}{v}}$
is the canonical extension of $\chi_{\indm v}$ with respect to $D$. Thus 
$\div(\chi_{\indbi{C}{v}})=\div_\ell(\chi_{\indm v};D)$.
\end{proof}

We next show that any two linearly equivalent $G$-admissible divisors
are connected to each other by a sequence of $G$-admissible chip
firing moves. In other words:

\begin{prop} Let $D$ be a $G$-admissible divisor on $H$. Then, for 
each $G$-admissible divisor $D'$ linearly equivalent to $D$, 
there exist an integer $N$ and a sequence
$v_1, \dots, v_N$ of vertices of $G$ such that
$D'=M_{v_N}(\dots (M_{v_1}(D))\dots)$.    
\end{prop}
             
\begin{proof}
Let $D'$ be a $G$-admissible divisor linearly equivalent to $D$. 
By Theorem~\ref{thm:admissible2}, there exists a function 
$f\colon V(G)\to\Z$ such that $D' = D+\div_\ell(f;D)$. 
Without loss of generality, we may assume that 
 the minimum value of $f$ is zero. Then there are $v_1,\dots,v_N\in
 V(G)$, possibly repeated, such that $f = \sum_{i=1}^N \chi_{\indmbi{v}{i}}$, where
 $N:=\sum_vf(v)$. Define $f_0:=0$, and for each $i=1, \dots,N$ put 
$f_i := \chi_{\indmbi{v}{1}}+\dots+\chi_{\indmbi{v}{i}}$.
 Define $D_i := D+\div_\ell(f_i;D)$. Since $f_{i+1} = f_i+ \chi_{\indmbi{v}{i+1}}$, 
 by Propositions \ref{cor:sumcan} and \ref{prop:2}, we have
 $$D_{i+1} = D_i +\div_\ell(\chi_{\indmbi{v}{i+1}}; D_{i}) = M_{\indmbi{v}{i+1}}(D_i),$$
which yields
 $$M_{\indmbi{v}{N}}(M_{\indmbi{v}{N-1}}(\dots M_{\indmbi{v}{1}}(D)\dots)) =D_{N} = D'.$$
\end{proof}

%%%Tilings1

\section{Tilings I: uniform edge length one}\label{sec:tiling1}

\subsection{Setting.} Let $G=(V,E)$ be a finite connected graph
without loops.  For the ring  $A= \R, \mathbb Q$ or $\Z$,
we denote by $C^0(G,A)$ the $A$-module of all
functions $f: V \rightarrow A$ and by 
$C^1(G,A)$ the $A$-module of all the 
functions $g: \mathbb E \rightarrow A$ that
verify $g(e) = -\, g(\ol e)$ for every 
oriented edge $e$ of $\E$. Clearly, $C^0(G,A)= A^V$. Also,
the choice of an orientation for $G$ gives an isomorphism
$C^1(G,,A)\simeq A^E$. 
The cochain complex 
$C^*(G,A): C^0(G,A)\stackrel {d}{\longrightarrow} C^{1}(G,A)$ is such
that the differential $d$ sends $f\in C^0(G, A)$ to $d(f)$ defined by 
\[
d(f)(e) := d_{uv}f:= f(v)-f(u)\,\,\forall\, e=uv\in\E.
\]
 
Similarly, we have the free $A$-module $C_0(G,A)$ generated by the
 set of vertices $V$, with generators denoted by $(v)$ for $v\in V$, and the $A$-module $C_1(G,A)$, which is the quotient of the free $A$-module generated by the
 oriented edges of $G$ by the submodule generated by elements of the form
$(e)+(\ol e)$ for each oriented edge $e\in\E$, where $(e)$ denotes the element
associated to $e$ both in the free module and in the quotient. As before, 
$C_0(G, A)\simeq A^V$, whereas $C_1(G,A)\simeq A^E$ under the choice
of an orientation. We have the chain complex 
$C_*(G,A): C_1(G,A)\stackrel{\partial}{\longrightarrow} C_0(G,A)$, where 
the boundary map $\partial$ is defined by
\[
  \partial (e) := (v)-(u) \,\,\forall \, e=uv\in\E.
\]

The spaces $C_i(G,A)$ and $C^i(G,A)$ are canonically dual for
$i=0,1$. More precisely, there are unique natural scalar products
$\langle\,,\rangle$ on
$C_0(G,A)$ and $C_1(G,A)$ satisfying $\langle (u),(v) \rangle=
\delta_{u,v}$ for all pairs of vertices $u,v \in V$, and 
\begin{align*}
\langle (e),(e') \rangle =
\left\{
\begin{array}{rl}
\pm1  & \mbox{if } (e') = \pm (e),\\
0 & \mbox{otherwise}
\end{array}
\right.
\end{align*}
for all pairs of oriented edges $e,e'\in\E$. 
These pairings naturally identify the space $C_i(G,A)$ with
$C^i(G,A)$, for $i=0,1$, in such a way that the adjoint $d^*$
of $d$ gets identified with $\partial$. We have
$\partial \circ d = d^*\circ\, d = \Delta$, 
where $\Delta$ is the Laplacian of the finite graph $G$, defined by 
\[
\Delta(f)(v) = \sum_{\substack{e\in\mathbb E\,\\ \te_e=v}}
\Bigl(\,f(v)-f(\he_e)\,\Bigr)  \,\, \forall \,f\in C^0(G,A)\,,\forall
\,v\in V(G).
\]
Because of the identification, we deliberately use $\alpha_e$
to denote the value of $\alpha\in C^1(G,A)$ at $e\in \E$.
Note that $\alpha_e = -\alpha_{\ol e}$ for each $e\in\E$.

Consider now the lattice $C^1(G, \Z) \subset C^1(G, \R)$ that upon the
choice of an orientation for the edges of $G$
we can identify with the lattice $\Z^E \subset \R^E$.
In this case, the scalar product defined above becomes the
natural Euclidean norm on $\R^E$. 
The Voronoi decomposition of $C^1(G, \Z) \subset C^1(G, \R)$ with
respect to the norm 
$\|\,\,\|$ associated to 
$\langle\,,\rangle$ is thus identified with the standard tiling of
$\R^E$ by hypercubes $\square_{\alpha}$ indexed by vectors 
$\alpha \in \Z^E$ and defined by 
\[
\square_{\alpha}:=\Bigl\{\,x\in \R^E\,\Bigl|\,
 \|x-\alpha\|_{\infty} \leq \frac
 12\,\Bigr\},\,\text{where }
 \|x-\alpha\|_{\infty} := \max_{e\in \E}|x_e-\alpha_e|.
\]

\subsection{Lattice of integer cuts and the Laplacian lattice.} Let $A$ be a ring.
The submodule $\F_A$ of $C^{1}(G,A)$ is defined by the image of $d$: 
$$
F_A:=\,d\bigl( C^0(G,A)\bigr)\subseteq C^{1}(G,A).
$$ 
For a subset $C \subset V$, we recall that the characteristic function 
$\chi_{\ind C}\in C^0(G,A)$ of $C$ takes value one 
at each $v\in C$ and zero elsewhere. Since any
function $f\in C^0(G, A)$ can be written as 
$$
f = \sum_{v\in V} f(v) \chi_{\indm v},
$$
it follows that $\F_A$ is generated by all the functions of the form
$d(\chi_{\ind C})$ for $C \subset V$. In fact, we have the following
result.

\begin{prop}
Let $G=(V,E)$ be a finite connected graph without loops, and let $S \subset V$ be a subset of cardinality 
 $|V|-1$. Then $\F_A$ is a free $A$-module of rank $|V|-1$ and the elements
 $d(\chi_{\indm v})$ for $v\in S$ form a basis of $F_A$.
\end{prop}

\begin{proof}
Since $\sum_{v\in V} d(\chi_{\indm v}) =0$, and since $\F_A$ is generated by
$d(\chi_{\indm v})$ for $v\in V$, 
the set of elements $d(\chi_{\indm v})$ for $v\in S$ generate $\F_A$.
In addition, these are linearly independent. 
In fact, a $f\in C^0(G, A)$ is in the kernel of $d$ if and only if
$f(u)=f(v)$ for all oriented edges $e=uv$ in $\E(G)$, whence  if and
only if $f$ is a constant function by
connectedness of $G$. Thus, if $f=\sum_{v\in V}a_v\chi_{\indm v}$ then
$d(f)=0$ if and only if $a_u=a_v$ for each $u,v\in V$. In particular,
since $S\neq V$, we have 
$$
\sum_{v\in S}a_vd(\chi_{\indm v})=0
\text{ if and only if }a_v=0 \textrm{ for all }v\in S.
$$
\end{proof}

\begin{defi} \rm Let $C_1$ and $C_2$ be disjoint subsets of $V$. Let
  $E(C_1,C_2)$ denote the set of all the edges of $G$ between a vertex
  of $C_1$ and a vertex of $C_2$, and by $\E(C_1,C_2)$ the set of all
  the oriented edges from a vertex of $C_1$ to a vertex of $C_2$. For a
  subset $C$ of $V$, its \emph{edge cut} is $E(C, V-C)$, and its
  \emph{oriented edge cut} is $\E(C,V-C)$. 
If there is no risk of confusion, we drop the word ``edge'' and simply
use ``cut'' and
``oriented cut'' when referring to these two edge sets.  
\end{defi}

Note that we have 
\begin{align*}
d(\chi_{\ind C})(e) :=
\left\{
\begin{array}{rl}
-1  & \quad \mbox{if } e\in \E(C,V-C),\\
1 &   \quad \mbox{if } \ol e \in \E(C , V-C),\\
0 & \quad \mbox{otherwise.}
\end{array}
\right.
\end{align*}
In other words, $d(\chi_{\ind C})$ is the \emph{characteristic function} of the
oriented cut $\E(V-C,C)$ in $C^1(G, A)$: We have 
\[
d(\chi_{\ind C}) = \sum_{e\in \E(V-C,C)} \Bigl(\,\chi_{\indm e} -\chi_{\indmbar{\ol
  e}}\,\Bigr),
\]
where for an oriented edge $e\in \E$, $\chi_{\indm e}$ is the characteristic
function of $e$ taking value 1 on $e$ and value zero on every other
oriented edge.

The following well-known result states that $F_A$
is the subspace of $C^1(G, A)$ orthogonal to the
first graph homology
$H_1(G, A) :=\ker(\partial) \subseteq C_1(G, A)$.
We provide a proof for the sake of completeness. 

\begin{prop} Let $\mu\in C^1(G, A)$. Then $\mu\in
  F_A$ if and only if $\sum_{e \in \gamma}\mu_e = 0$ for each oriented
  cycle $\gamma$ of 
$G$.
\end{prop}

\begin{proof} For each oriented cycle $\gamma$ and $v\in V$, we have 
$$
\sum_{e\in\gamma}
d(\chi_{\indm v})_e=d(\chi_{\indm v})_{e_+}+d(\chi_{\indm v})_{e_-}=-1+1=0,
$$
where $e_+$ (resp.~$e_-$) is the edge of $\gamma$ with tail
(resp.~head) $v$ if $v\in\gamma$. If $v\not\in\gamma$ the sum is also
zero. The ``only if'' statement follows. Conversely, let
  $\mu\in C^1(G, A)$ 
such that $\sum_{e \in \gamma}\mu_e = 0$ for each oriented
  cycle $\gamma$. Fix $v_0\in V$. For each $v\in V$, take an oriented path $P$  
 from $v_0$ to $v$ and define 
 $$f(v) :=\sum_{e\in \E(P)} \mu_e.$$
Note that $f(v)$ is independent of the choice of $P$: Indeed, for
each other oriented path $P'$ 
from $v_0$ to $v$, 
$P'-P$ can be decomposed as as a sum of oriented cycles in $G$, and thus
$$\sum_{e\in P}\mu_e -\sum_{e\in P'}\mu_e = 0.$$
This gives a function $f\in C^0(G,A)$ that verifies $d(f)=\mu$.
\end{proof}

\begin{defi}[{Cut lattice}] \rm Let $G$ be a connected graph without loops. The \emph{lattice of integer cuts},
or simply the \emph{cut lattice}, of $G$
 is the lattice $F_\Z \subset F_\R$. The space $F_\R$ is
 called the \emph{cut space} of $G$ and  an element of 
 $F_\Z$ of the form $d(\chi_{\ind C})$ for a subset $C\subset V$ is called a \emph{cut element}. 
\end{defi}

Note that the cut lattice and cut space live in $C^1(G, A)$ for $A = \Z$ and $\R$, respectively. We can use the operator $d^*$ to bring them down into $C^0(G, A)$. 
\begin{defi}\rm For each ring $A$, we denote by $\Lambda_A$ and
  $H_{0,A}$ the submodules of $C^0(G, A)$ defined by 
\begin{align*}
\Lambda_A :=& d^*\bigl(F_A\bigr) \subseteq C^0(G,A),\\
H_{0,A} :=& \Bigl\{f\in C^0(G, A)\,\Bigl|\, \sum_{v\in
  V}f(v)=0\Bigr\}.
\end{align*}
\end{defi}

\begin{prop} We have $\Lambda_A \subseteq H_{0,A}$. In addition, $d^*$ restricts
  to an isomorphism from $F_A$ onto $\Lambda_A$. 
\end{prop}
\begin{proof}
Since $\sum_{v\in V}f(v)=<f,\sum_{v\in V}\chi_{\indm v}>$ and $d(\sum_{v\in
  V}\chi_{\indm v})=0$, the first statement follows from the fact that 
$\mathrm{Im}(d^*)=\ker(d)^\perp$. The second statement follows from
the fact that $C^1(G, A) = \mathrm{Im}(d) \oplus \ker(d^*)$.
\end{proof}

\begin{defi}[Laplacian lattice]\rm The sublattice $\Lambda_\Z \subseteq H_{0,\Z}$ is called the \emph{Laplacian
lattice.}
\end{defi}
In general, the inclusion $\Lambda_\Z \subseteq H_{0,\Z}$ can be strict. In fact, we have the following proposition which characterizes trees as the only connected graphs for which the equality holds. 
\begin{prop}[Kirchhoff's matrix-tree theorem] The index of $\Lambda_\Z$ in $H_{0,\Z}$ is equal to the number of spanning trees of $G$. 
\end{prop}

 Let $f : V \rightarrow \mathbb R$. By definition,
 $$
 \|d(f)\|^2 = \langle d(f),d(f)\rangle =
 \langle f,d^*d(f) \rangle =
 \langle f,\Delta(f)\rangle,
 $$ 
 where, as before, $\Delta$ is the Laplacian of $G$. 
 So $d^*$ induces an isomorphism $(\F_\R, \|\|) \simeq
 (\Lambda_\R, q)$, where $q$ is the quadratic form on 
 $\Lambda_\R =H_{0,\R}$ induced by $\Delta$, that is

\begin{defi}[Quadratic form $q$ on $H_{0,\R}$]\rm
For each $h \in H_{0,\R}$, define $q(h) := \langle
 f,\Delta(f)\rangle$ for any element $f \in C^0(G, \R)$ with
 $\Delta(f) = h$.
\end{defi}

\begin{defi}[Voronoi decomposition of the Laplacian lattice]\rm The Voronoi decomposition of 
 $(\Lambda_\R,q)$ induced by the sublattice $\Lambda_\Z$ is denoted by
$\Vor_\Delta(\Lambda):=\Vor_q(\Lambda_\Z)$ 
 \end{defi}

 \begin{remark}\rm
  We observe that this decomposition is usually 
different from the decomposition of $H_{0,\R}$ equipped  with the 
restriction of 
  the Euclidean norm of $C^0(G, \R)$ and induced by the integer
  lattice $H_{0,\Z}$ (or by $\Lambda_\Z$). 
 \end{remark}

 From the preceding discussions, it follows that the Voronoi 
decomposition $\Vor_{\|.\|}(\F_\Z)$ of $(\F_\R, \|.\|)$ 
 induced by the cut lattice $\F_\Z$ is isomorphic to $\Vor_\Delta(\Lambda_\Z)$.

Our aim in this section is to provide a detailed description of the Voronoi decomposition of 
$\F_\R$ with respect to $\F_\Z$. As a consequence, we will provide a complete description of the union of 
hypercubes $\bigcup_{\beta\in \F_\Z} \square_\beta$ in $C^1(G, \R)$ in
terms of this Voronoi decomposition. Only later, in the next two
sections, will we generalize the results to the case of a graph $G$ 
equipped with an integer valued length function and the $G$-admissible 
setting of the last section.   
   
We start by giving the following basic characterization of the Voronoi cells in $\F_\R$.

\begin{prop}\label{prop:basic}
 Let $\beta \in \F_\Z$. A point of $\F_\R$ of the form 
$\beta +  \eta$ for $\eta\in \F_\R$ belongs to the $\beta$-centered
Voronoi cell 
$\Vor_{\F}(\beta)$ of $\Vor_{\|.\|}(\F_\Z)$ if and only if it satisfies the 
following set of inequalities:
\begin{equation} \label{eq:cut1}
\textrm{For every subset $S\subseteq V$,} \newline
\end{equation}
 $$-\frac 12 |\E(S, V-S)| \,\,\leq \sum_{e \in \E(S,V-S)} \eta_{e} \,\,\leq \,\,\frac 12 |\E(S, V-S)|.$$
\end{prop}

\begin{proof} If $\beta+\eta \in \Vor_\F(\beta)$, then, for each $S\subseteq V$, 
since $\beta \pm d(\chi_{\ind S}) \in \F_\Z$, 
we must have $\|\eta\|^2 \leq \|\eta \pm d(\chi_{\ind S})\|^2$, from which Inequality~\eqref{eq:cut1} 
follows for $S$.

Suppose now that an element $\eta \in \F_\R$ verifies the above set of
inequalities for subsets $S\subseteq V$. We need to prove that 
$\|\eta- d(f)\|^2 \geq \|\eta\|^2 $ for each $f\in C^0(G,\Z)$. For
such $f$, denote by $n_1<n_2< \dots< n_r$ all 
the values it takes. Without loss of generality, since
$d(f)=d(f+\chi_{\ind V})$, we may assume
that $n_r\leq 0$.  Set 
$n_{r+1}:=0$.
For each $i=1,\dots,r$, define the subset $S_i\subseteq V$ by
$S_i :=\bigcup_{j=1}^i f^{-1}(n_j)$. Then 
$f = \sum_{j=1}^r (n_j -n_{j+1}) \chi_{\indbi{S}{j}}$, which gives
$$
d(f) = \sum_{j=1}^r (n_j -n_{j+1}) d(\chi_{\indbi{S}{j}})
= \sum_{j=1}^r  (n_{j+1}-n_j)
\chi_{\indM{\E(S_{\hspace{-.04cm}\indm j}, V-S_{\hspace{-.04cm}\indm j})}}.
$$
Here
$\chi_{\indM{\E(S_{\hspace{-.04cm}\indm j}, V-S_{\hspace{-.04cm}\indm j})}}$ is the characteristic function of the oriented cut  $\E(S_{\hspace{-.04cm}\indmbar j}, V-S_{\hspace{-.04cm}\indmbar j})$ in $C^1(G, \R)$. 

Note that we have 
$$
\langle d(\chi_{\indbi{S}{j}}), d(\chi_{\indbi{S}{k}}) \rangle = \Bigl|\,\E(S_j,
V- S_j) \cap \E(S_k, V- S_k)\,\Bigr|\geq 0
$$
for each
$j,k=1,\dots,r$. This gives 
\[
\|\eta-d(f)\|^2 -\|\eta\|^2 \geq \sum_{j=1}^r \Bigl(\,(n_{j+1}-n_j)^2
\,\|d(\chi_{\indbi{S}{j}})\|^2- 2(n_{j+1}-n_j)(\sum_{e \in \E(S_j,
  V-S_j)} \eta_{e})\, \Bigr)
\]
Applying Inequality \eqref{eq:cut1} to each $S_j$, we have 
\begin{align*}
 (n_{j+1}-n_j)^2 \,\|d(\chi_{\indbi{S}{j}})\|^2 
&= (n_{j+1}-n_j)^2 |\E(S_j,V-S_j)| \\
 &\geq (n_{j+1}-n_j) |\E(S_j, V-S_j)|\\ 
&\geq  2(n_{j+1}-n_j)(\sum_{e \in \E(S_j, V-S_j)} \eta_{e}),
\end{align*}
which shows that $ \|\eta- d(f)\|^2 - \|\eta\|^2 \geq 0$. Since this holds for every $f\in C^0(G, \Z)$, we have 
$\beta+ \eta \in \Vor_F(\beta)$. 
\end{proof}

It will be convenient for what follows to make the following definition.

\begin{defi}\rm For each subset $A \subseteq \F_\Z$,
  we denote by $\square_A$ the closed subset of $C^1(G,\R)$
  defined by $\square_A:=\bigcup_{\alpha\in A} \square_\alpha$.
\end{defi}
  
 \subsection{The projection map from $\square_{F_\Z}$ to $\Vor_{\Delta}(\Lambda)$.} 
Consider the orthogonal decomposition of
$C^1(G, \mathbb R) = \im(d)\oplus \ker(d^*)$. 
The orthogonal projection $\pi_\F: C^1(G, \mathbb R) \rightarrow \F_\R
= \im(d)$ restricts to a 
map $\square_{\F_\Z}\rightarrow \F_\R$, that composed with $d^*$
yields the  
map  $\square_{\F_\Z} \rightarrow \Lambda_\R$, which is a
restriction of $d^*$. 
If there is no risk of confusion, we denote by $\pi$ the projection
map $\pi_\F$. We have the following theorem. 

 \begin{thm} \label{thm:projection1} The following statements hold true:
 \begin{itemize}
 \item[$(i)$] For each $\beta \in \F_\Z$, the map $d^*$ restricts to a
   linear projection 
from $\square_\beta$ 
 onto the Voronoi cell $\Vor_\Delta(d^*(\beta))$. In addition, $d^*$ maps the interior of   $\square_\beta$ onto the 
 interior of $\Vor_\Delta(d^*(\beta))$.
\item[$(ii)$] 
$(d^*)^{-1}(\Vor_\Delta(d^*(\beta)))\cap\square_{\F_\Z}=\square_\beta$.
 \item[$(iii)$] All the fibers of $d^*\: \square_{\F_\Z} \rightarrow \Lambda_\R$ are compact and contractible. In particular, $\square_{\F_\Z}$ 
 is contractible. 
 \item[$(iv)$] The dual complexes of $\square_{\F_\Z}$ and $\Vor_\Delta(\Lambda)$ are isomorphic. 
 \end{itemize}
\end{thm}

Later on, in Theorem~\ref{thm:projection1bis}, we
will in fact give a more refined statement.
Note that in $(iv)$, the dual complex of $\square_{\F_\Z}$
(resp. $\Vor_\Delta(\Lambda)$) refers to the (combinatorial)
Delaunay dual complex and means the simplicial complex
with the vertex set $\F_\Z$ (resp $\Lambda_\Z$) and with
simplices consisting of all subsets $S \subset \F_\Z$
(resp. $S\subset \Lambda_\Z$) such that the intersection
$\bigcap_{\beta \in S}\square_{\beta}$
(resp. $\bigcap_{\beta \in S}\Vor_{\Delta}(\beta)$) is non-empty. 

The rest of this subsection is devoted to the proof
of this theorem.  Since
$d^*: (\F_{\R},\|.\|) \rightarrow (\Lambda_\R,q)$ is an
isometry which maps $\F_\Z$ to $\Lambda_\Z$,
 it will be enough to prove the first three statements for the projection map 
 $\pi_\F : \square_{\F_\Z} \rightarrow \F_{\R}$.
 The last one is then a direct consequence.  
 
Let $f\in C^0(G, \Z)$ and $\beta := d(f) \in \F_\Z$.  A point $x \in \square_\beta$ is of the form 
 $\beta + \mu$, where for each $e\in \E$ the coordinate $\mu_e$ of $\mu$ verifies 
 $|\mu_e|\leq \frac 12$, that is, $\|\mu\|_\infty \leq \frac 12$.  
 In addition, we have the strict inequality  $\|\mu\|_\infty< \frac 12$ 
 if and only if  the point $x$ belongs to the relative interior of $\square_\beta$. 
 
\begin{lemma}\label{lem:k1} Notations as above,  for each $h\in C^0(G, \Z)$ we have 
  $$
  \|\pi_\F(\mu)\| \leq \|d(f-h) + \pi_\F(\mu)\|.
  $$
 The inequality is strict for each $h$ with $d(h)\neq \beta$, provided that
 $\|\mu\|_\infty < \frac 12.$
 \end{lemma}

\begin{proof} We have
  $$
  \|\pi_\F(\mu)\|^2 = \|\mu\|^2 -  \|\mu - \pi_\F(\mu)\|^2,
  $$ 
and similarly, 
$$
\|d(f-h) + \pi_\F(\mu)\|^2 =
\|d(f-h) + \mu\|^2 - \|\mu - \pi_\F(\mu)\|^2.
$$
So, in order to prove the first statement of the lemma, it will be enough to show that
$$
\forall h\in C^0(G, \Z), \quad
\|\mu\|^2 \leq \|d(f-h) + \mu\|^2.
$$
Now, for each oriented edge $e=uv$ in $\E$, since $d_{uv}(f-h) \in \mathbb Z$ and  
$|\mu_e|\leq \frac 12$, we have $|\mu_e| \leq |d_{uv}(f-h)+\mu_e|$, 
from which the inequality in the lemma 
follows. In addition, we observe that the inequality is strict provided that the 
inequality $|\mu_e| \leq |d_{uv}(f-h)+\mu_e|$ is strict for at least
one oriented edge $e=uv$ in $\E$.

Suppose now that $h\in C^0(G, \Z)$ satisfies $d(h)\neq \beta$, and 
$\|\mu\|_\infty <\frac 12$.   Since $d(h) \neq \beta = d(f)$, 
there exists at least one edge $e=uv$ with $d_{uv}(f-h) \neq 0$, and so 
we have  $|\mu_e| < |d_{uv}(f-g)+\mu_e|$,
given that $|\mu_e| < \frac 12$. Therefore, we get 
$$
\|\pi_\F(\mu)\| < \|d(f-h) + \pi_\F(\mu)\|.
$$
\end{proof}

\begin{proof}[Proof of Theorem~\ref{thm:projection1}(i)] From
  Lemma~\ref{lem:k1} we get directly that $\pi_F\bigl(\square_\beta\bigr)\subseteq 
\Vor_F(\beta)$, and that $\pi_F$ maps the interior of the hypercube
$\square_\beta$ into the interior of $\Vor_F(\beta)$.

We prove now the surjectivity of the projection map 
$\pi_\F: \square_\beta \rightarrow \Vor_{\F}(\beta)$. 
By Proposition~\ref{prop:basic}, a point of $\Vor_{\F}(\beta)$ is of the form 
$\beta +  \eta$ for $\eta\in F_{\R}$
satisfying the following set of inequalities:
\begin{equation}\label{eq:cut}
\forall \, S\subseteq V, \quad 
-\frac 12 |\E(S, V-S)| \,\,\leq \sum_{e \in \E(S, V-S)}
\eta_{e} \,\,\leq \,\,\frac 12 |\E(S, V-S)|\,.
\end{equation}
We want to prove the existence of $\mu \in C^1(G, \R)$ 
with $\|\mu\|_\infty \leq \frac 12$ and $\pi_\F(\mu) =\eta$.
Since the result is trivial for $\eta=0$, we may
suppose that $\eta\neq0$.
We prove this using the max-flow min-cut theorem
in graph theory, cf.~\cite[Chapter III]{Bol} or~\cite[Chapter 7]{BM}.

Let $h:= - d^*(\eta)$. Thus 
$h(v)=\sum_{e\,:\, \he_e=v} \eta_{e}$ for each vertex $v$.
Let $X$ be the set of vertices with $h(v) \geq 0$, and
$Y$ its complementary. Since $\eta \in \im(d) $ and
$\eta \neq 0$, we have $d^*(\eta) \neq 0$. It follows that $X,Y\neq \emptyset$. 

 Now, add a new vertex $s$, called source, and another new vertex $t$,
 called target, to $G$, and add new oriented edges $sv$ for all $v\in X$
 and $ut$ for all $u\in Y$ to obtain a new directed graph $\widetilde G$. 
 Define the capacity $c(sv)$ of $sv$ as $h(v)$, the capacity  $c(ut)$
 of $ut$ as $-h(u)$, and the capacity $c(e)$ of all the oriented
 edges $e\in \E$ as $\frac 12$. 
 By the max-flow min-cut theorem,
 the maximum amount of a flow from the source
 $s$ to the target $t$ in the corresponding network
 is equal to the minimum capacity of a cut separating
 $s$ from $t$ in $\widetilde G$. Recall that the capacity of a cut $S$
  is the sum of the values of $c$ over all the edges with one
  end-point in $S$ and another end-point out of $S$. 
  
  Let $S$ be a cut with $s\in S$ and $t \not\in S$.
  Write $S = \{s\}\sqcup C$ with $C = A \sqcup B$,
  for $A \subseteq X$ and $B \subseteq Y$. Then the capacity of $S$ is given by 
  \[
    c(S) = h(X- A) - h(B) + \frac 12 |\E(C, V-C)|,
  \]
  where $h(D):=\sum_{v\in D}h(v)$ for each $D\subseteq V$. 
We claim that for each such cut $S$, we have $c(S) \geq c(\{s\}) =
\sum_{v\in X} h(v)$. Indeed, writing
$c(\{s\}) = h(A) + h(X- A)$, this amounts to showing that 
  \begin{equation}\label{eq3}
    h(A)+h(B) \leq \frac 12 |\E(C, V-C)|.
  \end{equation}
  Note that by definition, 
  \[
    h(A) = \sum_{v\in A} h(v) =
    \sum_{v\in A}\sum_{e\,:\, \he_e=v} \eta_{e}
    = \sum_{e\in \E(A, V- A)} \eta_e.
  \]
  (For $e=vu$ with $v, u\in A$, we have $\eta_e + \eta_{\ol e}=0$.)
Similarly, we have 
\[
  h(B) = \sum_{w\in B} h(w) =
  \sum_{w\in B}\sum_{e\,:\, \he_e=w} \eta_{e} =
  \sum_{e\in \E(B,V-B)} \eta_e.
\]
  It follows that 
  \[
    h(A) + h(B) = \sum_{e\in \E(C, V-C)} \eta_e,
  \]
so that Inequality~\eqref{eq3} is precisely \eqref{eq:cut} applied to
the cut $C$. 

 What we have proved so far shows that the maximum flow $\alpha$ in
 $\widetilde G$ with the capacity function $c$ on edges is precisely
 $c(\{s\})$, and thus $\alpha$ uses maximum allowed capacity on each
 edge $sv$ and each edge $ut$. 
 The restriction of $\alpha$ to $\E(G)$ satisfies
\begin{center} 
$0\leq \alpha(e) \leq \frac 12$ for every $e\in\E(G)$, and
\end{center}
\begin{center}
 For each vertex $v\in V$, $\sum_{e\,:\,\he_e=v}\alpha(e) -\alpha(\ol e) = h(v).$ 
\end{center}

Define now $\mu(e) := \alpha(e) -\alpha(\ol e)$ for each $e\in\E(G)$. 
Note that $\mu(e) = -\mu(\ol e)$ for each $e\in \E(G)$, and so $\mu
\in C^{1}(G, \R)$. 
It follows from the first equation above that $\|\mu \|_{\infty} \leq
\frac 12$. In addition, for each vertex $v\in V$, we have 
$\sum_{e\,:\,\he_e=v} \mu(e) = h(v) = \sum_{e\,:\,\he_e=v} \eta(e)$,
which shows that $\mu - \eta \in \ker(d^*)$. Since $\eta\in \im(d)$, 
we finally conclude that $\pi_\F(\mu) =\eta$, which proves the
surjectivity.

Since $\pi_F$ maps two different
hypercubes 
onto two different Voronoi cells, surjectivity
of $\pi_F\:\square_{\F_\Z} \rightarrow \F_{\R}$ implies that  
the interior of a hypercube $\square_\beta$ is projected onto the relative interior of $\Vor_F(\beta)$ for each $\beta \in F_\Z$, and Part $(i)$ follows.   
\end{proof}

\begin{proof}[Proof of Theorem~\ref{thm:projection1}$(ii)$] Let 
$\alpha\in\square_{\F_\Z}$ such that
$\pi_F(\alpha)\in\Vor_F(\beta)$. By $(i)$, there is $\mu\in C^1(G,\R)$
with $\|\mu \|_{\infty} \leq \frac 12$ such that 
$\pi_F(\alpha)=\pi_F(\beta+\mu)$. Let $\eta:=\alpha-\beta-\mu$. Then 
$d^*(\eta)=0$. Since $\alpha = \beta +\eta+ \mu$, in order to prove the statement in the theorem, we need to 
show that $\|\eta+\mu\|_{\infty} \leq \frac 12$. 

Now, by assumption, $\alpha \in \square_{\F_\Z}$
and so there is $\beta'\in\F_\Z$ such that 
$\alpha\in\square_{\beta'}$. Set $\beta''=\beta'-\beta$. Then 
$\|\eta+\mu-\beta''\|_{\infty} = \|\alpha-\beta'\|_{\infty}\leq \frac 12$. We  
prove $\|\eta+\mu\|_{\infty} \leq \frac 12$ by showing
that for each oriented edge $e\in\E$,
either $\eta_e=0$ or $\beta''_e=0$. This implies
the result since in the first case
we get $|\eta_e + \mu_e| = |\mu_e|\leq \frac 12$, and in the second
case we get $|\eta_e + \mu_e| = |\eta_e + \mu_e -\beta''_e| \leq \frac 12$.

To prove the claim, note that since $\|\mu \|_{\infty} \leq \frac 12$ and 
$\|\eta+\mu-\beta''\|_{\infty} \leq \frac 12$, we have 
$|\eta_e-\beta''_e|\leq 1$ for each $e\in\E$. Since $\beta_e''\in\Z$, we 
have
\begin{equation}\label{be}
(\beta''_e)^2+(\eta_e-\beta''_e)\beta''_e\geq 0
\end{equation}
for every $e\in\E$. But 
$$
\langle \beta'',\beta''\rangle+\langle \eta-\beta'',\beta''\rangle=
\langle \eta,\beta''\rangle =0
$$
because $d^*(\eta)=0$ and $\beta'' \in \im(d)$. Thus equality holds in \eqref{be} for 
every $e\in\E$, which implies $\eta_e \beta_e'' =0$ and the claim follows. 
\end{proof}

\begin{proof}[Proof of Theorem~\ref{thm:projection1}$(iii)$] Since
$\pi_\F$ is a linear projection map, the fibers are all convex, thus 
contractible. They are all closed and bounded 
since each is included in a hypercube
of the form $\square_\beta$ for some $\beta\in \F_\Z$, by $(ii)$.
The projection map $\pi_F\:\square_{\F_\Z} \rightarrow \F_\R$ has closed 
convex fibers, so $\square_{\F_\Z}$ is contractible. 
\end{proof}

\begin{proof}[Proof of Theorem~\ref{thm:projection1}$(iv)$] This last
  assertion is a direct consequence of $(i)$ and $(iii)$.
\end{proof}

  \subsection{The combinatorics of the Voronoi decomposition 
$\Vor_{\|.\|}(\F_\Z) \simeq \Vor_\Delta(\Lambda)$.}
 Consider the Voronoi decomposition $\Vor(\F_\Z)$ of the lattice $\F_\Z$ in $\F_\R$ 
 under the natural distance defined by $\|.\|$
 and let $\Vor_\F(\beta)$ be the Voronoi cell of the 
 lattice element $\beta \in \F_\Z$. Denote by
 $\FP$ the face poset of the polytope $\Vor_F(O)$, where
 $O$ is the point of origin in $\F_\Z$.  Recall that,
 as a set, $\FP$ consists of all the faces of $\Vor_F(O)$,
 and the partial order is defined by the inclusion between faces.
 In this section, we give a combinatorial description of 
 this poset in terms of a particular class of orientations of subgraphs of $G$. 
 
 Recall that by a spanning subgraph $G'$ of $G$ we mean a subgraph
 with $V(G')=V(G)$.

\begin{defi}[Cut subgraph]\rm
A spanning subgraph $G'$ of $G$ is called a {\it cut subgraph} 
if there exists an integer $s$ and a
partition $V_1,\dots, V_s$ of $V$ such that $G'$ is the 
graph obtained  from $G$ by deleting all the edges which lie in any
one of the $V_i$, i.e., if $E(G') = E(G)-\bigcup_{i=1}^s E(V_i)$. 
 \end{defi}

 Recall that an orientation $D$ of a graph $G$ is
 called {\it acyclic} if $D$ does not contain any 
 oriented cycle. 

 \begin{defi}[Coherent acyclic orientations of subgraphs]\rm 
  An orientation $D$ of a subgraph $G'$ of $G$ is called 
  {\it coherent acyclic} if there exists a partition $V_1,\dots, V_s$ of $V(G)$ such that 
  \begin{itemize}
   \item[-] $G'$ is the cut subgraph of $G$ obtained by removing
     all the edges which lie in any of the $V_i$, and 
  \item[-] the orientation $D$ is that defined by the linear order 
  $1< \dots < s$, i.e., all the edges between 
$V_i$ and $V_j$ for $i<j$ are oriented  from $V_i$ to $V_j$.
  \end{itemize}
\end{defi}

We shall view a coherent acyclic orientation $D$ as a subset 
$\E(D)\subseteq\E$. Denote by $\AC$  the set consisting of all coherent acyclic 
orientations of all
  cut subgraphs of $G$.

  \begin{defi}[Poset structure of $\AC$]\rm 
   We define the following partial order on $\AC$. For 
  two elements $D_1$ and $D_2$ of $\AC$, we say  $D_1\preceq D_2$ 
  if $\E(D_2)\subseteq\E(D_1)$. 
We call the poset $\AC$ the \emph{poset of coherent acyclic orientations of cut subgraphs of $G$}. 
 \end{defi}

The main theorem of this section is the following:
 
\begin{thm}\label{thm:vor} The face poset $\FP$ of $\Vor_F(O)$ is
isomorphic to the poset $\AC$ 
of coherent acyclic orientations of cut subgraphs of $G$. 
\end{thm}

The rest of this section is devoted to the proof of Theorem~\ref{thm:vor}.

\begin{defi}[Positive support]\rm
 For an element $x \in C^1(G,\R)$, the \emph{positive support} of $x$,
 denoted $\supp^+(x)$, is 
 the set of all oriented edges $e\in \mathbb E$ with $x_e >0$.
\end{defi}

 Let $x\in \F_\R$ and $f\in C^0(G, \R)$ with $x =d(f)$. Denote by 
$a_1< \dots< a_s$ all the different values taken by $f$, and define the level sets 
$V_i = f^{-1}(a_i)$. By definition of $d(f)$, the positive support of 
$x$ is the coherent acyclic orientation defined by the ordered
partition $V_1, \dots, V_s$ of $V$.  
In addition, all the elements of $\AC$ are of this form, so:

\begin{prop}\label{prop:cac} $\AC = \{\supp^+(x)\,|\, x \in \F_\R\}$. 
\end{prop}

 \subsubsection{Intersecting Voronoi cells}
 We characterize first the intersecting Voronoi cells in the Voronoi
 decomposition of 
$\F_\R$ induced by the lattice
 $\F_\Z$.  We need the following definition.

 \begin{defi}[Generalized cut] \rm
 A \emph{generalized cut} $\mathfrak C$ in $G$ is an ordered partition
 of the vertex set $V$ 
  into sets $V_1,\dots, V_s$ for $s\geq 2$ 
such that there is no edge between $V_i$ and $V_j$ if $|i -j| \geq
 2$, or in other words, such that all the crossing edges between the
 $V_i$ are 
 between consecutive $V_{j}$ and $V_{j+1}$ for $j=1, \dots, s-1$. We denote by $\E(\mathfrak C)$ the set consisting of all  the oriented
edges with tail in $V_i$ and head in $V_{i+1}$ for any $i=1,\dots,s-1$.
 \end{defi}

\begin{defi}[Generalized cut elements of $F_\Z$]\rm Let $\mathfrak C$ be a generalized cut of $G$ with ordered partition $V_1, \dots, V_s$. 
  The \emph{characteristic function} of $\mathfrak C$
  is the function $\chi_{\ind{\mathfrak C}}\in C^0(G,\Z)$ 
 which takes value $i$ on $V_i$. An element $\beta$ of $\F_\Z$ 
 is called a \emph{generalized cut element} if it is of the form $d(\chi_{\mathfrak C})$ for 
 a generalized cut $\mathfrak C$ in $G$. 
\end{defi}

We have the following characterization of generalized cut elements.
 
 \begin{prop}\label{prop:generalizedcut}
  A $\beta\in\F_\Z$ is a generalized cut element if and only if 
 $\|\beta\|_\infty \leq 1$.
\end{prop}

\begin{proof} For a generalized cut element, we obviously have
  $\|\beta\|_\infty \leq 1$. 
To prove the reverse implication, let $f\in C^0(G, \Z)$
 with $\beta = d(f)$, and let $V_1, \dots, V_s$ be the partition of
 $V$ given by the level sets of $f$, that is, 
$V_i := f^{-1}(a_i)$ for each $i$, where
 $a_1<\dots<a_s$ are 
all the values taken by $f$. Since $\|\beta\|_\infty \leq 1$ and 
 $a_1, \dots, a_s\in \Z$, the ordered partition $V_1, \dots,
 V_s$ defines a generalized cut $\mathfrak C$ 
 of $G$, and we have
 $a_{i+1}-a_i =1$ for each $i=1,\dots, s-1$. 
 It follows that $f = \chi_{\ind{\mathfrak C}} + a_1-1$, and thus
 $\beta=d(\chi_{\ind{\mathfrak C}})$.
\end{proof}

\begin{remark}\rm Note that the ordered partition underlying a
generalized cut $\mathfrak C$ is uniquely determined
by $\E(\mathfrak C)$. Indeed, we have
$\beta = d(\chi_{\mathfrak C})$ for $\beta$ the element
of $F_\Z$ which takes value one on each oriented edge
in  $\E(\mathfrak C)$, and takes value zero on oriented
edges $e$ with neither $e$ nor $\ol e$ in $\E(\mathfrak C)$. This
shows that the characteristic function of the generalized cut
$\chi_{\mathfrak C}$ is uniquely determined, and thus so is the
ordered partition.
\end{remark}

We can now state the following characterization of intersecting Voronoi cells. 

 \begin{lemma}[Intersecting Voronoi Cells]\label{lem:inter2} The Voronoi cell $\Vor_F(\beta)$ for $\beta \in \F_\Z$ 
 intersects the Voronoi cell $\Vor_F(O)$ of the origin if and only if
 $\beta$ is a generalized cut element 
 of $\F_\Z$.  More generally, given $\beta, \lambda\in \F_\Z$, we have that
  $\Vor_F(\beta)\cap \Vor_F(\lambda)\neq \emptyset$ if and only if $\beta-\lambda$ is a generalized cut 
  element.
 \end{lemma}
 
  Before going through the proof, we first recall a basic useful
  property of Voronoi decompositions defined by  
  lattices.

  \begin{prop} For $\beta\in \F_\Z$, the following two conditions are equivalent:
  \begin{itemize}
   \item The two Voronoi cells $\Vor_F(\beta)$ and $ \Vor_F(O)$ intersect.
   \item $\beta/2 \in \Vor_F(\beta) \cap \Vor_F(O)$. 
  \end{itemize}
  \end{prop}

\begin{proof} One direction is obvious. To prove the
other direction, suppose
$\Vor_F(O)\cap \Vor_F(\beta)$ is not empty, and
suppose for the sake of contradiction that
$\beta/2 \notin \Vor_F(\beta) \cap \Vor_F(O)$.
Then there exists $\mu \in \F_\Z$ such that 
$\|\mu-\beta/2\| < \|\beta/2\|$,  or equivalently, 
$\langle \beta- \mu, \mu \rangle > 0$.
Consider  the four lattice points $O$, $\beta$, $\mu$ and 
$\beta-\mu$, and take a point $x \in \Vor_F(O)\cap \Vor_F(\beta)$.  By assumption $\|x\| = \|x-\beta\|$, which gives
$2\langle x,\beta \rangle = \|\beta\|^2$. Since
$\|x\| \leq \|x-\mu\|$, we have
$2\langle x,\mu\rangle \leq \|\mu\|^2$. Combining these two, we get 
$$
2 \langle x , \beta -\mu \rangle \geq \|\beta\|^2
- \|\mu\|^2 = \|\beta -\mu\|^2 + 2\langle \beta
-\mu, \mu \rangle > \|\beta -\mu\|^2.
$$
This shows that  $\|x\|^2 > \|x-\beta+\mu\|^2$,
which contradicts the assumption $x \in \Vor_F(O)$.  
 \end{proof}

 \begin{proof}[Proof of Lemma~\ref{lem:inter2}] The Voronoi cells  
$\Vor_F(O)$ and $\Vor_F(\beta)$ 
 intersect if and only if  $\beta/2 \in \Vor_F(O) \cap \Vor_F(\beta)$,
 or equivalently, 
if and only if  
$\|\beta/2\|\leq \|\beta/2-\mu\|$ for every $\mu\in\F_\Z$, that is, 
$\langle \beta, \mu \rangle \leq \|\mu\|^2$ for every $\mu\in\F_\Z$.

For a generalized cut element $\beta$, we
have $|\beta_e| \leq 1$ for all $e\in \E$. In particular, for $\mu \in \F_\Z$, 
we have
$\langle\beta,\mu\rangle \leq
\sum_{e\in E}|\mu_e| \leq
\sum_{e\in E} \mu_e^2 = \|\mu\|^2$.
This shows that 
$\Vor_F(\beta) \cap \Vor_F(O) \neq \emptyset$. This proves one
direction of the lemma.

To prove the other direction, suppose $\Vor_F(\beta) \cap
\Vor_F(O)\neq \emptyset.$ Let $D =
\supp^+(\beta)$. Proposition~\ref{prop:cac} asserts that $D$ is a coherent 
acyclic orientation, defined by a certain ordered partition $V_1, \dots, V_s$
of $V$. Notice that $\beta_e=0$ for each $e\in\E$ inside one of the
$V_i$. Furthermore, for each $i=1,\dots,s$,  
let  $C_i = V_1 \cup V_2 \cup \dots \cup V_i$ and 
$\mu_i := -d(\chi_{\indbi{C}{i}})$. Since $\E(C_i, V-C_i)
\subset D$, and the $\beta_e$ are all integers, 
we must have 
$$
\langle \beta, \mu_i \rangle = \sum_{e\in \E(C_i, V-C_i)} \beta_e \geq
|\E(C_i, V-C_i)| = \|\mu_i\|^2.
$$
But the reverse inequality holds as well, as $\Vor_F(O)$ and $\Vor_F(\beta)$ 
 intersect. Thus, equality holds, yielding that $\beta_e = 1$ for every
 $e\in \E(C_i, V-C_i)$. Therefore, $\|\beta\|_\infty \leq 1$, and hence
 $\beta$ is a generalized 
cut element by Proposition~\ref{prop:generalizedcut}.

The last assertion follows from the translation
invariance property of the Voronoi decomposition. 
\end{proof}

\subsubsection{Rank of a generalized cut} 
Let $G$ be a connected graph without loops. For each $C\subseteq V$
denote by $G[C]$ the induced subgraph of $G$ with vertex set $C$,
meaning that the
edges in $G[C]$ are those connecting two vertices of $C$.

Let $\mathfrak C$ be a generalized cut in $G$
given by an ordered partition $V_1, \dots, V_s$
of the vertex set $V$ and with the set of oriented edges $\E(\mathfrak C)$. We define the rank  $\kappa(\mathfrak C)$  of
$\mathfrak C$ as follows. 

For each $i=1,\dots, s-1$, let
$C_i :=V_1\cup\dots\cup V_i$ and 
$S_i := V_{i+1} \cup\dots \cup V_s$, and denote by 
$C_{i,1},\dots,C_{i,l_i}$ and $S_{i,1},\dots, S_{i,r_i}$ all the connected components of 
 $G[C_i]$ and $G[S_i]$, respectively. 

 \begin{defi}\rm \label{def:rank}
The \emph{rank} of $\mathfrak C$, denoted $\k(\mathfrak C)$, 
is defined by $\kappa(\mathfrak C) := \sum_{i=1}^{s-1} (l_i+r_i-1)$.  
\end{defi}

 \begin{defi}[Bond]\rm A {\it bond} in $G$ is a generalized cut $\C$
   with $\kappa (\mathfrak C)=1$. An element $\lambda \in F_\Z$ is
   called a {\it bond} if it is of the form 
 $d(\chi_{\ind C})$ for $C\subset V$  such that the cut defined by $C$ is a
 bond.  We denote by $\Bo$ the set of all bond elements of $F_\Z$. 
\end{defi}
 
\begin{prop}\label{bond3} The following three conditions are equivalent for each generalized cut $\mathfrak C$:
 \begin{enumerate}
 \item $\mathfrak C$ is a bond.
 \item $\mathfrak C$ is given by an ordered partition $V_1, V_2$ of $V$ such that $G[V_1]$ and $G[V_2]$ are both connected.
 \item There is no non-empty cut $\E(C, V-C)$ in $G$ properly
   contained in 
$\E(\mathfrak C)$.
 \end{enumerate}
 \end{prop}

\begin{proof} Let $\C$ be a generalized cut. We retain the terminology
  as in the paragraph preceding Definition~\ref{def:rank}, so
  $\mathfrak C$ is given by an 
ordered partition $V_1, \dots, V_s$ of $V$.

Suppose $\emph{(1)}$ holds, so $\k(\mathfrak C) =1$. Since all the terms $l_i+r_i-1$ in the definition of $\k(\mathfrak C)$ are at least 1 for 
each $i=1,\dots,s-1$, we must have $s=2$ and $l_1=r_1=1$. Thus $\emph{(1)}$ implies $\emph{(2)}$. 

We prove now that $\emph{(2)}$ implies $\emph{(3)}$. Suppose
$\mathfrak C$ is given by the ordered partition $V_1, V_2$ with
$G[V_1]$ and $G[V_2]$ connected. Let 
$C \subsetneq V$ be a proper non-empty subset of $V$ such that 
$\E(C,V-C) \subset \E(V_1, V_2)$. Then $C \cap V_1 \neq
\emptyset$. Moreover, since $G[V_1]$ is connected, we must have $V_1
\subseteq C$, because otherwise $\E(C, V-C)$ would contain the
non-empty set 
$\E(C\cap V_1, V_1- C\cap V_1)$, which is disjoint from $\E(V_1,V_2)$,
clearly not possible. 
Similarly, since 
$G[V_2]$ is connected, $V_2 \subseteq V-C$.
It follows that $C=V_1$, and hence $\E(C,V-C)=\E(\mathfrak C)$. 

Finally, to prove $\emph{(3)}$ implies $\emph{(1)}$ note that all the
cuts of the form 
$\E(C_{1, j}, V-C_{1,j})$, for $j=1,\dots, l_1$, are
included in $\E(\mathfrak C)$. This proves $s=2$ and
$\ell_1 =1$. By a similar reasoning, the cuts of the form 
$\E(V-S_{1,j}, S_{1,j})$ are included in $\E(\mathfrak C)$. This
proves $r_1=1$, and the result follows. 
\end{proof}

While not needed in the sequel, we state the following
result which gives lower and upper bounds on the rank function.

\begin{prop}
  Notations as above, let $\mathfrak C$ be a generalized
  cut with underlying ordered partition $V_1, \dots, V_s$. We have 
\[s-1 \leq \k(\mathfrak C) \leq \sum_{j=1}^s c_j-1\]
where $c_j$ denotes the number of connected components
of the induced subgraph $G[V_j]$ for $j=1, \dots, s$.
\end{prop}
\begin{proof}  The inequality on the left follows from
the observation that $l_j+r_j-1 \geq 1$ for each $j$.
To prove the inequality on the right, consider the
graph $H$ obtained from $G$ by contracting every
edge of $G$ with both extremities in $V_j$ for every
$j=1,\dots, s$. The quantity $\sum_{j=1}^s c_j$ coincides
then with the number of vertices of $H$. Now,
removing edges from $G$ can only result in a
possible increase in the value of $\k(\mathfrak C)$.
Let $G'$ be a subgraph of $G$ obtained by removing
some of the edges crossing between the parts of the
partition such that the resulting subgraph $H'$ of $H$
obtained after contraction of edges within the $V_j$
is a  tree. For the graph $G'$, the corresponding quantity $l'_j+r'_j-1$
is precisely the number of edges of $G'$ between
$C_j$ and $S_j$. Thus, we get 
\[
  \k(\mathfrak C) \leq \sum_{j=1}^{s-1} (l'_j+r'_j-1) = |E(H')| =
  |V(H')|-1 = \sum_{j=1}^s c_j -1.
\]
In the above equation, the first equality follows from
$\frak C$ being a generalized cut, which implies that
each edge of $E(H')$ is counted exactly once in the sum on the
left-hand side, and the second equality follows from $H'$ being a tree. 
\end{proof}

\subsubsection{Description of the codimension one faces of $\FP$}
    
  Let $\sigma$ be a face of $\Vor_F(O)$ contained in another Voronoi
  cell $\Vor_F(\beta)$. By Lemma~\ref{lem:inter2}, we have $\beta =
  d(\chi_{\ind{\mathfrak C}})$ for a generalized cut $\mathfrak C$ in $G$. 

\begin{lemma}\label{lem:codim} Let $\beta:= d(\chi_{\ind{\mathfrak C}})$ for
  a generalized cut $\mathfrak C$ in $G$.  Then the intersection
  $\Vor_F(O) \cap \Vor_F(\beta)$ is contained in an affine
  subspace of $F_\R$ of codimension 
$\kappa(\mathfrak C)$. In particular, if $\sigma$ is a facet of
$\Vor_F(O)$ in $\Vor_F(O)\cap
\Vor_F(\beta)$, then 
$\kappa(\mathfrak C) =1$, and so $\mathfrak C$ is a bond in $G$. 
\end{lemma}
 
\begin{proof} A point $x$ in the intersection $\Vor_F(O)\cap
  \Vor_F(\beta)$  is characterized by the following two properties:
\begin{itemize}
\item[I.] $\|x\| = \|x-\beta\|$, or equivalently, $2\langle x,\beta\rangle = \|\beta\|^2$.
\item[II.] $2\langle x,\mu\rangle \leq \|\mu\|^2$ for each 
$\mu \in F_\Z$. 
\end{itemize}

We start by making the following useful observation: Suppose we can
express $\beta= \sum_{i=1}^{k} \beta_i$, where
each $\beta_i$ is in $F_\Z$ and satisfies 
the following two conditions:
\begin{itemize}
\item[$(i)$] $\beta_i = d(\chi_{\indbi{X}{i}})$ for a cut $X_i$ in $V$,
\item[$(ii)$] $\supp^+(\beta_i) \subseteq \supp^+(\beta)$.
\end{itemize}
Then, since $\|\beta\|_\infty \leq 1$, we have $\supp^+(\beta_i) \cap
\supp^+(\beta_j) =\emptyset$ for different indices $i$ and $j$, and so 
$\|\beta\|^2 = \sum_{i=1}^k \|\beta_i\|^2$.
Applying  II to $\mu=\beta_i$ for each $i=1,\dots,k$, and using I,
it follows that 
$2\langle x,\beta_i\rangle = \|\beta_i\|^2$ for each $i=1,\dots, k$.

Assume now the generalized cut $\mathfrak C$ is given by an ordered
partition $V_1, \dots, V_s$ of $V$, for $s\geq 2$, and for each
$i=1,\dots,s-1$, let $C_{i,1},\dots,C_{i,l_i}$ and $S_{i,1}, \dots,
S_{i,r_i}$ denote the connected components of 
 $C_i = V_1\cup\dots\cup V_i$ and of $S_{i} = V_{i+1}\cup
 \dots\cup V_s$, respectively, as before.
 We have the following equalities:
 $$
 d(\chi_{\ind{\mathfrak C}}) = -\sum_{i=1}^s d(\chi_{\indbi{C}{i}}),
 $$
 $$
 d(\chi_{\indbi{C}{i}}) = \sum_{j=1}^{l_i} d(\chi_{\indbi{C}{i,j}}),
 \quad d(\chi_{\indbi{C}{i}}) =
 - \sum_{j=1}^{r_i} d(\chi_{\indbi{S}{i,j}}).
 $$  

Then there exists a decomposition of the form $\beta= \sum_{h} \beta_h$, where the
$\beta_h$ are the $-d(\chi_{\indbi{C}{i,j}})$  ordered arbitrarily, verifying
properties 
$(i)$ and $(ii)$ above. Similarly, there exists a decomposition of the form
$\beta= \sum_{h} \beta_h$, 
where the $\beta_h$ are the $d(\chi_{\indbi{S}{i,j}})$ ordered arbitrarily, verifying
properties 
$(i)$ and $(ii)$. It follows from the observation we made above that, for each $i=1,\dots,s-1$,
\begin{align*}
  (*) \qquad -2\langle x,d(\chi_{\indbi{C}{i,j}}) \rangle &=
                    \|d(\chi_{\indbi{C}{i,j}})\|^2 = |\E(C_{i,j}, V-C_{i,j}|,
                    \,\,\textrm{for } j=1,\dots l_i \quad \textrm{and} \\
  2\langle x,d(\chi_{\indbi{S}{i,j}}) \rangle &=
                    \|d(\chi_{\indbi{S}{i,j}})\|^2 =
                    |\E(V-S_{i,j}, S_{i,j})|,
                    \,\,\textrm{for } j=1,\dots r_i.
\end{align*}
Now, for each $i=1,\dots,s-1$, the cut elements
$d(\chi_{\indbi{C}{i,k}})$ for $k=1,\dots,l_i$ and $d(\chi_{\indbi{S}{i,j}})$ for
$j=1,\dots,r_i$ span a vector space $H_i$ of dimension $l_i+r_i-1$,
which is isomorphic to 
the cut space of 
the connected graph on $l_i+r_i$ vertices 
obtained from $G$ by contracting each of the subgraphs 
$C_{i,k}$ and $S_{i,j}$ in $G$. In addition, for different
$i,j\in\{1,\dots,s-1\}$, 
the spaces $H_i$ and $H_j$ are orthogonal.
It follows that a point $x$ satisfying all the equations in (*) above 
must lie in an affine plane of codimension
$\sum_{i=1}^{s-1} (l_i+r_i-1)$. The lemma follows.
\end{proof}

 Recall that $\Bo$ denotes the set of all the bond elements in $\F_\Z$.
 
 \begin{prop}\label{bondfacet} The following statements are true:
   \begin{itemize}
              \item[$(i)$]  The bond elements form a system of generators of $\F_\Z$. 
              \item[$(ii)$]   There is a bijection between the bond
                elements $\beta$ of $\F_\Z$ and the facets $\sigma$
                of $\Vor_F(O)$ given by
                $\beta \mapsto \sigma:= \Vor_F(O) \cap \Vor_F(\beta)$.
              \item[$(iii)$] The one-skeleton of the Delaunay dual of
                the Voronoi decomposition $\Vor_{\|.\|}(\F_\Z)$ is
                isomorphic to $\Cay(\F_\Z, \Bo)$, where $\Cay(\F_\Z, \Bo)$
                is the Cayley graph of $\F_\Z$ with respect to the system of 
              generators $\Bo$.
             \end{itemize}  
 \end{prop}
  
 \begin{proof}  $(i)$ The elements of the form $d(\chi_{\indm v})$
   for $v$ a vertex in $V$ form a system of generators for $\F_\Z$. If $G[V - v]$    consists of connected components $C_{v,1}, \dots, C_{v,k}$, then
   each element $d({\chi_{\indbi{C}{v,i}}})$ is a bond, and we have 
  $d(\chi_{\indm v}) = - \sum_{i} d(\chi_{\indbi{C}{v,i}})$. This
  shows that the bond elements form a system of generators.
  
\medskip
  
\noindent $(ii)$  We first prove that for any bond element $\beta =
d(\chi_{\ind C})$, for a subset $\emptyset \subsetneq C \subsetneq V$, the
intersection $\Vor_\F(\beta)\cap \Vor_\F(O)$ has codimension one in
$\Vor_\F(O)$, which means it is of dimension $|V|-2$. As we have
observed before, the intersection $\Vor_\F(\beta) \cap \Vor_\F(O)$ is
contained in the affine hyperplane $2\langle x, \beta \rangle =
\|\beta\|^2 = |\E(C, V-C)|$. 
  
Consider the two graphs $G_1:=G[C]$ and $G_2:=G[V-C]$. By 
Proposition~\ref{bond3}, they are connected. For $i=1, 2$,
define $\F_{i,\R} = \im(d\colon C^0(G_i,\R) \to C^1(G_i,\R))$,
and consider the corresponding sublattice $\F_{i,\Z}$ of full rank in $F_{i,\R}$. 
  Consider two points $x_1$ and $x_2$ in the interior of the Voronoi 
  cells $\Vor_{F_1}(O)$ and $\Vor_{F_2}(O)$ of the origin in
  $\F_{1,\R}$ and $\F_{2,\R}$,  respectively. Extending by zero on the
  edges outside the edges of $G_i$, we can identify the $C^1(G_i,\R)$ as subspaces of 
  $C^1(G, \R)$, and thus make sense of $\F_{1,\R}+\F_{2,\R}$.
  We claim that the point $\pi_F(x_1+x_2) + \frac 12d(\chi_{\ind C})$
  lies in the intersection $\Vor_F(O)\cap \Vor_F(\beta)$.
  (Recall that $\pi_F$ is the orthogonal projection to $\F_\R$
  in the orthogonal decomposition $C^1(G, \R) = \F_\R \oplus \ker(d^*)$.)
  
  To prove this, first note that since $x_i$ and $\beta$ have
  disjoint support, we have $\langle x_i, \beta \rangle =0$, and so we get
  $$
  \|\beta\|^2= 2\langle x_1+x_2 + \beta /2, \beta\rangle =
  2\langle \pi_F(x_1+x_2) + \beta/2, \beta \rangle.
  $$
  Thus, to prove the claim it will be enough to show that for any $\mu \in F_\Z$, 
  $$
  2 \langle  \pi_F(x_1+x_2) + \beta/2, \mu\rangle =
  2 \langle  x_1+x_2 + \beta/2, \mu\rangle \leq \|\mu\|^2.
  $$
  Each $\mu \in \F_\Z$ is of the form $d(f)$ for some
  $f\: V \rightarrow \mathbb Z$ with $\min_{v\in V} f(v)=0$.
  Let $N = \max_{v\in V}\{f(v)\}$.
  We can write $f= f_1+ \dots+f_N$ where each $f_i$ takes
  value $1$ on $f^{-1}([i,N])$ and zero elsewhere. Then
  $\mu =d(f) = \sum_{i=1}^N d(f_i)$ and
  $\sum_{i=1}^N \|d(f_i)\|^2 \leq \|\mu\|^2$,
  so the inequality
  $2 \langle  x_1+x_2 + \beta/2, \mu\rangle \leq \|\mu\|^2$
  will follow from the inequalities
  $$
  \forall X \subset V, \quad
  2 \langle  x_1+x_2 + \beta/2, d(\chi_{\ind X})\rangle
  \leq \|d(\chi_{\ind X})\|^2,
  $$
  that we prove now. 

Consider  $X \subset V$, and let $Y := X \cap C$ and $Z := X \cap
V-C$. We have 
$$
\E(X, V-X) = \E(Y, C - Y) \sqcup \E(Y, (V-C) - Z) \sqcup
\E(Z, (V-C) - Z) \sqcup \E(Z, C - Y).
$$
Since $x_1$ lies in the Voronoi cell of the origin in
$F_{1,\R}$ with respect to $\F_{1,\Z}$, it follows from
Proposition~\ref{prop:basic} that
$$
2 \sum_{e\in \E(C-Y, Y)} x_1(e) \leq |\E(Y, C - Y)|.
$$
Analogously,
$$
2 \sum_{e\in \E((V-C)-Z, Z)} x_2(e) \leq |\E(Z, (V-C) - Z)|.
$$
On the other hand, for every $e\in \E(V-C, C)$,
we have $x_1(e)+x_2(e) + \beta_e/2= 1/2$. These all together show that 
  \begin{align*}
   2 \langle  x_1+x_2 + \beta/2, d(\chi_{\ind X})\rangle =& 2 \sum_{e\in \E(C-Y, Y)} x_1(e) + 2 \sum_{e\in \E((V-C)-Z, Z)} x_2(e) \\
   &+  |\E(Y, (V-C)- Z)| - |\E(Z, C- Y)| \\
    \leq &|\E(Y, C - Y)|+|\E(Z, (V-C) - Z)|\\
    &+  |\E(Y, (V-C) - Z)| - |\E(Z, C - Y)| \\
   \leq &|\E(X, V-X)| = \|d(\chi_{\ind X})\|^2,
  \end{align*}
 and the claim follows.
  
 Finally to show that $\Vor_F(\beta)\cap \Vor_F(O)$
 has codimension one in $\Vor_F(O)$, note that the points of the form 
 $\pi_F(x_1+x_2)+\frac 12d(\chi_{\ind C})$ for
 $x_1$ and $x_2$ in the interior of the Voronoi cells
 $\Vor_{F_1}(O)$ and $\Vor_{F_2}(O)$, form an open set
 around $\beta/2$ in the affine hyperplane   
 $2\langle x, \beta \rangle = \|\beta\|^2 = |\E(C, V-C)|$.
 Indeed, $\pi_F(x_1+x_2) = d(h)$ where $h$ is a
 solution of the equation $\Delta(h) = d^*(x_1+x_2)$. 
 Since $\Delta : C^0(G, \R)/\R\chi_{\ind V} \rightarrow H_{0,\R}$
 is nondegenerate, and the points $d^*(x_1+x_2)$ 
  form an open subset of a hyperplane of $H_{0,\R}$ (consisting of those $f\colon V \rightarrow \R$ with $\sum_{v\in C} f(v) = \sum_{v\in V-C}f(v)=0$), the claim follows.   
  
  To conclude, it remains to show that the association 
$\beta\mapsto\Vor_F(O) \cap \Vor_F(\beta)$ is bijective. The
surjectivity follows by applying Lemma~\ref{lem:codim},
which shows that each facet of $\Vor_F(O)$
is shared with the Voronoi cell of a bond of $\F_\Z$. The injectivity follows from the fact that
$\Vor_F(O) \cap \Vor_F(\beta)$ spans the
affine hyperplane $2\langle x, \beta \rangle = \|\beta\|^2$, which
determines $\beta$. 
             
\medskip
               
\noindent $(iii)$ The last assertion follows from the definition
of the Delaunay dual, and the classification of the facets of Voronoi
cells in the Voronoi decomposition $\Vor_{\|.\|}(\F_\Z)$, given by $(ii)$.
\end{proof}

\subsubsection{Local hyperplane arrangement defined by bonds}
For each bond element $\beta \in \Bo$, let $\F_\beta$ be
the hyperplane in $\F_\R$ defined by the equation 
$$
\F_\beta:=\bigl\{\,x\in \F_\R\,|\, 2\langle x,\beta  \rangle
= \|\beta\|^2\,\bigr\}.
$$
More generally, for an arbitrary lattice point $\mu \in \F_\Z$
and a bond element  $\beta \in \Bo$, let $\F_{\mu,\beta}$
be the affine hyperplane $\mu+\F_\beta$. It follows from
Proposition~\ref{bondfacet} that for each $\mu \in \F_\Z$ the interior of the Voronoi cell
$\Vor_F(\mu)$ is the open cell containing $\mu$
of the hyperplane arrangement given by the
$\F_{\mu,\beta}$ for $\beta \in \Bo$.

\subsubsection{Definition of the map $\phi$ from $\FP$ to $\AC$.} 
Each face $\frak f$ in  the face poset $\FP$ of $\Vor(O)$ is the
intersection of some of the hyperplanes $\F_\beta$ for $\beta
\in\Bo$. This is the key to defining $\phi\colon\FP\to\AC$.

\begin{prop}\label{prop:sumofbonds}
  For each $\alpha \in \mathbb F_\Z$ there exist
  bond elements $\mu_1, \dots, \mu_k$ with
 $\mathrm{supp}^+(\mu_i) \subset \supp^+(\alpha)$ such that 
 $\alpha = \mu_1+\dots+\mu_k$. In particular,
 $\|\alpha\|^2=\sum\|\mu_i\|^2$. 
\end{prop}

\begin{proof}
  Let $f\in C^0(G, \Z)$ with $\alpha=d(f)$. We may
  suppose that $\min f=0$. Let $N = \max f$. 
 For each $i\in \N$ let $S_i :=f^{-1}([i,N])$. Then 
$f =\sum_{i=1}^{N} \chi_{\indbi{S}{i}}$, and hence $\alpha = \sum_{i=1}^{N}
d(\chi_{\indbi{S}{i}})$. Since $\supp^+(d(\chi_{\indbi{S}{i}})) \subseteq
\supp^+(\alpha)$ for every $i$, it is enough to prove the proposition
for
$\alpha=d(\chi_{\ind S})$ for a subset $S\subset V$.
 If $G[S]$ and $G[V - S]$ are connected, then $\alpha$ is a bond
 element. If $G[S]$ is not connected, denoting by $X_1, \dots, X_r$ 
 the connected components of $G[S]$, we have $\alpha =
 \sum d(\chi_{\indbi{X}{i}})$ and
 $\supp^+(d(\chi_{\indbi{X}{i}}))\subseteq
 \supp^+(d(\chi_{\ind S}))$ for each $i$.  We may thus assume $G[S]$
 is connected. Now, if $G[V-S]$ is not connected, denoting by $Y_1,
 \dots, Y_r$  the connected components of $G[V-S]$, we have
 $$
 d(\chi_{\ind S}) =-\sum d(\chi_{Y_i})=\sum d(\chi_{Z_i}),
 $$
 where $Z_i:=V-Y_i$ for each $i$, and the $d(\chi_{\indbi{Z}{i}})$ are bond
 elements satisfying $\supp^+(d(\chi_{\indbi{Z}{i}}))\subseteq
 \supp^+(d(\chi_{\ind S}))$.
\end{proof}

\begin{defi}\rm
 For two elements $\lambda$ and $\mu$ of $C^1(G,\Z)$, we say that 
$\supp^+(\lambda)$ and $\supp^+(\mu)$ are \emph{consistent in their orientations} if there is no oriented edge $e\in \E$ with $e\in \supp^+(\lambda)$ and $\ol e \in \supp^+(\mu)$.
\end{defi}

We have the following result.

\begin{lemma}\label{lem:cle2}
  Let $\frak f$ be a face of $\Vor_\F(O)$, and let $\lambda$
  and $\mu$ be two different bond elements of $\Bo$ such
  that $\frak f \subset \F_\lambda \cap \F_\mu$. Then
  $\supp^+(\lambda)$ and $\supp^+(\mu)$ are
  consistent in their orientations. 
\end{lemma}

\begin{proof}
  Let $\E_\lambda = \supp^+(\lambda)$ and
  $\E_\mu = \supp^+(\mu)$, and for the sake of
  contradiction suppose there exists an oriented edge
  $e\in \mathbb E$ with $e \in \E_\lambda$ and $\ol e \in \E_\mu$.  
 We claim that for each point $x$  in $\Vor_\F(O)\cap \F_\lambda \cap
 \F_\mu$ (and so for each point $x\in \frak f$),
 there exists an element $\beta \in \Bo$ with $\|x-\beta\|<\|x\|$,
 which gives a contradiction. 

Indeed, let $x \in \Vor_F(O) \cap \F_\lambda \cap \F_\mu$. We have 
\begin{equation*}
  2\langle x,\lambda  \rangle = \|\lambda\|^2=|\E_\lambda|
  \qquad \textrm{and} \qquad
  2\langle x,\mu  \rangle = \|\mu\|^2 = |\E_\mu|.
\end{equation*}
Thus we get
\begin{equation}\label{eq1}
2\langle x, \lambda + \mu  \rangle = |\E_\lambda| + |\E_\mu|.
\end{equation}
Let $D = \supp^+ (\lambda+\mu)$. Observe that neither $e$ nor $\ol e$
belongs to $D$ and so the sum $\sum_{f\in D} \lambda_f+\mu_f$ is
strictly less than $|\E_\lambda|+|\E_\mu|$. On the other hand, by
Proposition~\ref{prop:sumofbonds}, the element
$\lambda+\mu$ can be written as a sum $\beta_1+\dots+\beta_k$
such that each $\beta_i$ is a bond element of $F_\Z$ and
$\supp^+(\beta_i) \subset D$. This shows that 
\begin{equation}\label{eq2}
  \sum_i \|\beta_i\|^2 = \|\lambda+\mu\|^2=
  \sum_{e\in D} \lambda_e+\mu_e < |\E_\lambda| + |\E_\mu|.
\end{equation}
If $2\langle x,\beta_i \rangle \leq \|\beta_i\|^2$ for every $i$,
then by summing up all these inequalities and applying
Inequality~\eqref{eq2}, we would have
\begin{equation}
2\langle x,\lambda+\mu\rangle = \sum_i  2\langle x,\beta_i \rangle  \leq \sum_i \|\beta_i\|^2 < |\E_\lambda| + |\E_\mu|, 
\end{equation}
contradicting Equation~\eqref{eq1}. Thus, there exists an $i$ such that $2\langle x,\beta_i \rangle > \|\beta_i\|^2$, which implies $\|x\| > \|x-\beta_i\|$, contradicting the assumption that 
$x \in \Vor_\F(O)$. 
\end{proof}

 Let $\frak f$  be a face of $\Vor_F(O)$, and define 
$$\mathcal U_{\frak f} := \{\beta \in \Bo| \frak f \subset \F_\beta\}.$$ 
\begin{defi} \rm
 We define the map $\phi: \FP \rightarrow \AC $ by 
$$\forall\, \frak f \in \FP, \quad \phi(\frak f ) := \bigcup_{\beta \in \mathcal U_{\frak f}} \supp^+(\beta).$$
\end{defi}

To see that the image of $\phi$ lies in $\AC$, note that by Lemma~\ref{lem:cle2}, $\phi(\frak f)$ coincides with 
$\supp^+(\sum_{\beta \in \mathcal U_{\frak f}}\beta )$ which lies in
$\AC$ by Proposition~\ref{prop:cac}. In addition, $\phi$ is clearly
order preserving, that is, a homomorphism of posets, as
required in Theorem~\ref{thm:vor}.

\smallskip

The following proposition characterizes $\mathcal U_{\frak f}$ from
$\phi(\frak f)$.

\begin{prop}\label{prop:cle3}
 Let $\frak f$ be a face of $\FP$. For every bond element $\beta \in \Bo$, we have $\supp^+(\beta) \subseteq \phi(\frak f)$ 
 if and only if $\beta \in \mathcal U_{\frak f}$.  
\end{prop}

\begin{proof} One direction is trivial by the definition of
  $\phi(\frak f)$. Let $\beta \in \Bo$ be an element with
  $\supp^+(\beta) \subseteq \phi(\frak f)$. We need show
  that $\frak f \subset \F_\beta$. 
  Let $\mu = \sum_{\lambda \in \mathcal U_{\frak f}} \lambda$.
  Since $\supp^+(\mu) = \phi(\frak f)$ and
  $\supp^+(\beta) \subseteq \phi(\frak f)$, we have
  that $\supp^+(\mu-\beta)\subseteq\phi(\frak f)$. Applying
  Proposition~\ref{prop:sumofbonds}, there exists a decomposition 
  $\mu = \sum_{i=0}^k \mu_i$ with $\mu_0=\beta$ and $\mu_i \in \Bo$
  with $\supp^+(\mu_i) \subseteq \phi(\frak f)$ for each $i$.
For each $x\in \frak f$ and $\lambda \in \mathcal U_{\frak f}$, we have 
$$
2\langle x, \lambda \rangle =\|\lambda\|^2
= \sum_{e\in \supp^+(\lambda)} \lambda(e),
$$
which shows that 
$$
2\langle x, \mu \rangle = \sum_{e\in \supp^+(\mu)} \mu(e).
$$

On the other hand, since $\frak f\subset \Vor_F(O)$, we have for each $\mu_i$,  
$$
2\langle x, \mu_i \rangle \leq \|\mu_i\|^2 =
\sum_{e\in \supp^+(\mu_i)} \mu_i(e).
$$
Summing up all these inequalities, and using the fact that 
$\sum_{i=0}^k \mu_i(e) = \sum_{\lambda \in \mathcal U_{\frak f}}
\lambda(e) =\mu(e)$ for every $e\in\E$, we
obtain that all the inequalities above are indeed equalities. 
This shows that 
$$
2\langle x, \beta \rangle  = 2\langle x, \mu_0\rangle
= \|\mu_0\|^2=\|\beta\|^2,
$$ 
and thus $\beta \in \mathcal U_{\frak f}$.
\end{proof}

Let $D$ be a coherent acyclic orientation of a cut subgraph of $G$. Let 
$$
\mathcal U_D := \bigl\{\,\beta\in \Bo\,\,|\,\,
\supp^+(\beta) \subseteq \E(D)\,\bigr\},
$$
and define 
$$
F_D:= \cap_{\beta \in \mathcal U_D } \F_\beta.
$$

As a direct corollary of Proposition~\ref{prop:cle3}, we get:

\begin{prop}\label{Dff} If $D=\phi(\frak f)$ then $\mathcal U_D=\mathcal U_{\frak
    f}$ and $\mathfrak f=\Vor_F(O)\cap F_D$. In particular, $\phi$ is injective.
\end{prop}

\begin{proof} That $\mathcal U_D$ coincides with $\mathcal U_{\frak f}$ is clear from
  their definitions and Proposition~\ref{prop:cle3}. The equality
  $\mathfrak f=\Vor_F(O)\cap F_D$ follows now because
  $\frak f$ is the intersection of the facets of $\Vor_F(O)$ containing
  it, and these facets are of the form $\Vor_F(O)\cap F_\beta$ for
  $\beta\in\Upsilon_1$ by Proposition~\ref{bondfacet}.
  \end{proof}

\subsubsection{Surjectivity of the map $\phi$}
In this section we prove the surjectivity of the map $\phi$, thus finishing the proof of Theorem~\ref{thm:vor}. 

Let $D$ be a coherent acyclic orientation of a cut subgraph $H$ of
$G$.

\begin{defi}[Codimension of a cut subgraph] \rm 
  The \emph{codimension} of a cut subgraph $H$ of $G$,
  denoted $\mathrm{codim}(H)$, is by definition the number
  of connected components of $G - E(H)$ minus one. The
  codimension of a coherent acyclic orientation $D$ of a cut subgraph 
  $H$ of $G$, denoted $\mathrm{codim}(D)$,
  is defined by $\mathrm{codim}(D) :=\mathrm{codim}(H)$. 
\end{defi}

The terminology is justified by the following result.

\begin{lemma}\label{lem:codimension} The following statements are
  true:
\begin{itemize}
\item $\mathrm{codim}(D) = |V|-1$ if and only if $D$ is an
  acyclic orientation of the whole $G$. 
\item $\E(D)=\bigcup_{\beta\in\mathcal U_D}\supp^+(\beta)$.
\item  The affine plane $\F_D$ is of codimension $\mathrm{codim}(D)$ in $\F_\R$. 
\end{itemize}
\end{lemma}

\begin{proof} If $\mathrm{codim}(D)=|V|-1$ then
$G-E(H)$ has $|V|$ connected components, which consists thus of single
vertices. It follows that $E(H)=E(G)$ and thus $D$ is an acyclic
orientation of the whole graph $G$. Conversely, if $D$ is an acyclic
orientation of the whole graph $G$ then $D$ is a coherent acyclic
orientation of the cut subgraph obtained by partitioning completely
$V(G)$, whence $\mathrm{codim}(D)=|V|-1$.

As for the second statement,
let $d := \mathrm{codim}(D)$, and denote by $C_1, \dots,
  C_{d+1}$ all the connected components of $G - E(D)$, ordered
  so that all the oriented edges of $D$ are from $C_i$ to $C_j$ for $i<j$. 

  Denote by $\beta_i$ the cut element of $\F_\Z$ associated to the cut 
$S_i = C_{i+1}\sqcup \dots \sqcup C_{d+1}$ for
  $i=1,\dots,d$. Clearly,
  $\E(D)=\supp^+(\beta_1)\cup\cdots\cup\supp^+(\beta_d)$. By 
Proposition~\ref{prop:sumofbonds}, we can find a decomposition of $\beta_i$ as a sum of bond 
elements with positive support included in $\E(D)$, i.e., as a sum of
elements in $\mathcal U_D$. Thus the second statement follows.

Finally, denote by $\widetilde G$ the multigraph obtained by contracting
  each set $C_i$ to a single vertex for $i=1,\dots, d+1$. The multigraph
$\widetilde G$ is connected and has $d+1$ vertices. We have an injective map from $\F_\Z(\widetilde G)$, the cut lattice of $\widetilde G$, 
to $\F_\Z$, the cut lattice  of $G$, and
$\beta_1,\dots, \beta_{d}$ come from $\F_\Z(\widetilde G)$. In addition, they form a basis of  $\F_\Z(\widetilde G)$. In particular, they are linearly independent.

Denote by $\chi_{\ind D} \in C^1(G, \R)$ the element which takes value 1 at each oriented edge of $D$, and value zero at each edge of $C_i$, for all $i=1,\dots, d+1$.
By definition, we have $x\in F_D$ if and only if 
\[
  \forall \beta \in \mathcal U_D, \quad 2\langle x,\beta\rangle
  = \|\beta\|^2 = \langle \chi_{\ind D}, \beta\rangle,
  \quad \textrm{or equivalently, if and only if}
\]
\[
  \forall \beta \in \mathcal U_D, \quad
  \langle 2x-\chi_{\ind D}, \beta \rangle =0.
\]
Since, as before, $\beta_i$ is a sum of elements in $\mathcal U_D$, we get 
\[
  \forall i=1, \dots, d, \quad
  \langle 2x-\chi_{\ind D}, \beta_i \rangle = 0.
\]
This proves that $F_D \subseteq \bigcap_{i=1}^d F_{\beta_i}$, where
$F_{\beta_i} := \bigl\{x\,\bigl|\, \langle x, \beta_i\rangle =
\|\beta_i\|^2\bigr\}$ for each $i$.
Since $\beta_1,\dots, \beta_{d}$ are linearly independent,
this proves that $\F_D$ is of codimension at least $d=\mathrm{codim}(D)$. 
On the other hand,  each cut element $\beta$ of $\F_\Z$
with positive support in $D$ is in the image of the inclusion map 
$\F_\Z(\widetilde G) \hookrightarrow \F_\Z$, and so 
is a linear combination of
$\beta_1, \dots, \beta_{d}$. This implies that 
\[
  \forall\,x\in \bigcap_{i=1}^d F_{\beta_i}, \,\,
  \forall \beta \in \F_\Z \textrm{ with }
  \supp^+(\beta) \subseteq \E(D), \qquad
  \langle 2x-\chi_{\ind D}, \beta\rangle =0,
\]
or equivalently $\langle 2x, \beta\rangle = \langle \chi_{\ind D}, \beta\rangle.$
Since for $\beta$ with $\supp^+(\beta) \subset \E(D)$, we have
$\langle \chi_{\ind D}, \beta\rangle = \|\beta\|^2$, it follows
that $\bigcap_{i=1}^d F_{\beta_i} 
\subseteq \F_D$.  We conclude that
$F_D = \bigcap_{i=1}^d F_{\beta_i}$, which proves the second
statement of the lemma.
\end{proof}

We will prove the existence of a face $\frak f$ of
$\Vor_F(O)$ of codimension $\mathrm{codim}(D)$ which is included
in $\F_D$, and with $\phi(\frak f) =D$, which proves the surjectivity
of $\phi$. We start by giving a  classification of all the vertices of $\Vor_\F(O)$.

\begin{lemma}\label{vert}
  If $D$ is an acyclic orientation of $G$, then $\F_D$
  consists of the single point $\nu^D$ which is the vertex
  of $\Vor_\F(O)$ satisfying $\phi(\nu^D)=D$. Furthermore,
  all the vertices of $\Vor_\F(O)$ are of this form.
\end{lemma}

\begin{proof}
  We already know that $\F_D$ has codimension $|V|-1$, i.e.,
  $\F_D = \{\nu^D\}$. To prove that $\nu^D$ is a vertex of
  $\Vor_F(O)$, we have to show that $\|\nu^D\| \leq
  \|\nu^D-\beta\|$ for every $\beta\in F_\Z$.
  Since $D$ is an acyclic orientation of $G$, the elements
  in $\mathcal U_D$ generate $\F_\Z$. Thus, there exist
  coefficients $a_\mu\in \Z$, for $\mu\in \mathcal U_D$, such that 
 $$\sum_{\mu \in \mathcal U_D} a_\mu \mu = \beta.$$
 We have 
 \begin{align*}
  2\langle \nu^D, \beta\rangle &= \sum_{\mu \in \mathcal U_D} 2a_\mu \langle \nu^D, \mu\rangle\\
  &=\sum_{\mu \in \mathcal U_D} a_\mu \|\mu\|^2 = \sum_{\mu \in \mathcal U_D} a_\mu \sum_{e\in \supp^+(\mu)} \mu(e)\\
  &=\sum_{e\in\E(D)} \beta(e) \leq \sum_{e\in\E(D)} |\beta(e)| = \|\beta\|^2,
 \end{align*}
 and the inequality is strict if and only if
 $\supp^+(\beta)$ is not contained in $D$.
 This shows that $\nu^D$ is a vertex of $\Vor_F(O)$.

Furthermore, since $\{\nu^D\}\subseteq F_\beta$ for each $\beta\in\mathcal
U_D$, it follows that $\mathcal U_D\subseteq\mathcal U_{\nu^D}$. Since 
$$
\E(D)=\bigcup_{\beta\in\mathcal U_D}\supp^+(\beta)
$$
by Lemma~\ref{lem:codimension}, and since $D$ is already an orientation of
the whole $G$, it follows that $\phi(\nu^D)=D$.

To prove all the vertices are of the form $\nu^D$, consider a
vertex $\nu$ and the coherent acyclic 
orientation $D = \phi(\nu)$ of a cut subgraph of $G$. 
It follows from Proposition~\ref{Dff} that $\mathcal
U_D=\mathcal U_{\nu}$ and $\{\nu\}=\Vor_F(O)\cap F_D$. Let $D'$ be an
acyclic orientation of the whole $G$ extending $D$. Clearly,
$F_{D'}\subseteq F_D$. But then
$$
\{\nu_{D'}\}=\Vor_F(O)\cap F_{D'}\subseteq \Vor_F(O)\cap F_D=\{\nu\}
$$
and $\nu$ and $\nu_{D'}$  must coincide.
\end{proof}

Consider now the general case, and suppose that $D_0$ is a 
coherent acyclic orientation of a cut subgraph of $G$. Consider 
all the acyclic orientations $D$ of $G$ which contain $D_0$, 
and for each such $D$, let $\nu^D$ be the corresponding vertex of 
$\Vor(O)$, which lies on the affine plane 
$\F_{D_0}$. Define $\frak f_0$ as the convex hull of all the vertices 
$\nu^D$ of $\Vor_F(O)$. 

\begin{prop} Notation as above, $\frak f_0$ is a face of 
$\Vor_F(O)$ of codimension equal to $\mathrm{codim}(D_0)$, 
and we have $\phi(\frak f_0) = D_0$. In particular, the map $\phi$ is surjective.
\end{prop}

\begin{proof} The proof is a matter of putting all we have proved 
so far together. The intersection of $\F_{D_0}$ with $\Vor_F(O)$ is 
the convex hull of the vertices of $\Vor_F(O)$ which lie on
$\F_{D_0}$. By Lemma~\ref{vert}, these vertices are of the form
$\nu^D$ for acyclic orientations $D$ of $G$. Furthermore, since
$\nu^D\subseteq F_{D_0}$, we have that $\mathcal
U_{D_0}\subseteq\mathcal U_{\nu^D}$. As $\mathcal U_{\nu^D}=\mathcal
U_{\phi(\nu^D)}$ by 
Proposition~\ref{Dff}, it follows from Lemma~\ref{vert} that $\mathcal
U_{D_0}\subseteq\mathcal U_{D}$,  whence $\E(D_0)\subseteq\E(D)$ 
by Lemma~\ref{lem:codimension}. This shows that $\frak f_0 = \F_{D_0} \cap \Vor_F(O)$. 

We claim now that $\mathcal U_{D_0}=\mathcal U_{\frak f_0}$. Indeed,
if $\beta\in\mathcal U_{D_0}$ then $F_\beta\supseteq
F_{D_0}\supseteq\frak f_0$, whence $\beta\in\mathcal U_{\frak
  f_0}$. Conversely, suppose $\beta\not\in \mathcal U_{D_0}$. 
Since $\E(D_0) = \bigcap_D \E(D)$ with the intersection 
running over all the acyclic orientations $D$ of $G$ 
extending $D_0$, we would infer the existence of one such orientation $D$ such that 
$\beta\not\in\mathcal U_D$. Since $D=\phi(\nu^D)$ by
Lemma~\ref{vert}, it follows from Proposition~\ref{Dff} that
$\beta\not\in\mathcal U_{\nu^D}$. Since $\nu^D\in\frak f_0$, we have
that $\beta\not\in\mathcal U_{\frak f_0}$.

To prove the assertion on the codimension of $\frak f_0$, by
Lemma~\ref{lem:codimension} we need prove that 
$\mathrm{codim}(\frak f_0) = \mathrm{codim}(F_{D_0})$. 
Otherwise though, we would have 
$\mathrm{codim}(\frak f_0) > \mathrm{codim}(F_{D_0})$,
which means there would be a bond element $\beta$ with 
$\frak f_0 \subset F_{\beta}$ and $\beta \notin \mathcal U_{D_0}$,
contradicting the claim just proved.

Finally, since $\mathcal U_{\frak f_0}=\mathcal U_{\phi(\frak f_0)}$
by Proposition~\ref{Dff},  the claim yields $\mathcal U_{D_0}=\mathcal
U_{\phi(\frak f_0)}$, and thus $D_0=\phi(\frak f_0)$ 
by Lemma~\ref{lem:codimension}.
\end{proof}

\subsection{The projection map $\square_{\F_\Z}$ to $\Vor_F$ revisited}

By Theorem~\ref{thm:projection1}, the union of cubes $\square_{F_\Z}$ 
projects onto the Voronoi decomposition of $\F_\R$ with respect to 
the lattice $\F_\Z$.  Let $\mathfrak f_D$ be the face of $\Vor_F(O)$ 
corresponding to a coherent acyclic orientation $D$ of a cut subgraph of $G$ 
of codimension $d$. To such $D\in \AC$ corresponds 
the face $\square_{D,0}$ of $\square_0$ of dimension $|E| - |E(D)|$ 
described as follows: 

First we have the point $\chi_{\ind D}\in C^1(G, \Z)$ associated to $D$ 
which takes value $+1$ on every $e\in \E(D)$, value $-1$ on every $e$ with 
 $\ol e\in\E(D)$, and value 0 elsewhere. 
Denote by $X_1, \dots, X_{d+1}$ all the connected 
components of $G - E(D)$, and $\square_{1,0}, \dots, \square_{d+1,0}$ the 
hypercubes corresponding to $X_1, \dots, X_{d+1}$, so we have 
$\square_{i,0}\subset C^{1}(X_i, \R)$ for $i=1, \dots, d+1$.

For each $i=1, \dots, d+1$, extending by zero, 
we have an inclusion $C^1(X_i, \R) \subset C^1(G, \R)$ which sends 
the lattice  $C^1(X_i, \Z)$ to $ C^1(G, \Z)$. In this way, we identify
$C^1(X_i, \R)$ with its image in $C^1(G,\R)$.  Define
$$
\square_{D,0}=\frac 12 \chi_{\ind D} + \sum_{i=1}^{d+1} \square_{i,0},
$$
where the sums are taken in $C^1(G, \R)$. Note that
$\square_{D,0}$ is a cube of dimension $|E| - |E(D)|$ 
consisting of the points $x \in C^1(G, \R)$ whose $e^{\mathrm{th}}$-coordinate $x_e$ 
is equal to $1/2$ for each $e\in \E(D)$, and  satisfies $-1/2\leq
x_e\leq 1/2$ for every $e$ 
with $e,\ol e\notin \E(D)$. 

For each $i=1,\dots,d+1$ 
denote by $F_{i,\R}$ the image of $C^0(X_i, \R)$ by the coboundary map 
$d_i: C^0(X_i, \R) \rightarrow C^1(X_i, \R)$ and by $F_{i,\Z}$ the cut lattice of $X_i$. 
Let $\Vor_{F_i}$ be the Voronoi decomposition 
of $F_{i,\R}$ with respect to $F_{i,\Z}$, and denote by $\Vor_{F_i}(O)$ the Voronoi cell of 
the origin in $\F_{i,\R}$. We view $\Vor_{F_i}(O)$ as a convex
polytope in $C^1(G, \R)$, under the inclusions $F_{i,\R}\subseteq
C^1(X_i,\R)\subseteq C^1(G,\R)$.

\begin{thm} \label{thm:projection1bis} Notations as above, consider
  the orthogonal projection map  $\pi_F: C^1(G, \R) \to \F_\R$. The
  following three statements hold:
\begin{enumerate}
 \item The cube $\square_{D,0}$ is mapped onto the face $\mathfrak f_D$ of $\Vor_F(O)$.
 \item The interior of $\square_{D,0}$ is mapped onto the interior of the face $\mathfrak f_D$.
 \item The map $\pi_F$ induces an isomorphism between 
$\Vor_{F_1}(O) + \dots + \Vor_{F_{d+1}}(O) + \frac 12 \chi_{\ind D} $ and the 
 face $\mathfrak f_D$ of $\Vor_F(O)$.
\end{enumerate}
\end{thm}

\begin{proof}
\noindent (1) For each point $x\in \square_{D,0}$ and each $\beta\in \mathcal U_D$, we have 
 \[
\langle \pi_F(x), \beta\rangle =\langle x,\beta\rangle = 
\frac 12\sum_{e\in \supp^+(\beta)} \beta_e = \frac 12 \|\beta\|^2,
\]
which shows that $\pi_F(\square_{D,0}) \subseteq F_D = 
\bigcap_{\beta \in \mathcal U_D}F_{\beta}$. 
Since $\pi_F( \square_{D,0})\subseteq\pi_F(\square_0)=\Vor_F(O)$, it follows that 
 $\square_{D,0}$ is projected into $\f_D$. On other hand, we claim
 that each $x\in \square_0$ with $\pi_F(x) \in \frak f_D$ 
verifies $x_e =\frac 12$ for all $e\in \E(D)$. Indeed, for each
$e\in\E(D)$, by Lemma~\ref{lem:codimension} there is $\beta\in\mathcal U_D$ such that
$\beta_e>0$, and thus $\beta_e=1$ by
Proposition~\ref{prop:generalizedcut}. Since 
$\frak f_D\subseteq F_\beta\cap\Vor_F(O)=\Vor_F(\beta)\cap\Vor_F(O)$, 
it follows from Theorem~\ref{thm:projection1} that $x\in\square_0\cap\square_\beta$, and thus
$x_e=1/2$. By the claim, the surjectivity of $\pi_F : \square_{D,0} \to \frak f_D$ 
follows from that of  $\pi_F\:\square_0\to\Vor_F(O)$.
 
 \medskip
 
 \noindent (2) We have to show that for all the points $x\in\square_{D, 0}$ 
with $|x_e|< \frac 12$ for every $e$ with $e,\ol e\notin \E(D)$, we have 
that $\pi_F(x)$ lies in the interior of $\f_D$. This is equivalent 
to showing that for each $D'\in \AC$ with $\E(D)\subsetneqq \E(D')$, 
we have $\pi_F(x) \notin \f_{D'}$. For such $D'$, by
Lemma~\ref{lem:codimension}, there exists a bond 
$\beta'\in \mathcal U_{D'}-\mathcal U_D$, for which we have
 \[
\langle \pi_F(x) , \beta'\rangle = \langle x , \beta'\rangle < \frac 12 \sum_{e\in
  \supp^+(\beta')} \beta'_e = \frac 12 \|\beta'\|^2,
\]
which shows that $\pi_F(x) \notin\f_{D'}$.

  \medskip
  
  \noindent (3) Denote  
  $\Delta_D := \Vor_{F_1}(O) + \dots + \Vor_{F_{d+1}}(O) + \frac 12 \chi_{\ind D}$. 
  We prove first that the map $\pi_F$ restricted to 
$\Delta_D$ is injective. Consider two  points 
  $x=x_1+\dots+x_{d+1} + \frac 12 \chi_{\ind D}$ and 
$y=y_1+ \dots + y_{d+1}+ \frac 12 \chi_{\ind D}$ in 
  $\Delta_D$ with $x_i,y_i \in \Vor_{F_i}(O)$ for 
$i=1,\dots, d+1$. Suppose $\pi_F(x) =\pi_F(y)$.
  We need to prove that $x_i =y_i$ for every $i$. 
For each $i$ and each $S \subset X_i$, consider the cut 
  $\beta:=d(\chi_{\ind S})$ defined by $S$ in $G$. 
From $\pi_F(x)=\pi_F(y)$ we get 
\[
0=\langle \pi_F(x)-\pi_F(y), \beta \rangle = 
\langle x,\beta \rangle - \langle y,\beta \rangle = 
  \langle x_i, \beta_i \rangle_i -\langle y_i, \beta_i\rangle_i, 
\]
  where $\langle.\,,\,.\rangle_i$ denotes the pairing in $C^1(X_i,\R)$
  and $\beta_i:=d_i(\chi_{\ind S})$. 
Since the cut elements
  $d_i(\chi_{\ind S})$ generate $F_{i,\R}$, and $x_i,y_i\in\
  F_{i,\R}$, it follows that $x_i = y_i$, hence the injectivity follows.

To prove the surjectivity, write 
  $z\in \square_{D, 0}$ as the sum $z_1+\dots+z_{d+1} + \frac 12
  \chi_{\ind D}$ for $z_i\in \square_{i,0}$ for 
$i=1,\dots, d+1$. Put $x_i:=\pi_{F_i}(z_i)$ for each $i$. Then $x_i\in\Vor_{F_i}(O)$ by 
Theorem~\ref{thm:projection1}, and
$x_i-z_i$ is in the cycle space of $X_i$. As this space is included in
the cycle space of $G$, we have $\pi_F(x_i)=\pi_F(z_i)$. Thus 
$\pi_F(z)=\pi_F(x)$ for $x=x_1+\dots+x_{d+1}+\frac 12 \chi_{\ind D}
\in \Delta_D$, 
and then surjectivity follows from (1). 
\end{proof}

 %%%%%%%%%%%%%%%%%Section 4

 \section{Tilings II: general edge lengths}\label{sec:generaltiling}
 
In this section, we generalize the results of Section~\ref{sec:tiling1} 
to the case of graphs in the presence of a general integer 
valued length function on edges.

Consider the setting of Section~\ref{sec:admissible}. 
Let $G=(V, E)$ be a graph with integer edge lengths 
$\ell: E \to \mathbb N$, and let $H$ be the subdivision of $G$ 
where each edge $e$ is subdivided $\ell_e -1$ times.
 Consider the Voronoi decomposition of $C^1(G, \R)$ with 
respect to the lattice $C^1(G,\Z)$, which is the tiling by 
hypercubes $\square_\alpha$ for $\alpha \in C^1(G, \Z)$, considered in Section~\ref{sec:tiling1}.

We will define a subcomplex $\square_{H}$ of this Voronoi 
tiling which takes into account the structure of the $G$-admissible 
divisors on $H$, and generalize Theorem~\ref{thm:projection1} 
 to this setting. As a result, we will obtain a tiling of 
$H_{0,\R} = \bigl\{f\in C^0(G,\R) \,\bigl|\, \sum_{v\in V} f(v)=0\bigr\}$ 
 into convex polytopes where each polytope in the tiling is 
congruent to the Voronoi polytope $\Vor_{F_{G'}}(O)$ of a subgraph $G'$ of $G$.  

We follow the notation we used in Section~\ref{sec:admissible}. In
particular, for each $f\in C^0(G, \Z)$, and each edge $e=uv\in \E(G)$, we set $\delta_{e}(f) := \lfloor \frac{f(v)-f(u)}{\ell_{e}}\rfloor$. 
  
\begin{defi}[Point $\dl_f$, subgraph $G_f$, cube $\square_f$ 
and arrangement of cubes $\square_{H}$]\rm 
For each $f\in C^0(G, \Z)$, define:
\begin{itemize}
\item the \emph{point} $\dl_f \in C^1(G, \frac12 \Z)$ by
\[
\forall \,e=uv\in \E(G), \quad \dl_f(e) :=
\frac12\Bigl({\delta_{e}(f)-\delta_{\ol e}(f)}\Bigr) = 
\begin{cases} \delta_{e}(f) & \quad \textrm{ if } \ell_{e} \,\bigl|\, f(v) - f(u),\\
\delta_{e}(f)+\frac 12& \quad \textrm{ otherwise.}
\end{cases}
\]
\item the \emph{spanning subgraph} $G_{\dl_f}$ of $G$ by setting
  $\E(G_{\dl_f})$
  to be the set of all $e=uv\in\E(G)$ with  $\ell_{e} \, \Bigl |\,
  f(v) -f(u)$, i.e., with $\dl_{f}(e) \in \Z$.
\item the \emph{cube} of dimension $|E(G_{\dl_f})|$ by 
\[
\square_{\dl_f}  := \Bigl\{ \dl_f + \epsilon \,\,\,\Bigl| \,\, 
\epsilon \in C^1(G, \R) \textrm{ with } 
\|\epsilon\|_{\infty} \leq\frac 12, \,\, \textrm{and}\,
\epsilon(e) = 0 \textrm{ for all } e\in \E(G) -\E(G_{\dl_f}) \Bigr\}.
\]
\item the \emph{arrangement of hypercubes} $\square_{G, \ell}$, 
simply denoted $\square_H$ if $G$ is understood from the context, by 
\[
\square_{H} = \square_{G, \ell} := 
\bigcup_{f\in C^0(G, \Z)}\,\,\square_{\dl_f}.
\]
\end{itemize}
\end{defi}

If there is no risk of confusion, we simply denote $G_{\dl_f}$ and $\square_{\dl_f}$
by $G_f$ and $\square_f$.
Notice that different $f$ might give rise to the same hypercube
$\square_f$. In addition,
for each $f\in C^0(G, \Z)$, we have
$$
\square_f  := \dl_f + \iota_f(\square_{G_f,0}),
$$
where $\square_{G_f,0}$ is the hypercube of the origin 
in the Voronoi decomposition of $C^1(G_f, \R)$ 
with respect to the sublattice $C^1(G_f, \Z)$, and 
$\iota_f\: C^1(G_f, \R) \hookrightarrow C^1(G,\R)$ is the map 
obtained by extending functions by zero.

Notice that the arrangement of hypercubes for 
uniform unitary edge lengths, $\ell = 1$, is identified 
with the arrangement $\square_{\F_\Z}$ associated to  
the cut lattice $F_\Z$ of $G$ from Section~\ref{sec:tiling1}; in other words, 
$\square_G = \square_{F_\Z}$.

\subsection{The projection map $\theta:\square_{G} \rightarrow \square_H$} Let $\square_{G} = \square_{\F_\Z} = \bigcup_{\beta\in \F_\Z} \square_\beta$ be as in Section~\ref{sec:tiling1}. In this section, we construct a projection
$\theta: \square_G \rightarrow \square_H$. 

\begin{defi}[projection $\square_{d(f)} \to \square_{\dl_f}$]\rm
  Let $f\in C^0(G, \Z)$. Define the map
  $\theta_f : \square_{d(f)} \rightarrow \square_{\dl_f}$ as follows: 
 for a point $d(f)+ \epsilon$, with $\|\epsilon \|_{\infty} \leq \frac 12$, define 
 \[
   \theta_f(d(f)+\epsilon) := \dl_f +
   \iota_f(\epsilon_{|C^1(G_f,\R)}),
 \]
where $\iota_f: C^1(G_f, \R) \hookrightarrow C^1(G,\R)$ is the
extension by zero. 
\end{defi}

Note that the definition is independent of the choice
of $f$ which gives the hypercube $\square_{d(f)}$ and
that different hypercubes in $\square_G$ might project
down to the same cube $\square_{\dl_f}$ (this
phenomenon is dependent on the edge length function $\ell$).

\begin{prop}\label{prop:projection} The collection of maps
  $\theta_f$ for $f\in C^0(G, \Z)$ are consistent on intersections
  of hypercubes,  thus yielding a well-defined projection
  map $\theta: \square_G \rightarrow  \square_H$. 
\end{prop}

\begin{proof} Let $\beta = d(f)$ and $\lambda = d(h)$ for
two elements $f, h \in C^0(G, \Z)$. 
We need to prove that the projection maps
$\theta_\beta$ and $\theta_\lambda$ are consistent on the
points of the intersection $\square_\beta \cap \square_\lambda$. 
We may suppose that $\square_\beta \cap \square_\lambda \neq\emptyset$.  Thus, $|\beta_e-\lambda_e|\leq 1$ for each $e\in\E(G)$.

Let $x$ be a  point in the intersection $\square_\beta\cap
\square_\lambda$. We show that $\theta_f(x)$ is equal to $\theta_h(x)$ by showing
that their $e$-th coordinates
are the same for each $e\in \E(G)$.
Consider first the case $\beta_e = \lambda_e$.
In this case, we can write
$x_e = \beta_e+ \epsilon_e = \lambda_e + \epsilon_e$
where $|\epsilon_e| \leq 1/2$. Two cases can happen: If
$\ell_e \mid \beta_e$ then
$\theta_f(x)_e = \theta_h(x)_e = \beta_e/\ell_e + \epsilon_e$;
whereas if $\ell_e \nmid \beta_e$ then
$\theta_f(x)_e = \theta_h(x)_e = \lfloor \beta_e/\ell_e  \rfloor +
1/2$. In either case the coordinates are the same. 

Consider now the case $|\beta_e-\lambda_e|=1$, and assume,
using symmetry, that $\beta_e =\lambda_e +1$. In this case,
we must have $x_e =\lambda_e + 1/2 = \beta_e-1/2$. We have 
$\theta_h(x)_e = \lfloor \lambda_e/\ell_e\rfloor + 1/2$.
Two cases can happen: If $\ell_e \mid \beta_e$ then $x_e = \beta_e
-1/2$, whence
$\theta_f (x)_e  = \beta_e/\ell_e - 1/2 =
\lfloor \lambda_e/\ell_e\rfloor + 1/2= \theta_h(x)_e$; whereas if
$\ell_e \nmid \beta_e$ then
$\theta_f(x)_e = \lfloor \beta_e /\ell_e\rfloor + 1/2
= \lfloor (\beta_e-1) /\ell_e\rfloor + 1/2 = \theta_h(x)_e$. Again, in
either case, the coordinates are the same.
\end{proof}

\subsection{Fibers of the projection map $\theta$}
Next, we describe the fibers of  $\theta$.
Consider $f \in C^1(G, \Z)$, and let $G_f$ be the subgraph of
$G$ associated to $f$. Let $C_1,\dots,C_k$ be the connected components
of $G_f$. Consider an element $h\in C^1(G, \Z)$. 
If $\dl_f = \dl_h$ then $G_f = G_h$ and $\square_{\dl_f} =
\square_{\dl_h}$; thus
$\square_{d(f)}$ and $\square_{d(h)}$ are mapped to the same 
cube $\square_{\dl_f} =\square_{\dl_h}$.

\begin{prop}\label{claim1}
Notation as above, if $\dl_f = \dl_h$ then the difference $h-f$ is
constant on each connected component $C_i$ of $G_f$.
\end{prop}

 \begin{proof}
For each oriented edge $e=uv\in \E(C_i)$, we have
$\ell_e \mid f(v) - f(u)$ and $\ell_e \mid h(v) -h(u)$, and from
$\dl_f=\dl_h$ we get
$$
\frac{f(v)-f(u)}{\ell_{e}} =\dl_f(e)=\dl_h(e)= \frac{h(v)-h(u)}{\ell_{e}},
$$
which by the connectedness of $C_i$ shows that $h-f$ is constant on $C_i$.
\end{proof}

\begin{cor} If $\dl_f = \dl_h$ and $G_f$ is connected then
$d(f)=d(h)$.
\end{cor}

\begin{claim}\label{claim2} Assume $\dl_f=\dl_h$. Denote by
  $\eta:\{1,\dots, k\} \rightarrow \mathbb Z$ the function which takes
  value $h(u) - f(u)$ at $i$ for every $u \in C_i$. 
Let $e=uv \in \E(G)$ be an oriented edge
  with $u\in C_i$ and $v\in C_j$ for distinct $i,j$. Then
  $$
  \lfloor\dl_f(e)\rfloor\ell_{e} - d(f)(e) <
  \eta(j) - \eta(i) <
  \lceil\dl_f(e)\rceil \ell_{e} - d(f)(e).
  $$
\end{claim}

\begin{proof}
Since $i\neq j$, we have
$\lfloor\dl_h(e)\rfloor \ell_{e} < h(v)- h(u) <
\lceil\dl_h(e)\rceil \ell_{e}$.
Substituting $\dl_f$ for $\dl_h$, and using that 
$h(v) =f(v) +\eta(j)$ and $h(u) = f(u) + \eta(i)$ 
gives the result.
\end{proof}

 Note that the interval appearing in the above claim has extremities $-r_e$ and $\ell_e -r_e$ with $r_e$ the rest of the division of $d(f)(e)$ by $\ell_e$,

\smallskip

The above claim leads us to the following definition.

\begin{defi}\label{Iij}\rm For all distinct $i,j$ for which there is an edge of
  $G$ between $C_i$ to $C_j$,  define $I_{ij}$ as  
\[
I_{ij} := \bigcap_{e\in \E(C_i,C_j)}
\,\, \Bigl[\, \lfloor\dl_f(e)\rfloor\ell_{e} - d(f)(e) + 1\,,
\,\lceil\dl_f(e)\rceil \ell_{e} - d(f)(e)-1 \,\Bigr].
\]
If there is no edge between $C_i$ and $C_j$, set $I_{ij}:=\R$.
\end{defi}

Note in particular that $0\in I_{ij}$ and  $I_{ij} = -I_{ji}$ for all distinct $i,j$. 
Using this terminology, we can now reformulate Claim~\ref{claim2} as
follows:

\begin{claim} Notations as above, assume $\dl_f=\dl_h$ and 
let $\eta : \{1,\dots, k\} \to \Z$ be the function taking the value at
$i$ equal to the difference of $h-f$ at any vertex in $C_i$. For all
distinct $i,j$ we have $\eta(j)-\eta(i) \in I_{ij}$. 
\end{claim}

We infer the following proposition:

\begin{prop}\label{hfhf}
 A function $h: V \rightarrow \Z$ verifies $\dl_h = \dl_f$ if and only
 if the difference $\eta=h-f$ verifies the assertions 
of Proposition~\ref{claim1} and Claim~\ref{claim2}, namely:
 \begin{itemize}
  \item The function $\eta$ is constant on each connected component
    $C_i$ of $G_f$ for 
$i=1,\dots,k$; denote by $\eta(i)$ the common value taken on $C_i$.
  \item For all distinct $i,j$ we have $\eta(j) - \eta(i) \in I_{ij}$.
 \end{itemize}

\end{prop}
\begin{proof} Immediate from the above discussion.
\end{proof}

So in order to understand the fibers of the projection map, it will be
crucial to understand the structure of all the functions $\eta$ 
which verify the two properties stated in the above proposition. 
Let us make the following useful definition:

Denote by $K_k$ the complete graph 
on $k$ vertices: the vertices are labeled $1,\dots, k$ and each
two distinct vertices are connected by a unique edge. 
Let $\mathcal I$ be a collection of closed intervals $I_{ij} \subseteq
\R$ associated to the oriented edges $ij$ of $K_k$ such that $0\in I_{ij}$ and
$I_{ij} = -I_{ji}$ for all distinct $i,j$. Assume in addition that if
$I_{ij}$ is compact, then its 
endpoints are both integer. We denote by $K_k[\mathcal I]$ the spanning subgraph of $K_k$ with edges consisting of all the pairs $\{i,j\}$ with $I_{ij}$ a compact interval. In our setting, $K_k[\mathcal I]$ coincides with the graph obtained from $G_f$ by contracting all the components $C_1, \dots, C_k$ of $G_f$, and then removing the loops and parallel edges (making the contracted graph simple).

\begin{defi}\label{defi:rfcom}\rm Notations as above, 
define $C^0(K_k[\I], \Z; \I)$ as the subset of all functions
$\eta\in C^0(K_k[\I],\Z)$ which verify $\eta(j) - \eta(i) \in I_{ij}$ for
all distinct $i,j$ with $\{i,j\}$ and edge of $K_k[\I]$, i.e., with $I_{ij}$ compact.
\end{defi}

Since $0\in I_{ij}$ for all $i,j$, 
we have $0 \in C^0(K_k[\mathcal I], \Z; \I)$. In particular, this set is non-empty. 

\smallskip

Consider the subspace $F_{K_k[\I],\R} := d(C^0(K_k[\I], \R)) \subset C^1(K_k[\I],\R)$, and define the convex polyhedron $P_\I \subset C^1(K_k[\I], \R)$ as 
\[
P_\I := \Bigl\{x \in F_{K_k[\I],\R} \,\Bigl |\,  x_{ij} \in I_{ij} 
\textrm{ for all oriented edges  } ij\in \E(K_k[\I])\Bigr\}.
\]
By definition, $C^0(K_k[\I], \Z; \mathcal I)$ is the set of all 
$\eta \in C^0(K_k[\I],\Z)$ such that $d\eta \in F_{K_k[\I],\Z} \cap P_\I$.  
It follows that $C^0(K_k[\I], \Z; \mathcal I)$ is the set of integer 
points of  $d^{-1} (P_\I) \subseteq C^0(K_k[\mathcal I], \R)$, which is a convex polyhedron. 

Let $P_{\I, \Z} := P_\I \cap F_{K_k[\I],\Z}$, so we have $P_{\I, \Z} = d(C^0(K_k[\I], \Z; \mathcal I))$. From the above discussion we get the following proposition:

\begin{prop}\label{prop:fiber1} For each point $p$ in the relative 
interior of the cube $\square_{\dl_f}$, we have 
\[
\theta^{-1}(p) \simeq \bigcup_{\mu \in P_{\I,\Z}} \square_\mu,
\] 
where $\I$ is given in Definition~\ref{Iij}.
\end{prop}

To be more precise, let $\jmath: C^1(K_k[\I], \R) \hookrightarrow C^1(G,
\R)$ be defined as follows: For each $x\in C^1(K_k[\I], \R)$, the image 
$\jmath(x)$ in $C^1(G, \R)$ takes value $x_{ij}$ on each oriented edge $e=uv$ with 
$u\in C_i$ and $v \in C_j$ for all distinct $i,j$ with $\{i,j\}$ and edge of $K_k[\I]$, and takes value
zero elsewhere. Then, if $p=\dl_f+y$,
\[
\theta^{-1}(p) = 
d(f) + y +
\jmath\Bigl(\bigcup_{\mu \in P_{\I,\Z}} \square_\mu\Bigr) = 
d(f) + y +
\bigcup_{\mu \in P_{\I,\Z}} \jmath\bigl(\square_\mu\bigr) \simeq \bigcup_{\mu \in P_{\I,\Z}} \square_\mu.\]

 \subsection{The projection map 
$d^*:\square_H \rightarrow  H_{0,\R}$ 
and the mixed Voronoi tiling}
 
Consider the map $d^*: \square_H \rightarrow C^{0}(G, \R)$, 
which has image in $H_{0, \R}$. Recall that $H_{0,\R}$ is the 
hyperplane of $C^0(G, \R)$ which consists of all the functions  $f\in
C^0(G,\R)$ with 
$\sum_{v\in V(G)} f(v) =0$.

\smallskip

 In this section, we prove a generalization of Theorem~\ref{thm:projection1}. 
 
\begin{prop} \label{prop:proj2-1}
 For each $f \in C^0(G, \Z)$ such that $G_f$ is connected, 
the image $d^*(\square_f)$ is a polytope of dimension $|V|-1$ 
in $H_{0,\R}\subset C^0(G, \R)$ which is congruent to
$\Vor_{G_f}(O)$. More precisely, $d^*(\square_f) = d^*(\dl_f) + \Vor_{G_f}(O)$.
Furthermore, $d^*$ maps the relative interior of $\square_f$ onto the
interior of $d^*(\square_f) $.
\end{prop}

(Notice that $C^0(G_f,\R)=C^0(G,\R)$ naturally, so $\Vor_{G_f}(O)$ can
be seen in $C^0(G,\R)$.)

\begin{proof} 
Let $f\in C^0(G, \Z)$ with $G_f$ connected. We have 
$$
d^*(\square_f) = d^*(\dl_f) + d^*(\square_{G_f,0}),
$$
where $\square_{G_f,0}$ denotes the hypercube with center 
$0$ in $C^1(G_f, \R)$, consisting of all the points $\epsilon$ 
with $\|\epsilon \|_\infty \leq \frac 12$, viewed in $C^1(G,\R)$ by
extending by zero. The adjoint $d^*$ for $G$ restricts to that for
$G_f$ for this natural inclusion of $C^1(G_f,\R)$ in
$C^1(G,\R)$. Thus, by Theorem~\ref{thm:projection1}, since $G_f$ is connected, we infer that
$d^*(\square_{G_f,0})$ is the Voronoi cell of the origin 
$\Vor_{G_f}(O)$ in $C^0(G_f, \R) = C^0(G, \R)$. Furthermore, by the
same theorem, the interior of $\square_{G_f,0}$ is mapped onto the
interior of $\Vor_{G_f}(O)$.
\end{proof}

\begin{defi}\rm
 For each $f \in C^0(G, \Z)$ with connected $G_f$, we call
 $d^*(\square_f)$ the Voronoi polytope associated to $(G, \ell)$ and
 $f$, and denote it by $\Vor_{G,\ell}(f)$ or simply $\Vor_H(f)$
if $G$ is understood from the context.
\end{defi}

\begin{prop}\label{prop:proj2-2}
Let $f$ and $h$ be two elements of $C^0(G, \Z)$ such that 
$f- h$ is not constant and such that $G_f$ and $G_h$ are both
connected. Then the interiors of $\Vor_H(f)$ and $\Vor_H(h)$ are disjoint.
\end{prop}

\begin{proof}
We first note that $\dl_f\neq \dl_h$. Indeed, otherwise, 
if $\dl_f=\dl_h$, since $G_f$ is connected, Proposition~\ref{claim1} would
yield that $f-h$ would be constant.

For the sake of a contradiction, using Proposition~\ref{prop:proj2-1}, 
suppose there exist $\dl_f+ \epsilon$ and $\dl_h +
\varepsilon$ in the relative interiors of $\square_{f}$ and
$\square_h$, respectively, 
satisfying $d^*(\dl_f+\epsilon) = d^*(\dl_h + \varepsilon)$. By definition,
for each $e\in\E(G)$ we have:
\begin{itemize}
\item If $e\notin \E(G_f)$ then $\epsilon_e =0$, whereas if $e\in
  \E(G_f)$ then $|\epsilon_e|< \frac 12$.
\item If $e\notin \E(G_h)$ then $\varepsilon_e =0$, whereas if $e\in
  \E(G_h)$ then $|\varepsilon_e|< \frac 12$.
\end{itemize}
Rewriting the assumption, we have
\begin{equation}\label{eq4}
 d^*(\dl_f - \dl_h) = d^*(\varepsilon - \epsilon).
\end{equation}

Let $\eta := f-h$, and let $S$ be the set of all vertices of $G$ where $\eta$ takes its maximum value. Since $\eta$ is not constant, we have $S\subsetneq V(G)$.
Since $G_f$ and $G_h$ are both connected spanning subgraphs of $G$, 
we have $\E(G_f) \cap \E(S, V-S) \neq \emptyset$ and 
$\E(G_h) \cap \E(S, V-S) \neq \emptyset$. 

\medskip

\begin{claim}
 For each $e\in \E(S, V-S)$, we have 
$\dl_f(e) - \dl_h (e) \leq \varepsilon_e - \epsilon_e$. 
In addition, the inequality is  strict if $e$ belongs to either $G_f$ or $G_h$.
\end{claim}

Once this claim has been proved, we get a contradiction to
\eqref{eq4}. Indeed, we would have
\[
\sum_{v\in V-S} d^*(\dl_f - \dl_h)(v) = 
\sum_{e\in \E(S, V-S)} \Bigl(\dl_f(e) - \dl_h (e) \Bigr) 
< \sum_{e\in \E(S, V-S)} \bigl(\varepsilon_e - \epsilon_e\bigr) 
= \sum_{v\in V-S} d^*(\varepsilon - \epsilon)(v).
\]

\medskip
The proof of the claim is based on the following case by case analysis:

\begin{itemize}
\item If $e=uv \in \E(S, V-S) \cap \E(G_f)\cap \E(G_h)$ then both $f(v)-f(u)$ and $h(v) -h(u)$ are divisible by $\ell_e$, and since $\eta(v)< \eta(u)$, we get 
$$
\dl_f(e) - \dl_h (e) = 
\frac 1{\ell_e}\Bigl(f(v)-f(u)\Bigr)  
-\frac1{\ell_e}\Bigl(h(v)-h(u)\Bigr) =
\frac 1{\ell_e}\Bigl(\eta(v)-\eta(u)\Bigr)  
\leq -1\,.
$$
On the other hand, $- 1<\varepsilon_e -\epsilon_e < 1$, which proves the claim in this case.
\item Similarly, if $e=uv \in \E(S, V-S) \cap \E(G_f)$ but
  $e\not\in\E(G_h),$ then  
$$
\dl_f(e) - \dl_h (e) = \frac 1{\ell_e}\Bigl(h(v)-h(u) +
\eta(v)-\eta(u)\Bigr) - \Bigl\lfloor
\frac{h(v)-h(u)}{\ell_e}\Bigr\rfloor - 1/2 \leq -1/2,
$$
while $-\frac 12 < \varepsilon_e -\epsilon_e < \frac 12$.
\item If $e\in \E(S, V-S) \cap \E(G_h)$ but $e\not\in\E(G_f)$, the proof is similar.

\item Finally, if $e \in \E(S,V-S)$ but $e$ is neither in $G_f$ nor in
  $G_h$, then $\dl_f(e) - \dl_h(e) \leq 0 $, while $\varepsilon(e) -\epsilon(e)=0$.
\end{itemize}
This finishes the proof of the claim, and then that of the proposition.
\end{proof}

We can now state the main theorem of this section, 
which is the promised generalization of Theorem~\ref{thm:projection1}.

 \begin{thm}\label{thm:projection2} The $\Vor_H(f)$ with
   $G_f$ connected are the top dimensional cells of a tiling
   $\Vor_{G,\ell}$ of $H_{0,\R}$ into polytopes. Each top dimensional
   cell is congruent to the Voronoi 
 cell of a connected subgraph of $G$. More precisely, $\Vor_H(f)$ is
 congruent to $\Vor_{G_f}(O)$ for each $f\in C^0(G,\Z)$ with $G_f$ connected.
 \end{thm}

 \begin{defi}\rm
  The subdivision of $H_{0,\R}$ into Voronoi cell associated to subgraphs $G_f$ is referred to as the \emph{mixed Voronoi tiling} of $H_{0,\R}$ induced by the edge length function $\ell$ on $G$.
 \end{defi}
 \begin{remark}\rm We note that the connected subgraphs $G_f$
   determining the shape of polytopes that do appear in the tiling
   are special and dependent on the edge length function $\ell$.
   For example, they  all contain all the edges $e$ of $G$ for which $\ell_e =1$.
   Changing the length function can drastically change the tiling.
   On the one hand, only finitely many polytopes, up to translation,
   are used in all these tilings, and so the congruence classes of
   polytopes which appear form a finite set. On the other hand,
   there exist tilings which use all these polytopes. An example of
   such a tiling is the one associated with an edge length function
   $\ell$ which verifies that $\ell_e>1$ for all $e$ and that $\ell_e$
   and $\ell_{e'}$ are coprime for all pairs of distinct edges $e$ and
   $e'$. In this case, an application of the Chinese Remainder Theorem
   shows that every connected spanning subgraph $G'$ of $G$ is equal
   to $G_f$ for an element $f\in C^0(G, \Z)$, and so the Voronoi
   polytopes associated to all connected spanning subgraphs of $G$
   appear in the tiling.      
 \end{remark}
 
 \begin{remark}\rm
 Denote by $\mathrm{N}$ the least common multiple of all the $\ell_e$ for $e$ an edge of $G$. The tiling is periodic with respect to the sublattice $d^*(\dl(\mathrm{N}\,C^0(G,\Z )))$.  
\end{remark}

 We need some preparation before giving the
 proof of Theorem~\ref{thm:projection2}. Let $f\in C^{0}(G,\Z)$
 be an element with connected $G_f$. Let $S \subseteq V$ be a proper subset such
 that the characteristic function $\beta_f$ of the oriented cut
 $\E_{G_f}(V-S, S)$ is a bond, that is, such that
 $G_f[S]$ and $G_f[V-S]$ connected. Consider the functions
 $h_n:=f+ n\chi_{\ind S}$ for $n\in \N$. Note that
 $G_{h_n}[S] = G_f[S]$ and $G_{h_n}[V-S] = G_f[V-S]$, whence
 $G_f$ and $G_{h_n}$ can only differ in the edges in $G$
 which are in the cut $E(S,V-S)$. Let $h:=h_n$ for
 the smallest posititive integer $n$ such that 
 $G_{h_n}$ is connected. Since $G_f[S]$ and $G_f[V-S]$ are
 connected, $n$ is also the smallest positive integer for which there
 is $e\in\E(S,V-S)$ with $\dl_{h_n}(e)\in\Z$. In particular,
 $\dl_f(e)-1\leq\dl_h(e)\leq\dl_f(e)$ for each $e\in\E(S,V-S)$,
 the upper bound being achieved only if $\dl_f(e)\not\in\Z$ and the
 lower bound being achieved only if $\dl_f(e)\in\Z$. 

 Also because $G_f[S]$ and $G_f[V-S]$ are
 connected, the characteristic function $\beta_h$ of the oriented cut
 $\E_{G_h}(V-S, S)$ is also a bond. Let $\eta_f$ and $\eta_h$ be the
 extensions by zero of $\beta_f$ and $\beta_h$ to $C^1(G,\Z)$, and define
 \begin{equation}\label{etafh}
   \eta_{S,f}:= \frac 12 \eta_f + \frac 12 \eta_h\,\in\, C^1(G,\frac{1}{2}\Z).
\end{equation}

 \begin{lemma}\label{lem:sharedfacet}
  Notations as above, we have $\dl_h -\dl_f = \eta_{S,f}$.  In addition, the Voronoi cells $\Vor_{H}(f)$ and $\Vor_H(h)$
  intersect at the facet of $\Vor_H(f)$ (resp.~$\Vor_H(h)$) given by the bond
  $\beta_f$ of $G_f$ (resp.~by the bond $-\beta_h$ of $G_h$).
\end{lemma}

 \begin{proof} The first assertion is obtained from a case analysis:
   Let $e\in\E(G)$. If $e\in \E(S)$ or $e\in\E(V-S)$ then
   $\dl_h(e) =\dl_f(e) = \dl_f(e) + \eta_{S,f}(e)$, since $\eta_{S,f}(e)=0$. The
   same holds if $e\in\E(S,V-S)$ but $e\not\in\E(G_f)\cup\E(G_h)$. 
  Suppose now $e\in \E(S, V-S)$. If $e\in\E(G_f)\cap\E(G_h)$ then
   $\dl_h(e) = \dl_f(e)-1 = \dl_f(e)+ \eta_{S,f}(e)$, as $\eta_f(e) =
   \eta_h(e)= -1$. If $e\in\E(G_f)$ but $e\not\in\E(G_h)$ then 
$\dl_h(e) = \dl_f(e)-1/2 = \dl_f(e) + \eta_{S,f}(e)$, since $\eta_f(e)=-1$
but $\eta_h(e) =0$. 
The same holds if $e\in\E(G_h)$ but $e\not\in\E(G_f)$, the only
difference being that now $\eta_f(e)=0$
but $\eta_h(e) =-1$. The first assertion is proved.
 
\smallskip
   
As for the second assertion, by Theorem~\ref{thm:projection1bis} applied to $G_f$,
the facet of $\Vor_{G_f}(O)$ given by the bond $\beta_f$ is the image
under $d^*$ of 
$\frac 12 \eta_f + \square_{G_f[S],0} + \square_{G_f[V-S],0}$, where 
$\square_{G_f[S],0}$ and $\square_{G_f[V-S],0}$ are viewed as subsets
of $\square_0$. Similarly, by Theorem~\ref{thm:projection1bis}
applied to $G_h$, the facet of $\Vor_{G_h}(O)$ given by the bond
$-\beta_h$ is the image under $d^*$ of
$-\frac 12 \eta_h + \square_{G_h[S],0} +\square_{G_h[V-S],0}$. 
Since $G_{h}[S] = G_f[S]$ and $G_{h}[V-S] = G_f[V-S]$, it follows from
the first assertion that
\begin{equation}\label{facetfh}
  \dl_f + \frac 12 \eta_f + \square_{G_f[S],0} +\square_{G_f[V-S],0} =
  \dl_h - \frac 12 \eta_h + \square_{G_h[S],0} + \square_{G_h[V-S],0}.
\end{equation}
Since  $\Vor_H(f) = d^*(\dl_f)+ \Vor_{G_f}(O)$, it follows
that the facet of $\Vor_H(f)$ given by $\beta_f$ is the
image under $d^*$ of the left hand side of
\eqref{facetfh}. Analogously, the facet of $\Vor_H(h)$ given by
$-\beta_h$ is the image under $d^*$ of the right hand side of
\eqref{facetfh}. So they coincide. The lemma follows since the two Voronoi cells
$\Vor_H(f)$ and $\Vor_H(h)$ have disjoint interiors by
Proposition~\ref{prop:proj2-2}. 
 \end{proof}

We can now present the proof of Theorem~\ref{thm:projection2}.
 
 \begin{proof}[Proof of Theorem~\ref{thm:projection2}]
   We proved in Proposition~\ref{prop:proj2-2} that the
   polytopes $\Vor_H(f)$ for $f\in C^0(G, \Z)$ with $G_f$
   connected have disjoint interiors. In order to prove the
   theorem, it will be thus enough to show that the union 
   of the $\Vor_H(f)$ with $G_f$ connected covers $H_{0,\R}$,
   or equivalently, that the projection map
   $d^*: \square_H \rightarrow H_{0,\R}$ is surjective.

Now, each tile $\Vor_H(f)$ is a translation of the Voronoi cell
$\Vor_{G_f}(O)$ associated to a connected subgraph $G_f$ of $G$. As
there are finitely many such subgraphs, it follows from
Proposition~\ref{prop:proj2-2} that the union 
$$
\mathcal H = \bigcup_{\substack{f\in C^0(G, \Z)\\ G_f \textrm{
      connected }}} \Vor_{H}(f)
$$
is a closed subset of  $H_{0,\R}$.
Now consider a face $\f$ of a polytope $\Vor_{H}(f)$ in $\mathcal H$.
It will be enough to show that the star of $\f$ in the arrangement of the polytopes 
$\mathcal H$ is a complete fan. Suppose for the sake of a contradiction that this is not the case.
There exists then an element $f' \in C^0(G, \Z)$ such that $G_{f'}$ 
is connected, $\f$ belongs to $\Vor_{H}(f')$, and $\Vor_{H}(f')$ contains a facet which is not shared by any other 
polytope among the polytopes in the arrangement $\mathcal H$.  But
this is a contradiction. Indeed,  any facet of $\Vor_H(f') = \dl_{f'}
+ \Vor_{G_f}(O)$  corresponds to a bond in $G_f$ by
Theorem~\ref{thm:projection1}, and in Lemma~\ref{lem:sharedfacet} we
proved that such a facet  is shared with another polytope $\Vor_H(h)$
with $h\in C^0(G, \Z)$ such that $G_h$ is connected. 
 \end{proof}

For each $f: C^0(G, \Z)$ with $G_f$ connected and each proper subset
$S\subset V$ such that $G_f[S]$ and $G_f[V-S]$ are connected, recall 
$\eta_{S, f}\in C^1(G,\frac 12\Z)$, defined in \eqref{etafh}.

 \begin{thm}
  The one-skeleton of the dual of the mixed Voronoi tiling is
  isomorphic to the graph with vertex set $\dl_f$ for all 
  $f\in C^0(G, \Z)$ with $G_f$ connected, and with edge set
  $\{\dl_f, \dl_f+\eta_{S,f}\}$ for all $f$ as above and all proper subsets $S\subset V$
  such that $G_f[S]$ and $G_f[V-S]$ are connected.
 \end{thm}

 \begin{proof}
   This is a local assertion for each facet of $\Vor_{H}(f)$, and thus
   follows from Theorem~\ref{thm:projection2} and Lemma~\ref{lem:sharedfacet}.
 \end{proof}

\section{Tilings III: general edge lengths with twisting}\label{sec:moregeneraltiling}

In this section, we generalize the results of the
previous section to the case of general edge lengths
with a twisting. This will be connected to the admissible
divisors we considered in Section~\ref{sec:admissible}.

Let $G=(V, E)$ be a graph with integer edge
lengths $\ell: E \to \mathbb N$, and let $H$
be the subdivision of $G$ where each edge $e$
is subdivided $\ell_e -1$ times. 
Let $\mathfrak m\in C^1(G, \Z)$. We use
$\m_e$ to denote the value of $\m$ at the oriented edge $e\in \E(G)$.

As in the previous section, we will consider the
Voronoi decomposition of $C^1(G, \R)$ with respect to the
lattice $C^1(G,\Z)$, which is the tiling by hypercubes
$\square_\alpha$ for $\alpha \in C^1(G, \Z)$,
define a subcomplex $\square_{H}^\m$ of this tiling as the union of
certain specific cubes, and then take the projection by 
$d^*$ to obtain a tiling of
$H_{0,\R} = \bigl\{f\in C^0(G,\R) \,\bigl|\, \sum_{v\in V} f(v)=0\bigr\}$ 
in convex polytopes in which each polytope is congruent to the
Voronoi polytope $\Vor_{G'}(O)$ associated to a certain
connected spanning subgraph $G'$ of $G$. 

The subcomplex $\square_{H}^\m$  is defined as follows:
First, for each $f\in C^0(G, \Z)$ and each edge $e=uv\in \E(G)$, set
$\delta^\m_e(f) := \lfloor \frac{f(v)-f(u)+\m_e}{\ell_e}\rfloor$.  
 
\begin{defi}[Point $\dl^\m_f$, subgraph $G^\m_f$, cube $\square^\m_f$
  and arrangement of cubes $\square^\m_{H}$]\rm For each $f\in C^0(G, \Z)$, define:
\begin{itemize}
\item the point $\dl^\m_f \in C^1(G, \frac12 \Z)$ by
\[
\forall \,e=uv\in \E,  \,\,\dl^\m_f(e) :=  
\frac12\Bigl({\delta^\m_e(f)-\delta_e^\m(f)}\Bigr) = 
\begin{cases} 
\delta^\m_e(f) &  \textrm{ if } \ell_e \,\bigl|\, f(v) - f(u)+\m_e\\
\delta^\m_e(f)+\frac 12&  \textrm{ otherwise;}
\end{cases}
\]
\item the spanning subgraph $G_{\dl^\m_f}$ of $G$ by setting 
$\E(G_{\dl^\m_f})$ to be the set of all $e=uv\in\E(G)$ with  
$\ell_e \, \Bigl |\, f(v) -f(u)+\m_e$, i.e., with $\dl^\m_{f}(e) \in \Z$;
\item the cube of dimension $|E(G_{\dl^\m_f})|$ by 
\[
\square_{\dl^\m_f}  := \Bigl\{ \dl^\m_f + \epsilon \,\Bigl| \,\, 
\epsilon \in C^1(G, \R) \textrm{ with } \|\epsilon\|_{\infty} \leq
\frac 12, \,\, \textrm{and}\,\epsilon(e) = 0 \textrm{ for } 
e\in \E(G) - \E(G_{\dl^\m_f}) \Bigr\};\]
\item the arrangement of hypercubes $\square_{G, \ell}^\m$, simply
denoted $\square_H^\m$ if $G$ is understood from the context, by 
\[
\square^\m_{H} = \square^\m_{G, \ell} := 
\bigcup_{f\in C^0(G, \Z)}\,\,\square_{\dl^\m_f}.
\]
\end{itemize}
\end{defi}

If there is no risk of confusion, we simply denote $G_{\dl^\m_f}$ and
$\square_{\dl^\m_f}$ by $G^\m_f$ and $\square_f^\m$.

For each $f\in C^0(G, \Z)$ we denote by $\iota_f^\m$ 
the inclusion $C^1(G_f^\m, \R) \hookrightarrow C^1(G,\R)$ 
obtained by extending by zero an element $\alpha\in C^1(G^\m_f, \R)$ 
to the whole graph $G$. For each $f\in C^0(G, \Z)$, we have
$\square^\m_f  := \dl^\m_f + \iota_f^\m(\square_{G^\m_f,0})$, 
where $\square_{G^\m_f,0}$ is the hypercube of the origin 
in the Voronoi decomposition of $C^1(G^\m_f, \R)$ with 
respect to the sublattice $C^1(G^\m_f, \Z)$.

We note that the arrangement of hypercubes with $\m=0$ is identified with 
the arrangement $\square_{H}$ of the previous section. For $\ell=1$ and $\m=0$, we obtain 
$\square^\m_{H} = \square_G$.

Proceeding as in the previous section, 
we construct a projection
$\theta^\m: \square_G \rightarrow \square_H^\m$.

\begin{defi}[Projection $\square_{d(f)} \to \square_{\dl^\m_f}$]\rm 
Let $f\in C^0(G, \Z)$. Denote by $\theta_f^\m : \square_{d(f)}
\rightarrow \square_{\dl^\m_f}$ 
the map taking each point $d(f)+ \epsilon$ with 
$\|\epsilon \|_{\infty} \leq\frac 12$ to $\dl^\m_f +\iota^\m_f(\epsilon_{|G^\m_f})$.
\end{defi}

\begin{prop} \label{prop:projection3} There is a well-defined
  projection map 
$\theta^\m: \square_G \rightarrow  \square^\m_H$ restricting to the 
$\theta^\m_f$ for all $f\in C^0(G, \Z)$.
\end{prop}

\begin{proof} 
The proof is similar to that of Proposition~\ref{prop:projection} and is omitted. 
\end{proof}

Let $f \in C^1(G, \Z)$, and denote by $C_1, \dots, C_k$ 
the connected components of 
 $G^\m_f$.  Let $h\in C^1(G, \Z)$. We have the following extension of Proposition~\ref{claim1}.  

\begin{prop}\label{claim1bis} Notations as above, if $\dl^\m_h= \dl^\m_f$ 
then the difference $h-f$ is constant on each connected component $C_i$ of $G^\m_f$.
\end{prop}

 \begin{proof} The proof is similar to that of Proposition~\ref{claim1}. 
Namely, from $\dl^\m_f = \dl^\m_h$, we get $G^\m_f =G^\m_h$, which
gives for each $i$ and $e=uv\in \E(C_i)$: 
$$
\frac{f(v)-f(u)+\m_e}{\ell_e} = \dl^\m_f(e)=\dl_h^\m(e)=
\frac{h(v)-h(u) +\m_e}{\ell_e}.
$$ 
From the connectedness of $C_i$, we conclude that $h-f$ is constant on $C_i$.
\end{proof}

This then leads to the definition of $\eta:\{1,\dots, k\} \rightarrow
\mathbb Z$, the function which takes value $h(u) - f(u)$ at $i$ for
any $u \in C_i$, and the intervals  $I_{ij}$ for all distinct $i,j\in\{1,\dots,k\}$:
\[
I_{ij} :=  \bigcap_{e\in \E(C_i,C_j)} \,\, 
\Bigl[\, \lfloor\dl^\m_f(e)\rfloor\ell_e - d(f)(e) + 1\,, 
\,\lceil\dl_f^\m(e)\rceil \ell_e - d(f)(e)-1 \,\Bigr] 
\]
if $\E(C_i,C_j)$ is non-empty, and $I_{ij}:=\R$ otherwise. We 
get the following extension of Proposition~\ref{hfhf}:

\begin{prop}
 A function $h: V \rightarrow \Z$ verifies $\dl^\m_h = \dl^\m_f$ if
 and only if the difference $\eta=h-f$ verifies the following two properties:
 \begin{itemize}
  \item The function $\eta$ is constant on each connected component
    $C_i$ of $G^\m_f$; denote by $\eta(i)$ the common value taken on $C_i$.
  \item For all distinct $i,j$ we have $\eta(j) - \eta(i) \in I_{ij}$.
 \end{itemize}
\end{prop}

This leads verbatim as before to the characterization of the 
fibers of the projection map $\theta^\m$, as follows. Recall that for a collection of intervals $\I$ as above, 
we denote by $K_k[\I]$ the subgraph  of $K_k$ with edges $\{i,j\}$ such that $I_{ij}$ is compact, that 
$P_{\I}$ is the polytope of all the  $x\in F_{K_k[\I],\R}$ such that
$x_{ij}\in I_{ij}$ for every  oriented edge $ij$ $\in\E(K_k[I])$, and that
$P_{\I,\Z}:=\P_{\I}\cap F_{K_k,\Z}$ is the set of integer points of $P_\I$.

Let, as before, $\jmath: C^1(K_k[\I], \R) \hookrightarrow C^1(G, \R)$ be the 
embedding which sends the point $x\in C^1(K_k[\I], \R)$ to the element of
$C^1(G, \R)$ which takes value $x_{ij}$ on each $e=uv$ with $u\in C_i$
and $v \in C_j$ for distinct $i,j$ with $I_{ij}$ compact, and takes value zero elsewhere.

\begin{prop} For each point $p=\dl^\m_f+x$ in the relative interior of
  the cube $\square_{\dl^\m_f}$, we have 
\[
\theta^{\m^{-1}}(p) = d(f) + x+
\jmath\Bigl(\bigcup_{\mu \in P_{\I,\Z}} \square_\mu\Bigr) 
= d(f) + x+
\bigcup_{\mu \in P_{\I,\Z}} \jmath\bigl(\square_\mu\bigr) 
\simeq \bigcup_{\mu \in P_{\I,\Z}} \square_\mu.
\] 
\end{prop}

Consider now the map $d^*: \square_H^\m \rightarrow H_{0,\R}$. We have the following 
direct extensions of Propositions~\ref{prop:proj2-1} and~\ref{prop:proj2-2}:
 
\begin{prop} \label{prop:proj3-1} For each $f \in C^0(G, \R)$ with $G^\m_f$ connected, we have 
 $d^*(\square^\m_f) = d^*(\dl^\m_f) + \Vor_{G^\m_f}(O)$. In
 particular, $d^*(\square^\m_f)$ is a polytope of top dimension
 $|V|-1$ in $H_{0,\R}$. Furthermore,
 $d^*$ maps the relative interior of $\square^\m_f$ onto the interior
 of $d^*(\square^\m_f)$.
 \end{prop}

 \begin{proof}
  Similar to the proof of Proposition~\ref{prop:proj2-1}.
 \end{proof}

For each $f \in C^0(G, \R)$ with connected $G^\m_f$, define 
$\Vor_{H}^\m(f):=d^*(\square^\m_f)$.
 
\begin{prop} \label{prop:proj3-2} For each $f$ and $h$ in $C^0(G, \Z)$ such that 
$f- h$ is not constant and such that $G^\m_f$ and $G^\m_h$ are both connected, the interiors of 
$\Vor^\m_H(f)$ and $\Vor^\m_H(h)$ are disjoint.
\end{prop}

\begin{proof} The proof is similar to that of
  Proposition~\ref{prop:proj2-2}. Namely, by 
Proposition~\ref{claim1bis}, we have 
$\dl^\m_f\neq \dl^\m_h$. By Proposition~\ref{prop:proj3-1}, 
reasoning by absurd, suppose there are two points 
$\dl^\m_f+ \epsilon$ and $\dl^\m_h + \varepsilon$ in 
the relative interiors of $\square^\m_{f}$ and $\square^\m_h$,
respectively, with $d^*(\dl^\m_f+\epsilon) = d^*(\dl^\m_h +
\varepsilon)$. 
By definition, for each $e\in\E(G)$ we have:
\begin{itemize}
\item If $e\notin \E(G^\m_f)$ then $\epsilon_e =0$, whereas if $e\in
  \E(G^\m_f)$ then $|\epsilon_e|< \frac 12$.
\item If $e\notin \E(G^\m_h)$ then $\varepsilon_e =0$, whereas if
  $e\in \E(G^\m_h)$ then $|\varepsilon_e|< \frac 12$.
\end{itemize}

Let $\eta := f-h$, and let $S$ be the set of all vertices of $G$ where
$\eta$ takes its maximum value. Since $f-h$ is not constant, $S$ is a
proper subset of $V(G)$. Since $G^\m_f$ and $G^\m_h$ are
connected spanning subgraphs of $G$, we have 
$\E(G^\m_f) \cap \E(S, V-S) \neq \emptyset$ and 
$\E(G^\m_h) \cap \E(S, V-S) \neq \emptyset$. 

By assumption, 
\begin{equation}\label{eq4bis}
 d^*(\dl^\m_f - \dl^\m_h) = d^*(\varepsilon - \epsilon).
\end{equation}
We claim that for each edge $e\in \E(S, V-S)$ we have 
$\dl^\m_f(e) - \dl^\m_h (e) \leq \varepsilon_e - \epsilon_e$, with
strict inequality if $e$ belongs either to $G^\m_f$ or $G^\m_h$. 
Once this has been proved, we get a contradiction to \eqref{eq4bis} as in
the proof of Proposition~\ref{prop:proj2-2}.

To prove the claim, we proceed by a case by case analysis:
\begin{itemize}
\item If $e=uv \in \E(S, V-S) \cap \E(G^\m_f)\cap \E(G^\m_h)$ then 
both $f(v)-f(u)+\m_e$ and $h(v) -h(u)+ \mu_e$ 
are divisible by $\ell_e$, and since $\eta(v)< \eta(u)$, we get 
\begin{align*}
\dl^\m_f(e) - \dl^\m_h (e) =& 
\frac 1{\ell_e}\Bigl(f(v)-f(u) +\m_e\Bigr)  
-\frac1{\ell_e}\Bigl(h(v)-h(u)+\m_e\Bigr)\\
=& \frac 1{\ell_e}\Bigl(\eta(v)-\eta(u)\Bigr) 
\leq -1\,.
\end{align*}
On the other hand, $- 1<\varepsilon_e -\epsilon_e < 1$, proving the claim in this case.
\item Similarly, if $e=uv \in \E(S, V-S) \cap \E(G^\m_f)$ but $e\not\in\E(G^\m_h),$, then
$$
\dl^\m_f(e) - \dl^\m_h (e) =
\frac 1{\ell_e}\Bigl(h(v)-h(u) + \m_e+\eta(v)-\eta(u)\Bigr) 
- \Bigl\lfloor \frac{h(v)-h(u)+\m_e}{\ell_e}\Bigr\rfloor - 1/2,
$$
whence $\dl^\m_f(e) - \dl^\m_h (e)\leq \frac 12$, while 
$-\frac 12 < \varepsilon_e -\epsilon_e < \frac 12$.
\item If $e\in \E(S, V-S) \cap \E(G^\m_h)$ but $e\not\in\E(G^\m_f)$, the proof is similar. 
\item Finally, if $e \in \E(S, V-S)$ but $e$ is neither in $G^\m_f$
  nor in $G^\m_h$, then $\dl^\m_f(e) - \dl^\m_h(e) \leq 0 $, 
while  $\varepsilon(e) -\epsilon(e)=0$.
\end{itemize}
This proves the claim, and finishes the proof of the proposition.
\end{proof}
 
Here is the main theorem of this section:

 \begin{thm}\label{thm:projection3} The set of polytopes
   $\Vor^\m_H(f)$ with $G^\m_f$ 
connected provide a tiling $\Vor^\m_{G,\ell}$ 
 of $H_{0,\R}$. Each top-dimensional cell in this tiling is congruent to the Voronoi 
 cell of a connected subgraph of $G$. More precisely, $\Vor^\m_H(f)$
 is congruent to $\Vor_{G^\m_f}(O)$ for each $f\in C^0(G,\Z)$ with $G^\m_f$ 
connected.
 \end{thm}
 
 \begin{remark}\rm This is the most general type of tiling we can get. Note that it might happen for particular choices of $\ell$ and $\m$ that none of the subgraphs $G_f^\m$ is equal to $G$.
 \end{remark}

The proof goes as in the proof of Theorem~\ref{thm:projection2}, so we
need first to extend Lemma~\ref{lem:sharedfacet}. For future use in~\cite{AE2}, we also provide a complete characterization of intersecting Voronoi cells.

\subsection{Characterization of intersecting Voronoi cells and proof of Theorem~\ref{thm:projection3}} We start by the extension of Lemma~\ref{lem:sharedfacet}. Let $f\in C^{0}(G,\Z)$ with connected $G^\m_f$. Let $S\subseteq V$ be
a proper subset such that the characteristic function $\beta_f$ of the
oriented cut $\E_{G^\m_f}(V-S, S)$ is a bond. Consider the functions
$h_n:=f+ n\chi_{\ind{S}}$ for $n\in\N$ and let $h:=h_n$ for the
smallest positive integer $n$ such that $G_h$ is connected. Note that
$G^\m_h[S]=G^\m_f[S]$ and $G^\m_h[V-S]=G^\m_f[V-S]$. Also, as before, 
$\dl^\m_f(e)-1\leq\dl^\m_h(e)\leq\dl_f(e)$ for each $e\in\E(S,V-S)$,
the upper bound being achieved only if $\dl^\m_f(e)\not\in\Z$ and the
lower bound being achieved only if $\dl^\m_f(e)\in\Z$.

The characteristic funtion $\beta_h$ of the oriented cut
$\E_{G^\m_h}(V-S,S)$ is also a bond. Let $\eta_f$ and $\eta_h$ be the
extensions by zero of $\beta_f$ and $\beta_h$ to $C^1(G,\Z)$, and
define
$$
\eta^\m_{S,f} = \frac 12 \eta_f + \frac 12 \eta_h.
$$ 

 \begin{lemma}\label{lem:sharedfacet2}
  Notations as above, we have $\dl^\m_h -\dl^\m_f = \eta_{S,f}^\m$. In
  addition, the Voronoi cells $\Vor^\m_{H}(f)$ and $\Vor^\m_H(h)$
  intersect at the the facet of $\Vor^\m_H(f)$ (resp.~$\Vor^\m_H(h)$) 
given by the bond $\beta_f$ of $G^\m_f$ (resp. by the bond $-\beta_h$ of $G^\m_h$). 
 \end{lemma}

 \begin{proof} As before, the first assertion follows from a case
   analysis. Let $e\in\E(G)$. If $e\in\E(S)$ or $e\in\E(V-S)$, then 
 $\dl^\m_h(e) =\dl^\m_f(e) = \dl^\m_f(e) + \eta^\m_{S,f}(e)$, 
since $\eta^\m_{S,f}(e)=0$. The same holds if $e\in\E(S,V-S)$ 
but $e\not\in\E(G^\m_f)\cup\E(G^\m_h)$.  Suppose now
$e\in\E(S,V-S)$. If $e\in \E(G^\m_f)\cap\E(G^\m_h)$, then 
$\dl^\m_h(e) = \dl^\m_f(e)-1 = \dl^\m_h(e)+ \eta^\m_{S,f}(e)$, as
$\eta_f(e) = \eta_h(e)= -1$.  If $e\in \E(G^\m_f)$ but
$e\not\in\E(G^\m_h)$, then 
$\dl^\m_h(e) = \dl^\m_f(e)-1/2 = \dl^\m_f(e) + \eta^\m_{S,f}(e)$, 
since $\eta_f(e)=-1$ and $\eta_h(e) =0$. The same holds if $e\in \E(G^\m_h)$ but
$e\not\in\E(G^\m_f)$, the only difference been that now $\eta_f(e)=0$
and $\eta_h(e) =-1$. The first assertion is proved.
  
As for the second assertion, by Theorem~\ref{thm:projection1} applied
to $G^\m_f$, 
the facet of $\Vor_{G^\m_f}(O)$ given by the bond $\beta_f$ is the
image under $d^*$ of 
$\frac 12 \eta_f + \square_{G^\m_f[S],0} + \square_{G^\m_f[V-S],0}$, where 
  $\square_{G^\m_f[S],0}$ and $\square_{G^\m_f[V-S],0}$ are viewed as
  subsets of $\square_0$. 
Similarly, the facet of $\Vor_{G^\m_h}(O)$ 
  given by the bond $-\beta_h$ is the
image under $d^*$ of  
  $-\frac 12 \eta_h + \square_{G^\m_h[S],0} +
  \square_{G^\m_h[V-S],0}$. 
Since $G^\m_h[S]=G^\m_f[S]$ and $G^\m_h[V-S]=G^\m_f[V-S]$, it follows
from the first assertion that 
\begin{equation}\label{facetfhm}
\dl^\m_f + \frac 12 \eta_f + \square_{G^\m_f[S],0} + \square_{G^\m_f[V-S],0} =  
\dl^\m_h - \frac 12 \eta_h + \square_{G^\m_h[S],0} +\square_{G^\m_h[V-S],0}.
\end{equation}
Since  $\Vor^\m_H(f) = d^*(\dl^\m_f)+ \Vor_{G^\m_f}(O)$, 
  it follows that the facet of $\Vor^\m_H(f)$ defined by $\beta_f$ is
  the image under $d^*$ of the left hand side of \eqref{facetfhm}. 
Analogously, the facet of $\Vor^\m_H(h)$ given by the bond
$-\beta_h$ is the image under $d^*$ of the right hand side of
\eqref{facetfhm}. So they coincide. The lemma follows since the two Voronoi cells 
$\Vor^\m_H(f)$ and $\Vor^\m_H(h)$ have disjoint interiors by 
Proposition~\ref{prop:proj3-2}.
 \end{proof}
The following proposition provides a characterization of intersecting
Voronoi cells, and will be used in~\cite{AE2}.

\begin{prop} Let $f_1,f_2\in C^0(G,\Z)$. Let $X_1,\dots,X_q\subseteq V$ be the level subsets
of $f_2-f_1$, in increasing order. Let $D_1$ be the coherent acyclic orientation of
the cut subgraph of $G^\m_{f_1}$ induced by the ordered partition
$X_1,\dots,X_q$ of $V$, and $D_2$ that of the cut subgraph of
$G^\m_{f_2}$ 
induced by the same partition in the reverse order,
$X_q,\dots,X_1$. Then $\Vor^\m_H(f_1)$ intersects
$\Vor^\m_H(f_2)$ if and only if 
\begin{equation}\label{f1f1f2f2}
\dl^\m_{f_1}+\frac 12 \chi_{\indbi{D}{1}} = \dl^\m_{f_2}+\frac 12 \chi_{\indbi D2},
\end{equation}
where $\chi_{\indbi Di}$ is the characteristic function of $D_i$, taking
value $+1$ at $e\in\E(D_i)$, value $-1$ at $e\in\E$ with $\bar
e\in\E(D_i)$, and value $0$ elsewhere, for $i=1,2$. Furthermore, in
this case, letting $\mathfrak f_1$ and $\mathfrak f_2$ be the
corresponding faces of $\Vor^\m_H(f_1)$ and $\Vor^\m_H(f_2)$ to $D_1$
and $D_2$, respectively, we have
$$
\mathfrak f_1=\mathfrak f_2=\Vor^\m_H(f_1)\cap \Vor^\m_H(f_2).
$$
\end{prop}

\begin{proof} If $\Vor^\m_H(f_1)$ and $\Vor^\m_H(f_2)$
  intersect then there are $x,y\in C^1(G,\R)$ with
$\dl^\m_{f_1}+x\in\square^\m_{f_1}$ and
$\dl^\m_{f_2}+y\in\square^\m_{f_2}$ such that $d^*(\alpha)=0$, where
$$
\alpha:=\dl^\m_{f_2}-\dl^\m_{f_1}+y-x.
$$

\smallskip

{\bf Claim 1:} Equation \eqref{f1f1f2f2} holds and 
\begin{equation}\label{xy}
x\in \frac 12\chi_{\indbi D1}+\sum_{i=1}^q\square_{G^\m_{f_1}[X_i],0}\quad\text{and}\quad
y\in \frac 12\chi_{\indbi D2}+\sum_{i=1}^q\square_{G^\m_{f_2}[X_i],0}.
\end{equation}
Equivalently, for each $e\in\E(G)$, we have that 
\begin{equation}\label{c1c2}
\dl^\m_{f_1}(e)+\frac 12 \chi_{\indbi D1}(e) = \dl^\m_{f_2}(e)+\frac 12
\chi_{\indbi D2}(e)
\end{equation}
and 
\begin{equation}\label{xye}
\Big|x_e-\frac 12\chi_{\indbi D1}(e)\Big|\leq \frac 12 b_{1,e}\quad\text{and}\quad
\Big|y_e-\frac 12\chi_{\indbi D2}(e)\Big|\leq\frac 12 b_{2,e},
\end{equation}
where $b_{j,e}=1$ if $e\in\E(G^\m_{f_j}[X_i])$ for some $i$ and
$b_{j,e}=0$ otherwise, for $j=1,2$.

Indeed, for each $i=1,\dots,q$, since $f_1|_{X_i}$ and $f_2|_{X_i}$
differ by a constant, for each $e\in\E(G[X_i])$ we have
$\dl^\m_{f_1}(e)=\dl^\m_{f_2}(e)$ and
$\chi_{\indbi D1}(e)=\chi_{\indbi D2}(e)=0$. Thus \eqref{c1c2} and \eqref{xye} hold.

Now, for each
$j=2,\dots,q$, let $S_j:=X_1\cup\cdots\cup X_{j-1}$
and $T_j:=X_j\cup\cdots\cup X_q$. 

\smallskip

{\bf Claim 2:} We have that $\alpha_e\geq 0$ for
each $e\in\E(S_j,T_j)$, with equality only if \eqref{c1c2} and \eqref{xye} hold.

Indeed, given
$e=uv\in\E(S_j,T_j)$ we have that $h(v)>h(u)$, and thus
\begin{equation}\label{f2f1>}
\frac{f_2(v)-f_2(u)+\m_e}{\ell_e}>
\frac{f_1(v)-f_1(u)+\m_e}{\ell_e}.
\end{equation}
There are four cases to consider, whose analysis finishes the
proof of Claim 2:
\begin{enumerate}
\item If $e\not\in\E(G^\m_{f_1})\cup\E(G^\m_{f_2})$ then
  $\chi_{\indbi D1}(e)=\chi_{\indbi D2}(e)=0$ and $x_e=y_e=0$, whence \eqref{xye}
  holds. Furthermore, $\alpha_e=\dl^\m_{f_2}(e)-\dl^\m_{f_1}(e)$, whence $\alpha_e\geq 0$
from \eqref{f2f1>}, with equality only if \eqref{c1c2} holds.
\item If $e\in\E(G^\m_{f_1})$ but
$e\not\in\E(G^\m_{f_2})$ then $\chi_{\indbi D1}(e)=1$ and $x_e\leq 1/2$,
whereas $\chi_{\indbi D2}(e)=0$ and $y_e=0$. Thus 
$\alpha_e\geq\dl^\m_{f_2}(e)-\dl^\m_{f_1}(e)-1/2$. But
$\dl^\m_{f_2}(e)\geq\dl^\m_{f_1}(e)+1/2$ 
from \eqref{f2f1>}, whence $\alpha_e\geq 0$, with equality only if
$x_e=1/2$ and $\dl^\m_{f_2}(e)=\dl^\m_{f_1}(e)+1/2$, that is, only if \eqref{c1c2}
and \eqref{xye} hold.
\item If $e\in\E(G^\m_{f_2})$ but
$e\not\in\E(G^\m_{f_1})$ then $\chi_{\indbi D1}(e)=0$ and $x_e=0$, whereas 
$\chi_{\indbi D2}(e)=-1$ and $y_e\geq -1/2$. Thus 
$\alpha_e\geq\dl^\m_{f_2}(e)-\dl^\m_{f_1}(e)-1/2$. But 
$\dl^\m_{f_2}(e)\geq\dl^\m_{f_1}(e)+1/2$ from \eqref{f2f1>}, whence 
$\alpha_e\geq 0$, with equality only if $y_e=-1/2$ and
$\dl^\m_{f_2}(e)-1/2=\dl^\m_{f_1}(e)$, that is, only if \eqref{c1c2}
and \eqref{xye} hold.
\item Finally, if
$e\in\E(G^\m_{f_1})\cap\E(G^\m_{f_2})$ then $\chi_{\indbi D1}(e)=1$ and
$x_e\leq 1/2$, whereas 
$\chi_{\indbi D2}(e)=-1$ and $y_e\geq -1/2$. Thus 
$\alpha_e\geq\dl^\m_{f_2}(e)-\dl^\m_{f_1}(e)-1$. 
But $\dl^\m_{f_2}(e)\geq\dl^\m_{f_1}(e)+1$ from \eqref{f2f1>}, whence
$\alpha_e\geq 0$, with equality only if $x_e=1/2$, $y_e=-1/2$ and 
$\dl^\m_{f_2}(e)-1/2=\dl^\m_{f_1}(e)+1/2$, that is, only if \eqref{c1c2}
and \eqref{xye} hold.
\end{enumerate}

To finish the proof of Claim 1, observe now that, since
$d^*(\alpha)=0$, we have that
$$
\sum_{e\in\E(S_j,T_j)}\alpha_e=0.
$$
It thus follows from Claim 2 that $\alpha_e= 0$ and hence \eqref{c1c2}
and \eqref{xye} hold for each $e\in\E(S_j,T_j)$. As this holds for
each $j=2,\dots,q$, the proof of Claim 1 is finished.

Assume now that \eqref{f1f1f2f2} holds. Now, by
Theorem~\ref{thm:projection1bis}, we have
\begin{equation}\label{fff}
\mathfrak f_1=d^*\Big(\dl^\m_{f_1}+\frac 12\chi_{\indbi D1}
+\sum_{i=1}^q\square_{G^\m_{f_1}[X_i],0}\Big)
\quad\text{and}\quad
\mathfrak f_2=d^*\Big(\dl^\m_{f_2}+\frac 12\chi_{\indbi D2}
+\sum_{i=1}^q\square_{G^\m_{f_2}[X_i],0}\Big).
\end{equation}
Now, if \eqref{f1f1f2f2} holds then $G^\m_{f_1}[X_i]=G^\m_{f_2}[X_i]$
for each $i$. Indeed, if $e\in\E(G[X_i])$ for some $i$ then
$\chi_{\indbi{D}1}(e)=\chi_{\indbi{D}2}(e)=0$,  and hence \eqref{f1f1f2f2} yields that 
$e\in\E(G^\m_{f_1})$ if and only if $e\in\E(G^\m_{f_2})$. Thus, if
\eqref{f1f1f2f2} holds then it follows from \eqref{fff} that 
$\mathfrak f_1=\mathfrak f_2$ and
hence $\Vor^\m_H(f_1)$ intersects $\Vor^\m_H(f_2)$. 

Furthermore, a point on $\Vor^\m_H(f_1)\cap\Vor^\m_H(f_2)$ can be
expressed as $d^*(\dl^\m_{f_1}+x)$ for
$\dl^\m_{f_1}+x\in\square^\m_{f_1}$ and as $d^*(\dl^\m_{f_2}+y)$ for 
$\dl^\m_{f_2}+y\in\square^\m_{f_2}$; as it is the same point, $d^*(\alpha)=0$, where
$$
\alpha:=\dl^\m_{f_2}-\dl^\m_{f_1}+y-x.
$$
It follows then from Claim 1 and from \eqref{fff} that that point is
on $\mathfrak f_1$ and on $\mathfrak f_2$. To conclude, 
$$
\Vor^\m_H(f_1)\cap\Vor^\m_H(f_2)=\mathfrak f_1=\mathfrak f_2.
$$
The proof is finished.
\end{proof}

 \begin{proof}[Proof of Theorem~\ref{thm:projection3}]
   Proceed as in the proof of Theorem~\ref{thm:projection2} by applying 
  Lemma \ref{lem:sharedfacet2}.
  \end{proof}

\vspace{.5cm}

\subsection*{Acknowledgements.}  This project benefited very much from
the hospitality of the \'Ecole Normale Superieure in Paris and the
Instituto Nacional de Matem\'atica Pura e Aplicada in Rio de Janeiro
during the mutual visits of both authors. We thank both 
institutions and their members for making them possible.
We are specially grateful to the Brazilian-French Network in
Mathematics for providing financial support for a visit of E.E. at ENS Paris
and a visit of O.A. at IMPA.

\bibliographystyle{alpha}
\bibliography{bibliography}

\begin{thebibliography}{BTW87}

\bibitem[AE20a]{AE2}
Omid Amini and Eduardo Esteves.
\newblock Voronoi tilings, toric arrangements and degenerations of line bundles
  {II}.
\newblock {\em preprint}, 2020.

\bibitem[AE20b]{AE3}
Omid Amini and Eduardo Esteves.
\newblock Voronoi tilings, toric arrangements and degenerations of line bundles
  {III}.
\newblock {\em preprint}, 2020.

\bibitem[AH99]{AH99}
Klaus Altmann and Lutz Hille.
\newblock Strong exceptional sequences provided by quivers.
\newblock {\em Algebras and Representation Theory}, 2(1):1--17, 1999.

\bibitem[AM10]{AM}
Omid Amini and Madhusudan Manjunath.
\newblock {Riemann-Roch for Sub-Lattices of the Root Lattice $A_n$}.
\newblock {\em Electronic Journal of Combinatorics}, 17:Research Paper 124,
  50p, 2010.

\bibitem[Ami]{Ami}
Omid Amini.
\newblock Lattice of integer flows and poset of strongly connected
  orientations.
\newblock {\em preprint}.

\bibitem[Ami14]{Ami-W}
Omid Amini.
\newblock Equidistribution of {W}eierstrass points on curves over
  non-{A}rchimedean fields.
\newblock {\em preprint}, (revised version 2020), 2014.

\bibitem[Ami18]{Ami-hdr}
Omid Amini.
\newblock {\em Geometry of graphs and applications}.
\newblock M\'emoire d'HdR, Sorbonne Universit\'e, 2018.

\bibitem[AN20]{AN20}
Omid Amini and Noema Nicolussi.
\newblock Moduli of hybrid curves and variations of canonical measures.
\newblock {\em arXiv preprint arXiv:2007.07130}, 2020.

\bibitem[Art06]{Art}
Igor~Vadimovich Artamkin.
\newblock Discrete {T}orelli theorem.
\newblock {\em Sbornik: Mathematics}, 197(8):1109--1120, 2006.

\bibitem[BD09]{BD}
Felix Breuer and Aaron Dall.
\newblock Viewing counting polynomials as {H}ilbert functions via {E}hrhart
  theory.
\newblock {\em arXiv preprint arXiv:0911.5109}, 2009.

\bibitem[Big99]{Biggs}
Norman Biggs.
\newblock Chip-firing and the critical group of a graph.
\newblock {\em Journal of Algebraic Combinatorics}, 9(1):25--45, 1999.

\bibitem[BJ16]{BJ}
Matthew Baker and David Jensen.
\newblock Degeneration of linear series from the tropical point of view and
  applications.
\newblock In {\em Nonarchimedean and Tropical Geometry}, pages 365--433.
  Springer, 2016.

\bibitem[BLN97]{BHN97}
Roland Bacher, Pierre~de {La Harpe}, and Tatiana Nagnibeda.
\newblock The lattice of integral flows and the lattice of integral cuts on a
  finite graph.
\newblock {\em Bulletin de la soci{\'e}t{\'e} math{\'e}matique de France},
  125(2):167--198, 1997.

\bibitem[BLS91]{BLS}
Anders Bj{\"o}rner, L{\'a}szl{\'o} Lov{\'a}sz, and Peter Shor.
\newblock Chip-firing games on graphs.
\newblock {\em European Journal of Combinatorics}, 12(4):283--291, 1991.

\bibitem[BM08]{BM}
John~Adrian Bondy and U.~S.~Ramachandra Murty.
\newblock {\em Graph Theory}, volume 244.
\newblock Springer, 2008.

\bibitem[BN07]{BN06}
Matthew Baker and Serguei Norine.
\newblock {R}iemann--{R}och and {A}bel--{J}acobi theory on a finite graph.
\newblock {\em Advances in Mathematics}, 215(2):766--788, 2007.

\bibitem[Bol98]{Bol}
B{\'e}la Bollob{\'a}s.
\newblock {\em Modern graph theory}, volume 184.
\newblock Springer, 1998.

\bibitem[BS12]{BS}
Felix Breuer and Raman Sanyal.
\newblock Ehrhart theory, modular flow reciprocity, and the tutte polynomial.
\newblock {\em Mathematische Zeitschrift}, 270(1-2):1--18, 2012.

\bibitem[BTW87]{BTW}
Per Bak, Chao Tang, and Kurt Wiesenfeld.
\newblock Self-organized criticality: An explanation of the 1/f noise.
\newblock {\em Physical review letters}, 59(4):381, 1987.

\bibitem[BZ06]{BZ}
Matthias Beck and Thomas Zaslavsky.
\newblock The number of nowhere-zero flows on graphs and signed graphs.
\newblock {\em Journal of Combinatorial Theory, Series B}, 96(6):901--918,
  2006.

\bibitem[CC19]{CC}
Lucia Caporaso and Karl Christ.
\newblock Combinatorics of compactified universal {J}acobians.
\newblock {\em Advances in Mathematics}, 346:1091--1136, 2019.

\bibitem[Che10]{Chen}
Beifang Chen.
\newblock Orientations, lattice polytopes, and group arrangements {I}:
  {C}hromatic and tension polynomials of graphs.
\newblock {\em Annals of combinatorics}, 13(4):425--452, 2010.

\bibitem[Chv75]{Chv}
Va{\v{s}}ek Chv{\'a}tal.
\newblock On certain polytopes associated with graphs.
\newblock {\em Journal of Combinatorial Theory, Series B}, 18(2):138--154,
  1975.

\bibitem[CKV13]{CKV}
Sebastian {Casalaina-Martin}, Jesse Kass, and Filippo Viviani.
\newblock The geometry and combinatorics of cographic toric face rings.
\newblock {\em Algebra \& Number Theory}, 7(8):1781--1815, 2013.

\bibitem[CKV14]{CKV2}
Sebastian {Casalaina-Martin}, Jesse~Leo Kass, and Filippo Viviani.
\newblock The singularities and birational geometry of the universal
  compactified jacobian.
\newblock {\em Algebraic Geometry}, 4(3):353--393, 2014.

\bibitem[CS82]{CS2}
John Conway and Neil Sloane.
\newblock Voronoi regions of lattices, second moments of polytopes, and
  quantization.
\newblock {\em IEEE Transactions on Information Theory}, 28(2):211--226, 1982.

\bibitem[CS84]{CS}
John Conway and Neil Sloane.
\newblock On the {V}oronoi regions of certain lattices.
\newblock {\em SIAM Journal on Algebraic Discrete Methods}, 5(3):294--305,
  1984.

\bibitem[CV10]{CV}
Lucia Caporaso and Filippo Viviani.
\newblock Torelli theorem for graphs and tropical curves.
\newblock {\em Duke Mathematical Journal}, 153(1):129--171, 2010.

\bibitem[DG16]{DG}
Zsuzsanna Dancso and Stavros Garoufalidis.
\newblock A construction of the graphic matroid from the lattice of integer
  flows.
\newblock {\em arXiv preprint arXiv:1611.06282}, 2016.

\bibitem[Dha90]{Dhar}
Deepak Dhar.
\newblock Self-organized critical state of sandpile automaton models.
\newblock {\em Physical Review Letters}, 64(14):1613, 1990.

\bibitem[DM69]{DM69}
Pierre Deligne and David Mumford.
\newblock The irreducibility of the space of curves of given genus.
\newblock {\em Publications Math{\'e}matiques de l'IHES}, 36:75--109, 1969.

\bibitem[EM02]{EM}
Eduardo Esteves and Nivaldo Medeiros.
\newblock Limit canonical systems on curves with two components.
\newblock {\em Inventiones mathematicae}, 149(2):267--338, 2002.

\bibitem[ES07]{ES07}
Eduardo Esteves and Parham Salehyan.
\newblock Limit {W}eierstrass points on nodal reducible curves.
\newblock {\em Transactions of the American Mathematical Society},
  359(10):5035--5056, 2007.

\bibitem[Est01]{Esteves01}
Eduardo Esteves.
\newblock Compactifying the relative {J}acobian over families of reduced
  curves.
\newblock {\em Transactions of the American Mathematical Society},
  353(8):3045--3095, 2001.

\bibitem[Gab93a]{Gab2}
Andrei Gabrielov.
\newblock Abelian avalanches and tutte polynomials.
\newblock {\em Physica A: Statistical Mechanics and its Applications},
  195(1-2):253--274, 1993.

\bibitem[Gab93b]{Gab}
Andrei Gabrielov.
\newblock Avalanches, sandpiles and {T}utte decomposition.
\newblock In {\em The Gelfand Mathematical Seminars, 1990--1992}, pages 19--26.
  Springer, 1993.

\bibitem[Ger82]{Ger}
Lothar Gerritzen.
\newblock Die {J}acobi-{A}bbildung {\"u}ber dem {R}aum der {M}umfordkurven.
\newblock {\em Mathematische Annalen}, 261(1):81--100, 1982.

\bibitem[LP09]{LP}
L{\'a}szl{\'o} Lov{\'a}sz and Michael~D Plummer.
\newblock {\em Matching theory}, volume 367.
\newblock American Mathematical Soc., 2009.

\bibitem[NT74]{NT}
George Nemhauser and Leslie~Earl Trotter.
\newblock Properties of vertex packing and independence system polyhedra.
\newblock {\em Mathematical programming}, 6(1):48--61, 1974.

\bibitem[OS79]{OS79}
Tadao Oda and Conjeerveram Seshadri.
\newblock Compactifications of the generalized {J}acobian variety.
\newblock {\em Transactions of the American Mathematical Society}, pages 1--90,
  1979.

\bibitem[Pad73]{pad1}
Manfred~W Padberg.
\newblock On the facial structure of set packing polyhedra.
\newblock {\em Mathematical programming}, 5(1):199--215, 1973.

\bibitem[Pad74]{pad2}
Manfred~W Padberg.
\newblock Perfect zero--one matrices.
\newblock {\em Mathematical Programming}, 6(1):180--196, 1974.

\bibitem[Sch03]{Sch}
Alexander Schrijver.
\newblock {\em Combinatorial optimization: polyhedra and efficiency},
  volume~24.
\newblock Springer Science \& Business Media, 2003.

\bibitem[Whi32]{Whi1}
Hassler Whitney.
\newblock {Congruent Graphs and the Connectivity of Graphs}.
\newblock {\em American Journal of Mathematics}, 54(1):150--168, 1932.

\bibitem[Whi33]{Whi2}
Hassler Whitney.
\newblock {2-Isomorphic Graphs}.
\newblock {\em American Journal of Mathematics}, 55(1):245--254, 1933.

\end{thebibliography}
\end{document}